\documentclass[11pt,english,reqno]{extarticle}
\usepackage{amsmath}
\usepackage{amsthm}
\usepackage{libertineRoman}
\usepackage[scaled=0.92]{helvet}
\usepackage{courier}
\usepackage[libertine]{newtxmath}
\usepackage[T1]{fontenc}
\usepackage[utf8]{inputenc}
\usepackage{color}
\usepackage{babel}
\usepackage{cprotect}
\usepackage{mathrsfs}
\usepackage{mathtools}
\usepackage{url}
\usepackage{enumitem}
\usepackage{stmaryrd}
\usepackage{geometry}
\usepackage{wasysym}
\usepackage{microtype}
\usepackage[pdfusetitle,
 bookmarks=true,bookmarksnumbered=false,bookmarksopen=false,
 breaklinks=true,pdfborder={0 0 0},pdfborderstyle={},backref=false,colorlinks=true]
 {hyperref}
\hypersetup{
 allcolors=blue}

\makeatletter
\numberwithin{equation}{section}
\numberwithin{figure}{section}
\numberwithin{table}{section}
\usepackage{cleveref}
\theoremstyle{plain}
\newtheorem{thm}{\protect\theoremname}[section]
\theoremstyle{definition}
\newtheorem{defn}[thm]{\protect\definitionname}
\theoremstyle{plain}
\newtheorem{prop}[thm]{\protect\propositionname}
\theoremstyle{plain}
\newtheorem{cor}[thm]{\protect\corollaryname}
\theoremstyle{plain}
\newtheorem{lem}[thm]{\protect\lemmaname}
\theoremstyle{remark}
\newtheorem{rem}[thm]{\protect\remarkname}
\newlist{casenv}{enumerate}{4}
\setlist[casenv]{leftmargin=*,align=left,widest={iiii}}
\setlist[casenv,1]{label={{\itshape\ \casename} \arabic*.},ref=\arabic*}
\setlist[casenv,2]{label={{\itshape\ \casename} \roman*.},ref=\roman*}
\setlist[casenv,3]{label={{\itshape\ \casename\ \alph*.}},ref=\alph*}
\setlist[casenv,4]{label={{\itshape\ \casename} \arabic*.},ref=\arabic*}
\newlist{thmstepnv}{enumerate}{4}
\setlist[thmstepnv]{leftmargin=*,align=left,wide,labelwidth=0pt,labelindent=0pt}
\setlist[thmstepnv,1]{label={\itshape {\thmstepname} \arabic*.},ref=\arabic*}
\setlist[thmstepnv,2]{label={\itshape {\thmstepname} {\thethmstepnvi\alph*}.},ref=\thethmstepnvi\alph*}
\setlist[thmstepnv,3]{label={\itshape {\thmstepname\ \alph*.}},ref=\alph*}
\setlist[thmstepnv,4]{label={\itshape {\thmstepname} \arabic*.},ref=\arabic*}

\usepackage{tikz}
\hypersetup{final}
\let\prettyref\cref
\crefname{figure}{Figure}{Figures}
\crefname{section}{Section}{Sections}

\crefdefaultlabelformat{#2\textup{#1}#3}
\crefformat{equation}{\textup{(#2#1#3)}}
\crefformat{enu}{#1}
\crefmultiformat{equation}{(#2#1#3)}%
{ and~(#2#1#3)}{, (#2#1#3)}{ and~(#2#1#3)}
\crefrangeformat{equation}{(#3#1#4--#5#2#6)}
\usepackage{mleftright}
\mleftright

\usepackage{caption}
\setlist[casenv]{leftmargin=*,align=left,wide,labelwidth=0pt,labelindent=0pt}
\setlist[casenv,1]{label={{\itshape \casename} \arabic*.},ref=\arabic*}
\setlist[casenv,2]{label={{\itshape \casename} \roman*.},ref=\roman*}
\setlist[casenv,3]{label={{\itshape \casename\ \alph*.}},ref=\alph*}
\setlist[casenv,4]{label={{\itshape \casename} \arabic*.},ref=\arabic*}

\usepackage{aligned-overset}
\let\oversetx\overset
\usepackage{placeins}

\usepackage{csquotes}
\AtBeginDocument{\AtEveryBibitem{\clearfield{eprintclass}\clearfield{primaryclass}}}

\ifdefined\showcaptionsetup
 \PassOptionsToPackage{caption=false}{subfig}
\fi
\usepackage{subfig}
\makeatother

\usepackage[style=trad-plain,isbn=false,doi=false,url=false]{biblatex}
\providecommand{\casename}{Case}
\providecommand{\corollaryname}{Corollary}
\providecommand{\definitionname}{Definition}
\providecommand{\lemmaname}{Lemma}
\providecommand{\propositionname}{Proposition}
\providecommand{\remarkname}{Remark}
\providecommand{\theoremname}{Theorem}
\providecommand{\thmstepname}{Step}

\begin{document}
\global\long\def\supp{\operatorname{supp}}%

\global\long\def\Uniform{\operatorname{Uniform}}%

\global\long\def\dif{\mathrm{d}}%

\global\long\def\e{\mathrm{e}}%

\global\long\def\ii{\mathrm{i}}%

\global\long\def\Cov{\operatorname{Cov}}%

\global\long\def\dist{\operatorname{dist}}%

\global\long\def\Var{\operatorname{Var}}%

\global\long\def\pv{\operatorname{p.v.}}%

\global\long\def\e{\mathrm{e}}%

\global\long\def\p{\mathrm{p}}%

\global\long\def\Law{\operatorname{Law}}%

\global\long\def\supp{\operatorname{supp}}%

\global\long\def\image{\operatorname{image}}%

\global\long\def\dif{\mathrm{d}}%

\global\long\def\eps{\varepsilon}%

\global\long\def\sgn{\operatorname{sgn}}%

\global\long\def\tr{\operatorname{tr}}%

\global\long\def\Hess{\operatorname{Hess}}%

\global\long\def\Re{\operatorname{Re}}%

\global\long\def\Im{\operatorname{Im}}%

\global\long\def\Dif{\operatorname{D}}%

\global\long\def\divg{\operatorname{div}}%

\global\long\def\arsinh{\operatorname{arsinh}}%

\global\long\def\sech{\operatorname{sech}}%

\global\long\def\erf{\operatorname{erf}}%

\global\long\def\Cauchy{\operatorname{Cauchy}}%

\global\long\def\artanh{\operatorname{artanh}}%

\global\long\def\Left{\mathrm{L}}%

\global\long\def\Right{\mathrm{R}}%

\global\long\def\LR{\mathrm{LR}}%

\global\long\def\ssS{\mathsf{S}}%

\global\long\def\GS{\mathsf{GS}}%

\global\long\def\clA{\mathcal{A}}%

\global\long\def\ssN{\mathsf{N}}%

\global\long\def\ssM{\mathsf{M}}%

\global\long\def\ssH{\mathsf{H}}%

\global\long\def\ssQ{\mathsf{Q}}%

\global\long\def\sfs{\mathsf{s}}%

\global\long\def\sfr{\mathsf{r}}%

\global\long\def\sff{\mathsf{f}}%

\global\long\def\sfY{\mathfrak{Y}}%

\global\long\def\TT{\mathbb{T}}%

\global\long\def\RR{\mathbb{R}}%

\global\long\def\ZZ{\mathbb{Z}}%

\global\long\def\PP{\mathbb{P}}%

\global\long\def\NN{\mathbb{N}}%

\global\long\def\scM{\mathscr{M}}%

\global\long\def\st{\mathrel{:}}%

\title{Simultaneous global inviscid Burgers flows\\
with periodic Poisson forcing}
\author{Alexander Dunlap\thanks{Department of Mathematics, Duke University, Durham, NC 27708, USA. Email: \href{mailto:dunlap@math.duke.edu}{\url{dunlap@math.duke.edu}}.}}
\maketitle
\begin{abstract}
We study the inviscid Burgers equation on the circle $\mathbb{T}\coloneqq\mathbb{R}/\mathbb{Z}$
forced by the spatial derivative of a Poisson point process on $\mathbb{R}\times\mathbb{T}$.
We construct global solutions with mean $\theta$ simultaneously for
all $\theta\in\RR$, and in addition construct their associated global
shocks (which are unique except on a countable set of $\theta$).
We then show that as $\theta$ changes, the solution only changes
through the movement of the global shock, and give precise formulas
for this movement. This can be seen as an analogue of previous results
by the author and Yu Gu in the viscous case with white-in-time forcing,
which related the derivative of the solution in $\theta$ to the density
of a particle diffusing in the Burgers flow.
\end{abstract}
\tableofcontents{}

\renewcommand{\listfigurename}{List of figures}

\listoffigures

\vspace{\baselineskip}
An animation is also included in the supplementary material\cprotect\footnote{available online at \url{https://arxiv.org/src/2406.06896v1/anc/moving_global_shock.mp4}}
and described in \prettyref{subsec:illustrations}.

\section{Introduction}

\subsection{Background}

Let $\TT\coloneqq\RR/\ZZ$ and let $\pi\colon\RR\to\TT$ be the projection
map. Let $\mu$ be a purely atomic measure on $\mathbb{R}\times\mathbb{T}$
such that $\ssN\coloneqq\supp\mu$ is a discrete set. We will think
of $\mu$ as a random measure, and our main example will be when $\mu$
is a realization of a homogeneous Poisson point process. We are interested
in global solutions to the forced inviscid Burgers equation formally
given by
\begin{subequations}
\label{eq:uproblem}
\begin{align}
\partial_{t}u_{\theta}(t,x)+\partial_{x}\left(\frac{1}{2}u_{\theta}^{2}+\mu\right)(t,x) & =0,\qquad\theta,t\in\RR,x\in\TT;\label{eq:ueqn}\\
\int_{\TT}u_{\theta}(t,x)\,\dif x & =\theta,\qquad\theta,t\in\RR.\label{eq:uintegral}
\end{align}
\end{subequations}
This type of Burgers equation with discrete forcing has been considered
in the whole-line setting, and briefly in the periodic setting, in
\cite{Bak13,BCK14}. As is usual for the inviscid Burgers equation,
we make sense of the problem \prettyref{eq:uproblem} via entropy
solutions, which we define via Lagrangian minimizers. For $X\in H^{1}([s,t];\TT)$,
and $\theta\in\RR$, we define the action
\begin{equation}
\mathcal{A}_{\theta,s,t}[X]\coloneqq\frac{1}{2}\int_{s}^{t}(X'(r)-\theta)^{2}\,\dif r-\mu(\{(r,X(r))\st r\in[s,t)\}),\label{eq:Adef}
\end{equation}
where $X'$ denotes the derivative of $X$. Roughly speaking, entropy
solutions to \prettyref{eq:uproblem} with initial data $u_{\theta}(s,x)=\theta+\partial_{x}G(x)$
are given by $u_{\theta}(t,x)=X'(t)$, where $X$ minimizes $\mathcal{A}_{\theta,s,t}[X]+G(X(s))$
over all paths satisfying $X(t)=x$. (See \prettyref{def:minimizers}
below.) At some points, the slope of the minimizer is not unique,
and at these ``shock'' points, the solution $u$ is discontinuous
in space. This minimization problem resembles that defining the Hammersley
process \cite{Ham72,AD95}, but here we impose a quadratic penalty
rather than a hard cutoff on the slope of the minimizers.

The integral on the left side of \prettyref{eq:uintegral} is preserved
by the Burgers evolution, as can be seen formally by integrating \prettyref{eq:ueqn}
in space. In terms of the minimization problem \prettyref{eq:Adef},
the role of $\theta$ is to encourage potential minimizers to have
average slope $\theta$. Previous works in the mathematics literature
(see e.g. \cites{EKMS00}{IK03}{GIKP05}{Bak13}) on the forced Burgers
equation in the periodic setting have mostly considered the case $\theta=0$.
This is essentially equivalent to considering the problem for any
fixed value of $\theta$, since if $u$ solves \prettyref{eq:ueqn},
then $\tilde{u}(t,x)=u(t,x-\theta t)+\theta$ solves
\begin{equation}
\partial_{t}\tilde{u}(t,x)+\frac{1}{2}\partial_{x}\tilde{u}^{2}(t,x)+\partial_{x}\mu(t,x-\theta t)=0.\label{eq:veqn}
\end{equation}
If the forcing $\mu$ is taken to be random and shear-invariant in
law, then the solution maps of the problems \prettyref{eq:veqn} and
\prettyref{eq:ueqn} have the same probability distributions.

However, equivalence in law of the solutions with different values
of $\theta$ does not tell us about the behavior of the equation solved
for all values of $\theta$ \emph{simultaneously} (i.e.~with the
same realization of the forcing $\mu$). This has been a topic of
significant recent interest for various models in the KPZ universality
class on the whole line; see e.g.~\cite{FS20,BFS23,BSS23,BSS24,GRASS23}.
In the Burgers setting on the torus, a time-periodic version of the
problem has been discussed at a physics level in \cite{BK07}. Important
questions in the multi-$\theta$ context include (1) how to construct
jointly invariant measures for the randomly forced Burgers equation
simultaneously for all values of $\theta$; and (2) how solutions
sampled from these invariant measures change as the parameter $\theta$
changes. On the whole line, it has been found in all examples that
have been considered that the spatial integral of the process (i.e.~the
process of spatial increments of the KPZ equation) almost surely exhibits
a countable dense set of discontinuities in $\theta$. However, the
proofs of this property in each case are relying on some exact computation
that can be performed for the model, and the physical phenomenology
behind these discontinuities remains to be fully understood.

The author and Yu Gu have studied a related problem for the \emph{viscous}
Burgers equation on the torus, forced by a white-in-time Gaussian
noise, in \cite{DG23}. The equation we considered in that setting
is
\begin{subequations}
\label{eq:vproblem}
\begin{align}
\dif v_{\theta}(t,x) & =\frac{1}{2}[\partial_{x}^{2}v_{\theta}+\partial_{x}(v_{\theta}^{2})](t,x)\dif t+\dif V(t,x), & \theta,t\in\RR,x\in\TT;\label{eq:veqn-1}\\
\int_{\TT}v_{\theta}(t,x)\,\dif x & =\theta, & \theta,t\in\RR,\label{eq:vintegral}
\end{align}
\end{subequations}
where $\dif V$ is the noise. We showed that, at stationarity, we
have
\begin{equation}
u_{\theta_{2}}(t,x)-u_{\theta_{1}}(t,x)=\int_{\theta_{1}}^{\theta_{2}}g_{\theta}(t,x)\,\dif\theta,\label{eq:viscouscase}
\end{equation}
where $g_{\theta}$ is a statistically stationary solution to the
associated coupled PDE
\begin{subequations}
\label{eq:gproblem}
\begin{align}
\partial_{t}g_{\theta}(t,x) & =\frac{1}{2}\Delta g_{\theta}(t,x)+\partial_{x}(u_{\theta}g_{\theta})(t,x), & \theta,t\in\RR,x\in\TT;\label{eq:geqn}\\
\int_{\TT}g_{\theta}(t,x)\,\dif x & =1, & \theta,t\in\RR.\label{eq:gintegral}
\end{align}
\end{subequations}
The problem \prettyref{eq:gproblem} can be obtained by differentiating
\prettyref{eq:vproblem} in $\theta$ and setting $g_{\theta}=\partial_{\theta}u_{\theta}$.
On the other hand, from the form of \prettyref{eq:geqn}, we see that
$g_{\theta}(t,\cdot)$ is the density of a particle diffusing in the
flow given by $u_{\theta}$, with unit diffusivity. A particular consequence
of \prettyref{eq:viscouscase} (along with a moment bound on $g_{\theta}$
proved in \cite{DG23}) is that $u_{\theta}(t,x)$ is continuous in
$\theta$, which is in sharp contrast to the behavior that has been
observed on the real line.

The purpose of the present paper is to study the jointly invariant
measures of the \emph{inviscid} periodic Burgers equation \prettyref{eq:uproblem}
simultaneously for all $\theta$. The main results we will prove can
be seen as analogues of those in \cite{DG23}. In particular, we will
show that jointly stationary solutions for \prettyref{eq:uproblem}
exist, that a one-force-one-solution principle holds, and that an
inviscid analogue of the relation \prettyref{eq:viscouscase} holds.
It will no longer be the case that $u_{\theta}(t,x)$ is continuous
in $\theta$ because the analogue of $g_{\theta}$ is less regular
in our setting. However, the spatial integral of $g_{\theta}$ will
be continuous in $x$, and hence $\int_{0}^{x}u_{\theta}(t,y)\,\dif y$
will still be continuous in $\theta$. Actually, this is a rather
generic feature of the periodic case: if we define $h_{\theta}(t,x)=\int_{0}^{x}u_{\theta}(t,x)\,\dif x$,
then by the comparison principle we also have for $\theta_{2}\ge\theta_{1}$
that $u_{\theta_{2}}\ge u_{\theta_{1}}$ (assuming a similar ordering
for the initial condition), and so for $x\in[0,1]$ we have 
\begin{equation}
0\le h_{\theta_{2}}(t,x)-h_{\theta_{1}}(t,x)\le h_{\theta_{2}}(t,1)-h_{\theta_{1}}(t,1)\overset{\prettyref{eq:uintegral}}{=}\theta_{2}-\theta_{1}.\label{eq:hsqueeze}
\end{equation}

In the inviscid setting we consider here, the particle diffusing in
the Burgers flow should be replaced by a particle simply moving in
the flow, without diffusion. It is well-known that such a particle
will eventually join a \emph{shock }of the Burgers flow. This suggests
that, at least formally, the density $g_{\theta}$ should be replaced
by a delta mass at a single ``global'' shock, and \prettyref{eq:viscouscase}
then suggests that the change in $u$ as $\theta$ is varied should
occur only at the location of this global shock. We will prove a precise
version of this statement in our main theorem \prettyref{thm:mainthm}
below. We also point the reader to the video included the supplementary
material for a visualization of the shocks in the Burgers flow as
$\theta$ is changed; see \prettyref{subsec:illustrations} for a
description.

The proof techniques in the present setting are entirely different
from those of \cite{DG23}. In particular, we study the minimizers
of the functional \prettyref{eq:Adef}, and the associated shocks
as mentioned above, rather than using the stochastic analysis tools
of \cite{DG23}. This leads us to a fine study of the structure of
one-sided minimizers and global shocks for the inviscid Burgers equation.
One-sided minimizers and global shocks have been studied extensively
in the literature on the stochastic Burgers equation. We refer to
the survey \cite{BK18} for an illuminating heuristic discussion.
The use of minimizers rather than polymers makes many features of
the problem more explicit than in the viscous case, and certain aspects
of the phenomenology are clearer. In particular, we will see how the
topology of the minimizers plays an important role in the analysis.
The importance of the topology of the minimizers has been previously
observed in the physics literature in \cite{BIK02,BK07} for the Burgers
equation with time-periodic forcing that is smoother in time than
we consider here.

A consistent theme of work on stochastic Burgers with multiple means
considered simultaneously is the presence of exceptional values of
$\theta$ at which behavior is observed that happens with probability
zero for any fixed $\theta$. This holds as well in our setting. Thus,
our study of minimizers and shocks will go beyond that of previous
work in that we will prove the behavior at these exceptional values
of $\theta$ as well. In particular, an important feature of the study
of minimizers for the Burgers equation is that distinct minimizers
from the same point do not typically cross each other. At exceptional
values of $\theta$, a certain amount of crossing is possible, and
thus the picture that we describe exhibits significant additional
topological complexity compared to the fixed-$\theta$ case.

\subsection{Main results}

We now state precisely the main results of our study. Although we
are primarily interested in the setting in which $\mu$ is a Poisson
point process, we can actually state simple deterministic conditions
on $\mu$ under which our results hold. We will then check (see \prettyref{thm:poisson-ok})
that these conditions are satisfied with probability $1$ by any homogeneous
compound Poisson point process. To facilitate this, we introduce the
following space of forcing measures.
\begin{defn}
Let $\Omega$ be the space of purely atomic measures $\mu$ on $\RR\times\TT$
such that $\ssN\coloneqq\supp\mu$ is a discrete set.
\end{defn}

Since much of our work concerns the behavior of Lagrangian minimizers,
we first define them precisely.
\begin{defn}[Lagrangian minimizers]
\label{def:minimizers}Let $\theta\in\RR$ and $s<t$.
\begin{enumerate}
\item We define the set $\scM_{s,y\mid t,x}^{\theta}$ comprising all paths
$X\in H^{1}([s,t];\TT)$ with $X(s)=y$ and $X(t)=x$ such that if
$Y$ is another such path, then $\mathcal{A}_{\theta,s,t}[X]\le\mathcal{A}_{\theta,s,t}[Y]$.
\item We define the set of \emph{one-sided minimizers}
\[
\scM_{t}^{\theta}\coloneqq\{X\in H_{\mathrm{loc}}^{1}((-\infty,t];\TT)\st X|_{[s,t]}\in\scM_{s,X(s)\mid t,X(t)}^{\theta}\text{ for all }s\le t\}.
\]
\item For $x\in\TT$, we define
\[
\scM_{t,x}^{\theta}\coloneqq\{X\in\scM_{t}^{\theta}\st X(t)=x\}.
\]
\item We define the partial order $\preccurlyeq$ on $\scM_{t,x}^{\theta}$
by
\begin{equation}
X_{1}\preccurlyeq X_{2}\text{ if }\int_{r}^{t}X_{1}'(s)\,\dif s\ge\int_{r}^{t}X_{2}'(s)\,\dif s\text{ for all }r\le t.\label{eq:partialorderdef}
\end{equation}
\end{enumerate}
\end{defn}

The partial order $\preccurlyeq$ represents the ordering of minimizers
when lifted to the universal cover $\RR$ of $\RR/\ZZ$. The reason
for the ``$\ge$'' sign in \prettyref{eq:partialorderdef} is that
we say that $X_{1}\preccurlyeq X_{2}$ if the graph of $X_{1}$ lies
(not necessarily strictly) to the left of that of $X_{2}$ when lifted
to the universal cover, which corresponds to the integral of the derivative
of $X_{1}$ being at least that of $X_{2}$. This ordering will play
an important role in our topological arguments.

Our first theorem will state that one-sided minimizers exist. To state
it, we first need to introduce another assumption on the noise field:
\begin{defn}
\label{def:smallnoisezones}We define $\widetilde{\Omega}_{1}$ to
be the set of all $\mu\in\Omega$ such that there exists an $M=M(\mu)>0$
such that for all $t\in\RR$, there is an $s=s(t)\le t-M$ and a $y\in\TT$
such that 
\begin{equation}
\ssN\cap([s-M,s+M]\times\TT)=\{(s,y)\}\label{eq:intersectiononly}
\end{equation}
and
\begin{equation}
\mu(\{(s,y)\})>\frac{1}{4M},\label{eq:youshouldtake}
\end{equation}
and moreover that we can choose $s(t)$ in such a way that
\begin{equation}
\lim_{t\to+\infty}s(t)=+\infty.\label{eq:stttoinfty}
\end{equation}
\end{defn}

\prettyref{def:smallnoisezones} encodes the notion of \emph{small-noise
zones}: regions of space-time containing only a single, sufficiently
strong forcing point. All minimizers started sufficiently far in the
future and extending sufficiently far into the past must pass through
such points. Indeed, since there are no other forcing points in the
vicinity, there is nothing to be lost by passing through the forcing
point; see \prettyref{prop:gothroughy} below. These small-noise zones
act as regeneration times for the dynamics, as the behavior of the
minimizers before and after them is decoupled. They have already been
mentioned for this problem in \cite{Bak13}, and the simple new observation
here is that the decoupling happens simultaneously over all $\theta$.

Given the existence of small-noise zones, we can prove the existence
of global solutions to \prettyref{eq:uproblem}:
\begin{thm}
\label{thm:minimizers-exist}Suppose that $\mu\in\widetilde{\Omega}_{1}$.
For each $\theta,t\in\RR$ and $x\in\TT$, the set $\scM_{t,x}^{\theta}$
is nonempty, consists of piecewise-linear paths connecting $(t,x)$
and points of $\ssN$, and has unique minimal and maximal elements
$X_{\theta,t,x,\Left}$ and $X_{\theta,t,x,\Right}$, respectively,
under the partial order $\preccurlyeq$.
\end{thm}

Again, we note that, in the case when $\mu$ is a Poisson point process,
\prettyref{thm:minimizers-exist} has been proved for a single $\theta$
at a time in \cite{Bak13}. The novelty here (even in the case when
$\mu$ is a Poisson point process) is that we prove it for all $\theta$
on a single set $\widetilde{\Omega}_{1}$ (which will have probability
$1$ when $\mu$ is a Poisson point process). \prettyref{thm:minimizers-exist}
is proved as part of \prettyref{prop:globalminimizersexist} in \prettyref{sec:Existence}.

Let $\LR=\{\Left,\Right\}$. We define, for $\Square\in\LR$,
\begin{equation}
u_{\theta,\Square}(t,x)=X_{\theta,t,x,\Square}'(t-),\label{eq:udef}
\end{equation}
the derivative from below of the process $X_{\theta,t,x,\Square}$
at $t$. Both $u_{\theta,\Left}$ and $u_{\theta,\Right}$ are global-in-time
entropy solutions to \prettyref{eq:uproblem}. It is clear from the
definitions that $u_{\theta,\Left}(t,x)\ge u_{\theta,\Right}(t,x)$
for all $(\theta,t,x)\in\RR\times\RR\times\TT$. They differ only
on the set
\begin{equation}
\ssS\coloneqq\{(\theta,t,x)\in\RR\times\RR\times\TT\st u_{\theta,\Left}(t,x)>u_{\theta,\Right}(t,x)\},\label{eq:Sdef}
\end{equation}
which is the set of \emph{shocks}. We define
\begin{equation}
\ssS_{\theta}\coloneqq\{(t,x)\in\RR\times\RR\times\TT\st(\theta,t,x)\in\ssS\}\label{eq:Sthetadef}
\end{equation}
and
\begin{equation}
\ssS_{\theta,t}\coloneqq\{x\in\TT\st(t,x)\in\ssS_{\theta}\}.\label{eq:Sthetatdef}
\end{equation}

\prettyref{thm:minimizers-exist} describes a picture (previously
observed in \cite{Bak13}) in which, for each fixed $\theta$, the
time-space cylinder $\RR\times\TT$ is tessellated by regions of points
$(t,x)$ for which the last forcing point on $X_{\theta,t,x,\Square}$
is a given forcing point. The shock set $\ssS_{\theta}$ is formed
from the boundaries of these regions (except that points of $\ssN$
are generally on the boundaries of these regions but not in $\ssS_{\theta}$).
The dynamics of shocks are also well-understood. Each forcing point
creates a pair of shocks starting at that point. (See \prettyref{prop:advanceint-N}.)
A shock at position $(t,x)$ moves (as time advances) with velocity
$\frac{1}{2}(u_{\theta,\Left}(t,x)+u_{\theta,\Right}(t,x))$ (the
well-known \emph{Rankine--Hugoniot condition}; see \prettyref{prop:advanceint}
for a proof in our setting). When two shocks collide with one another,
they merge to form a single shock. See \prettyref{fig:shocks_picture}
for an illustration.
\begin{figure}[t]
\begin{centering}
\input{figures/shocks_picture.pgf}
\par\end{centering}
\caption[Forcing points, minimizers, shocks, and the Burgers solution.]{The bottom plot shows the forcing points (yellow), a sample of minimizers
(black dashed lines), and the shocks (blue solid lines), which are
points at which minimizers extend in multiple directions. The top
plot shows $u_{\theta}(t,\cdot)=u_{\theta,\Square}(t,\cdot)$, where
$\Square\in\protect\LR$ is arbitrary from the point of view of plotting.
Note that the discontinuities of $u_{\theta}(t,\cdot)$ correspond
to the locations of the blue shock curves at time $t$.\label{fig:shocks_picture}}
\end{figure}

It has been observed in \cite{EKMS00} for a fixed value of $\theta$
(with a differnt type of forcing) that, with probability $1$, there
is a unique \emph{global minimizer} and a unique \emph{global shock}
(the latter also known as the \emph{main shock} or \emph{topological
shock}). All minimizers merge with the global minimizer as time goes
to $-\infty$, and all shocks merge with the global shock as time
goes to $+\infty$. (In the forcing considered in \cite{EKMS00},
the shocks converge towards each other exponentially fast rather than
literally merging.) The term \emph{topological shock }is particularly
illuminating; it refers to the fact that the global shock is characterized
by the presence of minimizers that, at the point that they merge back
together, have between them completed a nontrivial winding about the
torus. See \cite[Theorem 5.2]{EKMS00}.

Let us now state this topological characterization of global shocks
precisely. For $(\theta,t,x)\in\ssS$, we define
\begin{equation}
T_{\vee}(\theta,t,x)\coloneqq\sup\{s<t\st X_{\theta,t,x,\Left}(s)=X_{\theta,t,x,\Right}(s)\}.\label{eq:Tveedef}
\end{equation}

\begin{defn}
\label{def:globalshocks}We define the set $\GS$ of \emph{global
shocks} comprising all shocks $(\theta,t,x)\in\ssS$ such that
\begin{equation}
\int_{T_{\vee}(\theta,t,x)}^{t}X_{\theta,t,x,\Left}'(s)\,\dif s>\int_{T_{\vee}(\theta,t,x)}^{t}X_{\theta,t,x,\Right}'(s)\,\dif s.\label{eq:shockslopesdiffernt}
\end{equation}
We also define $\GS_{\theta}\coloneqq\{(t,x)\in\ssS_{\theta}\st(\theta,t,x)\in\GS\}$
and $\GS_{\theta,t}\coloneqq\{x\in\ssS_{\theta,t}\st(\theta,t,x)\in\GS\}$.
\end{defn}

\begin{figure}[t]
\begin{centering}
\input{figures/global_shock_picture.pgf}
\par\end{centering}
\caption[The difference between global and non-global shocks]{The global shock curve is shown as a thicker light blue line. The
point $(\theta,t,y)$ is a global shock, since the left and right
minimizer coming from $(t,y)$, considered up until their first meeting
point, together complete a wrap around the torus before meeting again.
The point $(\theta,t,x)$ is a shock, since minimizers come from $x$
in multiple directions, but not a global shock, since the minimizers
do not accumulate a nontrivial winding before meeting again.\label{fig:global_shock_picture}}
\end{figure}
See \prettyref{fig:global_shock_picture} for an illustration of the
definition of global shock. Additional pictures are available in \cite{BIK02,BK07},
in particular in the higher-dimensional setting, which is also of
significant interest but which we do not consider here.

When we consider all values of $\theta$ simultaneously, it is \emph{not
}the case that there exists a unique global shock for each $\theta,t$.
In fact, this is impossible, since the global shock must have asymptotic
slope $\theta$, and so it cannot vary continuously as $\theta$ is
varied. We can nonetheless make a strong statement about the structure
of the global shock set $\GS$ if we make the following assumption
on the noise, which holds with probability $1$ for the Poisson point
process.
\begin{defn}
\label{def:Omegatilde2}We define $\widetilde{\Omega}_{2}$ to be
the set of all $\mu\in\Omega$ such that, for each $\theta\in\RR$,
we have
\begin{equation}
\#(\GS_{\theta}\cap\ssN)\le1\label{eq:GSNlt1}
\end{equation}
and
\begin{equation}
\ssS_{\theta}\cap\ssN\setminus\GS_{\theta}=\varnothing.\label{eq:allshockNsglobal}
\end{equation}
We define
\begin{equation}
\Theta_{\otimes}\coloneqq\{\theta\in\RR\st\#(\GS_{\theta}\cap\ssN)=1\},\label{eq:Theta1def}
\end{equation}
and we define maps $s_{\otimes}\colon\Theta_{\otimes}\to\RR$ and
$y_{\otimes}\colon\Theta_{\otimes}\to\TT$ by letting $(s_{\otimes}(\theta),y_{\otimes}(\theta))$
be the unique element of $\GS_{\theta}\cap\ssN$ for each $\theta\in\Theta_{\otimes}$.

We further define 
\[
\widetilde{\Omega}=\widetilde{\Omega}_{1}\cap\widetilde{\Omega}_{2}.
\]
\end{defn}

\begin{thm}
\label{thm:poisson-ok}Let $\mathbb{P}$ be the probability measure
associated to a homogeneous compound Poisson point process on $\mathbb{R}\times\mathbb{T}$.
Then $\mathbb{P}(\widetilde{\Omega})=1$.
\end{thm}

We prove \prettyref{thm:poisson-ok} in \prettyref{sec:poisson}.
Now we can state our second result on the structure of the global
shock set.
\begin{thm}
\label{thm:shockstructure}Suppose that $\mu\in\widetilde{\Omega}$.
There are unique functions $\sfs_{\Left},\sfs_{\Right}\colon\RR\times\RR\to\TT$
(which we call the \emph{left} and \emph{right global shocks}, respectively)
such that the following properties hold:
\begin{enumerate}
\item \label{enu:ssquareshocks}For each $\theta,t\in\RR$ and $\Square\in\LR$,
we have $\GS_{\theta,t}=\{\sfs_{\Left}(\theta,t),\sfs_{\Right}(\theta,t)\}$.
\item \label{enu:shocksleftright}For each fixed $t\in\RR$ and $\Square\in\LR$,
the function $\theta\mapsto\sfs_{\Square}(\theta,t)$ is piecewise
continuous. In particular, for each fixed $t\in\RR$,
\begin{equation}
\text{the set \ensuremath{\{\theta\in\RR\st\sfs_{\Left}(\theta,t)\ne\sfs_{\Right}(\theta,t)\}} is discrete,}\label{eq:sLnotsRdiscrete}
\end{equation}
and
\begin{equation}
\sfs_{\Left}(\theta,t)=\lim_{\theta'\uparrow\theta}\sfs_{\Diamond}(\theta,t)\qquad\text{and}\qquad\sfs_{\Right}(\theta,t)=\lim_{\theta'\downarrow\theta}\sfs_{\Diamond}(\theta,t)\label{eq:sLsRareleftrightlimits}
\end{equation}
for each $\theta\in\RR$ and $\Diamond\in\LR$.
\end{enumerate}
Moreover, these functions have the following additional properties
(which are not necessary for the uniqueness statement):
\begin{enumerate}[resume]
\item \label{enu:shockcontinuity}For each fixed $\theta\in\RR$ and $\Square\in\LR$,
the map $t\mapsto\sfs_{\Square}(\theta,t)$ is continuous.
\item \label{enu:splitatforcing}If $\theta\in\RR\setminus\Theta_{\otimes}$,
then $\sfs_{\Left}(\theta,t)=\sfs_{\Right}(\theta,t)$ for all $t\in\RR$.
On the other hand, if $\theta\in\Theta_{\otimes}$, then there exists
an $s_{\wedge}(\theta)\in(s_{\otimes}(\theta),\infty)$ such that
\[
\{t\in\RR\st\sfs_{\Left}(\theta,t)\ne\sfs_{\Right}(\theta,t)\}=(s_{\otimes}(\theta),s_{\wedge}(\theta)).
\]
\end{enumerate}
\end{thm}

The functions $\sfs_{\Left}$ and $\sfs_{\Right}$ are constructed
in \prettyref{def:sLsRdefs}, which gives a characterization that
is in many ways easier to work with than the characterization given
in \prettyref{thm:shockstructure}, but requires some additional definitions
to state. Part~\ref{enu:ssquareshocks} of \prettyref{thm:shockstructure}
is implicit in that definition. The statement \prettyref{eq:sLnotsRdiscrete}
is proved as \prettyref{prop:shocksdiffdiscrete}, and \prettyref{eq:sLsRareleftrightlimits}
is proved simultaneously with \prettyref{prop:relatetooneoverderiv}
in \prettyref{subsec:Movement-of-shocks}. Part~\ref{enu:shockcontinuity}
of the theorem statement is proved as \prettyref{prop:smovement}(\ref{enu:scts}),
and part~\ref{enu:splitatforcing} is proved in \prettyref{subsec:movement-global}.
The uniqueness statement in \prettyref{thm:shockstructure} holds
because part~\ref{enu:ssquareshocks} and \prettyref{eq:sLnotsRdiscrete}
characterize, for each $t$, $\mathrm{s}_{\mathrm{L}}(\theta,t)$
and $\mathrm{s}_{\mathrm{R}}(\theta,t)$ at all except a discrete
set of $\theta$, which means that the limits in \prettyref{eq:sLsRareleftrightlimits}
are well-defined and characterize $\mathrm{s}_{\mathrm{L}}(\theta,t)$
and $\mathrm{s}_{\mathrm{R}}(\theta,t)$ for all $\theta$.

The discontinuity of $\sfs_{\Square}(\theta,t)$ is in accordance
with the topological obstruction to the continuity of $\sfs_{\Square}(\theta,t)$
in $\theta$ mentioned above. The phenomenon we observe is that, generically,
we have $\sfs_{\Left}(\theta,t)=\sfs_{\Right}(\theta,t)$: a single
global shock for each $\theta$ and $t$. However, for $\theta\in\Theta_{\otimes}$,
there is a time $s_{\otimes}(\theta)$ at which the global shock hits
a forcing point and splits into two global shocks $\sfs_{\Left}(\theta,t)\ne\sfs_{\Right}(\theta,t)$.
These global shocks then re-merge at a later time $s_{\wedge}(\theta)$,
but they may have accumulated a nontrivial winding relative to one
another by this merging time. By this last statement we mean that
the union of the two branches is not contractible; see \prettyref{fig:global_shock_split_picture}.
\begin{figure}[t]
\begin{centering}
\input{figures/global_shock_split_picture.pgf}
\par\end{centering}
\caption[The global shock splitting and re-merging]{The global shock splitting into two when it hits a forcing point.
Later, the two global shocks merge back together. The two branches
accumulate different a nontrivial winding relative to one another
by the time they re-merge.\label{fig:global_shock_split_picture}}
\end{figure}
 This motivates the conditions in \prettyref{def:Omegatilde2}, which
state that for each $\theta$, there can be at most one forcing point
that lies on a shock. For fixed $\theta$, this happens with probability
$0$, but with probability $1$, it will happen for some values of
$\theta$.

The global shock set is the analogue of $g_{\theta}$ defined in \prettyref{eq:gproblem}
for the viscous problem. Indeed, $g_{\theta}$ represents the density
of a passive particle that has been diffusing in the Burgers flow
since time $-\infty$. A similar particle moving in the inviscid Burgers
flow (without diffusivity, in accordance with the inviscidity) will
end up in the global shock set, since all particles eventually merge
with a shock and all shocks eventually merge with the global shock
(up to the fact that the global shock may itself split at most once
for each $\theta$).

We are now ready to state our main theorem, which describes how the
global solutions to the Burgers equation change as $\theta$ is varied.
\begin{thm}
\label{thm:mainthm}Suppose that $\mu\in\widetilde{\Omega}$ and fix
$t\in\RR$.
\begin{enumerate}
\item \label{enu:derivintheta}The functions $\sfs_{\Left}(\cdot,t)$ and
$\sfs_{\Right}(\cdot,t)$ are left- and right-differentiable, respectively.
For any $\theta_{1}<\theta_{2}$ and any $\Square\in\LR$, we have
\begin{equation}
u_{\theta_{2},\Square}(t,x)-u_{\theta_{1},\Square}(t,x)=\sum_{\substack{\theta\in\llbracket\theta_{1},\theta_{2}\rrbracket_{\Square}\\
\sfs_{\Square}(\theta,t)=x
}
}\frac{1}{\partial_{\theta}\sfs_{\Square}(\theta\pm_{\Square},t)},\label{eq:utheta2theta1}
\end{equation}
where we have defined
\begin{equation}
\llbracket\theta_{1},\theta_{2}\rrbracket_{\Square}=\begin{cases}
(\theta_{1},\theta_{2}] & \text{if }\Square=\Left;\\{}
[\theta_{1},\theta_{2}) & \text{if }\Square=\Right
\end{cases}\qquad\text{and}\qquad\pm_{\Square}=\begin{cases}
- & \text{if }\Square=\Left;\\
+ & \text{if }\Square=\Right.
\end{cases}\label{eq:pmbookkeeping}
\end{equation}
\item \label{enu:switching}If $\sfs_{\Left}(\theta,t)\ne\sfs_{\Right}(\theta,t)$,
then there is an $\eps=\eps(\theta,t)>0$ such that 
\begin{equation}
(\theta',t,\sfs_{\Left}(\theta,t))\in\ssS\qquad\text{for all }\theta'\in[\theta,\theta+\eps)\label{eq:leftstillshock}
\end{equation}
and
\begin{equation}
(\theta',t,\sfs_{\Right}(\theta,t))\in\ssS\qquad\text{for all }\theta'\in(\theta-\eps,\theta].\label{eq:rightstillshock}
\end{equation}
\end{enumerate}
\end{thm}

\begin{figure}[p]
\centering{}%
\foreach\mytheta in {-0.51,-0.49,-0.47,-0.45,-0.43,-0.41,-0.39,-0.37,-0.35,-0.33,-0.31,-0.29,-0.27,-0.25,-0.23,-0.21,-0.19,-0.17}{\input{figures/global_shock_moving_picture-\mytheta.pgf}\hfill{}}\caption[Plots of the global shock set for a sequence of values of $\theta$]{\label{fig:global_shock_moving_picture}Plots of the global shock
set for a sequence of values of $\theta$, which increase as the figures
are read like English text, top to bottom and left to right.}
\end{figure}
Let us now describe how \prettyref{eq:utheta2theta1} is an inviscid
analogue of \prettyref{eq:viscouscase}. Suppose for sake of illustration
that there is just a single global shock $\sfs(\theta,t)$ for each
$\theta$ and $t$, and that it is differentiable as a function of
$\theta$. Under this (incorrect) assumption, the analogy of \prettyref{eq:viscouscase}
would be
\begin{equation}
u_{\theta_{2}}(t,x)-u_{\theta_{1}}(t,x)=\int_{\theta_{1}}^{\theta_{2}}\delta(x-\sfs(\theta,t))\,\dif\theta.\label{eq:viscouscase-delta}
\end{equation}
Here we have used the fact that the inviscid analogue of $g_{\theta}$
is a delta mass at $\sfs(\theta,t)$, as discussed above. Formally
performing the change of variables $y=\sfs(\theta,t)$, $\dif y=\partial_{\theta}\sfs(\theta,t)\dif\theta$
in \prettyref{eq:viscouscase-delta}, we get
\begin{equation}
u_{\theta_{2}}(t,x)-u_{\theta_{1}}(t,x)=\int_{\sfs(\theta_{1},t)}^{\sfs(\theta_{2},t)}\sum_{\substack{\theta\in[\theta_{1},\theta_{2}]\\
\sfs(\theta,t)=y
}
}\frac{\delta(x-y)}{\partial_{\theta}\sfs(\theta,t)}\,\dif y=\sum_{\substack{\theta\in[\theta_{1},\theta_{2}]\\
\sfs(\theta,t)=x
}
}\frac{1}{\partial_{\theta}\sfs(\theta,t)},\label{eq:viscouscase-delta-1}
\end{equation}
which is almost but not precisely correct. Our main result \prettyref{eq:utheta2theta1}
is a corrected version of the statement: it takes into the account
that in the inviscid case, neither $\sfs_{\Square}(\theta,t)$ nor
$u_{\theta,\Square}(t,x)$ will be continuous in $\theta$, and selects
the left or right versions of these functions as appropriate.

\prettyref{thm:mainthm} implies that, if $t$ and $x$ are held fixed
and $\theta$ varies, then the solutions $u_{\theta,\Square}(t,x)$
only change when $\sfs_{\Square}(\theta,t)=x$. \prettyref{def:globalshocks}
makes transparent the reason for this phenomenon. Indeed, if we differentiate
\prettyref{eq:Adef} in $\theta$, we get
\[
\frac{\dif}{\dif\theta}\mathcal{A}_{\theta,s,t}[X]=-\int_{s}^{t}X'(r)\,\dif r+\theta(t-s).
\]
Therefore, roughly speaking, we see a change in the slope of the minimizer
only when there are multiple minimizers with different values of the
integral of $X'$ from $t$ until the point at which the minimizers
meet, which is exactly what is encoded in \prettyref{eq:shockslopesdiffernt}.

Part~\ref{enu:switching} of \prettyref{thm:mainthm} further develops
the behavior of the jumps of $\sfs_{\Square}(\theta,t)$. Indeed,
it shows that when $\theta\mapsto\sfs_{\Square}(\theta,t)$ has a
jump, it jumps to a preexisting non-global shock, which then becomes
a global shock and begins to move. See \prettyref{fig:global_shock_moving_picture}
and also the animation in the supplementary material (described in
\prettyref{subsec:illustrations}).

\subsection{Organization of the paper}

In \prettyref{sec:Existence}, we establish basic facts about the
minimizers, and in particular prove \prettyref{thm:minimizers-exist}.
In \prettyref{sec:Continuity-of-minimizers}, we set the stage for
our perturbative arguments by using the discreteness of $\ssN$ to
show that, as the parameters $\theta,t,x$ are changed, the points
used in the minimizer can only change when there are multiple minimizers
coming from the same point. In \prettyref{sec:global-shocks}, we
study the topological features of global shocks \prettyref{def:globalshocks}.
In \prettyref{sec:Movement-of-shocks}, we study how shocks, and in
particular global shocks, move as $t$ is varied. In \prettyref{sec:udependenceontheta},
we study how the flow changes as $\theta$ is changed. Finally, in
\prettyref{sec:poisson}, we prove \prettyref{thm:poisson-ok}, showing
that our assumptions are satisfied almost surely for a compound Poisson
process.

\subsection{Key to illustrations}

\label{subsec:illustrations}\typeout{(./anc/moving_global_shock.mp4)}The
paper features several figures and is also accompanied by one animation.
In each of the figures, minimizers are drawn as dotted black lines,
forcing points as yellow dots with black borders, non-global shocks
by dark blue solid lines, and global shocks by thicker light blue
lines. Time increases along the vertical axis and space is drawn along
the horizontal axis.

In the animation, which is included in the supplementary material,
the value of $\theta$ starts negative and is increased as the animation
progresses. The shocks, global shocks, and a sampling of minimizers
are drawn, and the graph of $u_{\theta}(t,\cdot)$ for the last plotted
time $t$ is also shown above as in \prettyref{fig:shocks_picture}.
The reader will note that, as proved in \prettyref{thm:mainthm},
the only movement in the picture is through the movement of the global
shock. When the global shock hits a forcing point, it jumps to the
shock extending from the other side of that forcing point.

\subsection{Acknowledgments}

The author is very grateful to Yuri Bakhtin, Yu Gu, and Evan Sorensen
for many inspiring discussions and important feedback throughout the
project, and also to an anonymous referee for helpful remarks on the
manuscript.

\section{Existence of one-sided minimizers}

\label{sec:Existence}In this section we prove \prettyref{thm:minimizers-exist}
on the existence of one-sided minimizers. First we must establish
some basic properties of minimizers. If $X\in H^{1}([s,r];\TT)$ and
$Y\in H^{1}([r,t];\TT)$ are such that $X(r)=Y(r)$, we define the
concatenation $X\odot_{r}Y\colon[s,t]\to\TT$ by
\begin{equation}
(X\odot_{r}Y)(q)\coloneqq\begin{cases}
X(q), & \text{if }q\le r;\\
Y(q), & \text{if }q\ge r.
\end{cases}\label{eq:odotdef}
\end{equation}

\begin{prop}
\label{prop:basicproperties}Suppose that $\mu\in\Omega$. Let $\theta\in\RR$,
$-\infty<s<t<+\infty$, and $x,y\in\TT$.
\begin{enumerate}
\item \label{enu:straightlines}The set $\scM_{s,y\mid t,x}^{\theta}$ is
nonempty. Every $X\in\scM_{s,y\mid t,x}^{\theta}$ consists of straight
line segments connecting points of $\{(t,x),(s,y)\}\cup\ssN$.
\item \label{enu:subpathalsominimizes}If $\varnothing\ne[s',t']\subseteq[s,t]$
and $X\in\mathscr{M}_{s,y\mid t,x}^{\theta}$, then $X|_{[s',t']}\in\scM_{s',X(s')\mid t',X(t')}^{\theta}$
as well.
\item \label{enu:pathsurgery}If $r\in(s,t)$, $X\in\scM_{s,y\mid t,x}^{\theta}$
and $Y\in\scM_{s,y\mid r,X(r)}^{\theta}$, then $Y\odot_{r}X\in\scM_{s,y\mid t,x}^{\theta}$
as well.
\end{enumerate}
\end{prop}

\begin{proof}
The first point is a standard property of the convexity of the Dirichlet
energy in \prettyref{eq:Adef}. For the second point, we note that
if not, then we could modify $X$ on $[s',t']$ to improve the value
of $\mathcal{A}_{\theta,s,t}$, contradicting the assumption that
$X\in\scM(\theta,s,t)$. To see the third point, note that 
\[
\mathcal{A}_{\theta,s,t}[Y\odot_{r}X]=\mathcal{A}_{\theta,s,r}[Y]+\mathcal{A}_{\theta,r,t}[X]\le\mathcal{A}_{\theta,s,r}[X]+\mathcal{A}_{\theta,s,r}[X]=\mathcal{A}_{\theta,s,t}[X]
\]
by the definitions and part~\ref{enu:subpathalsominimizes}.
\end{proof}
To prove \prettyref{thm:minimizers-exist}, we use the existence of
small-noise zones described in \prettyref{def:smallnoisezones} to
achieve decoupling. The point is that when a small-noise zone occurs,
the behavior of polymers inside the small-noise zone is independent
of what happens outside of the small zone. This then implies that
the behaviors of the polymer before and after the small-noise zone
are conditionally independent. The following proposition, whose statement
appeared already in \cite{Bak13} in the case when $\mu(\{(t,x)\})=1$
for all $(t,x)\in\ssN$, is the reason for the definition of $\widetilde{\Omega}_{1}$.
\begin{prop}
\label{prop:gothroughy}Suppose that $\mu\in\Omega$ and that $s\in\RR$,
$M>0$, and $y\in\TT$ are such that \prettyref{eq:intersectiononly}
and \prettyref{eq:youshouldtake} hold. Then, for any $\theta\in\RR$,
$z_{1},z_{2}\in\TT$, and $X\in\scM_{s-M,z_{1}\mid s+M,z_{2}}^{\theta}$,
we have $X(s)=y$.
\end{prop}

\begin{proof}
Assume towards a contradiction that we have $\theta\in\RR$ and $X\in\scM_{s-M,z_{1}\mid s+M,z_{2}}^{\theta}$
such that $X(s)\ne y$. By the assumptions on $s$ and $y$, along
with \prettyref{prop:basicproperties}(\ref{enu:straightlines}),
we see that $X$ consists of a single straight line segment. Let 
\begin{equation}
\xi\in(-1/2,1/2]\label{eq:xirange}
\end{equation}
be such that $X(s)+\xi=y$, and define
\[
Y(r)\coloneqq X(r)+\xi\cdot\begin{cases}
\frac{s+M-r}{M} & \text{if }r\in[s,s+M];\\
\frac{r-s+M}{M} & \text{if }r\in[s-M,s];\\
0 & \text{otherwise}.
\end{cases}
\]
In particular, this means that $Y(s-M)=X(s-M)=z_{1}$, $Y(s+M)=X(s+M)=z_{2}$,
and $Y(s)=X(s)+\xi=y$. Then we have
\begin{align*}
 & \mathcal{A}_{\theta,s-M,s+M}[Y]-\mathcal{A}_{\theta,s-M,s+M}[X]\\
 & \qquad=M\left[\frac{1}{2}(X'(s)-\theta+\xi/M)^{2}+\frac{1}{2}(X'(s)-\theta-\xi/M)^{2}-(X'(s)-\theta)^{2}\right]-\mu(\{(s,y)\})\\
 & \qquad=\frac{\xi^{2}}{M}-\mu(\{(s,y)\})\overset{\prettyref{eq:xirange}}{\le}\frac{1}{4M}-\mu(\{(s,y)\})\overset{\prettyref{eq:youshouldtake}}{<}0.
\end{align*}
But since $Y(s\pm M)=X(s\pm M)$, this contradicts the assumption
$X\in\scM_{s-M,z_{1}\mid s+M,z_{2}}^{\theta}$.
\end{proof}
Now we can make the following important definition.
\begin{defn}
\label{def:Tstar}Suppose that $\mu\in\widetilde{\Omega}_{1}$ and
let $M(\mu)$ be as in \prettyref{def:smallnoisezones}. For $t\in\RR$,
we define $T_{*}(t)$ to be the supremum of all $s<t$ such that there
exists a $y\in\TT$ such that $(s,y)\in\ssN$ and $X(s)=y$ for all
$X\in\bigcup\limits_{\substack{\theta\in\RR\\
x,z\in\TT
}
}\scM_{s-M(\mu),z\mid t,x}^{\theta}$.
\end{defn}

We note that it is an immediate consequence of \prettyref{def:Tstar}
and the definition \prettyref{eq:Tveedef} of $T_{\vee}$ that
\begin{equation}
T_{*}(t)\le T_{\vee}(\theta,t,x)\qquad\text{for any }\theta,t\in\RR\text{ and }x\in\TT.\label{eq:TstarTvee}
\end{equation}
The point of \prettyref{prop:gothroughy} is that $T_{*}(t)$ is finite: 
\begin{prop}
\emph{\label{prop:Tstartfinite}}Suppose that $\mu\in\widetilde{\Omega}_{1}$
and $t\in\RR$.
\begin{enumerate}
\item We have $T_{*}(t)>-\infty$.
\item There is a $y_{*}(t)\in\TT$ such that
\begin{equation}
X(T_{*}(t))=y_{*}(t)\text{ for all }X\in\bigcup\limits_{\substack{\theta\in\RR\\
x,z\in\TT
}
}\scM_{T_{*}(t)-M(\mu),z\mid t,x}^{\theta}.\label{eq:allpassthrough}
\end{equation}
\item Finally, we have
\begin{equation}
\lim_{t\to+\infty}T_{*}(t)=+\infty.\label{eq:Tstarttoinfinity}
\end{equation}
\end{enumerate}
\end{prop}

\begin{proof}
Fix $t\in\RR$. Let $M>0$, $s\in(-\infty,t-M]$, and $y\in\TT$ be
as in the definition of $\widetilde{\Omega}$. We claim that in fact
$T_{*}(t)\ge s$. \prettyref{prop:basicproperties}(\ref{enu:subpathalsominimizes})
implies that if $\theta\in\RR$, $x,z\in\TT$, and $X\in\scM_{s-M,z\mid t,x}^{\theta}$,
then $X\in\scM_{s-M,z\mid s+M,X(s+M)}^{\theta}$ as well, and so by
\prettyref{prop:Tstartfinite} we have $X(s)=y$. Since $(s,y)\in\ssN$
by definition, we have $T_{*}(t)\ge s>-\infty$, and the first assertion
is proved. The second assertion is tantamount to asserting that the
supremum in \prettyref{def:Tstar} is achieved; this follows from
the discreteness of $\ssN$ and the compactness of $\TT$. Finally,
\prettyref{eq:Tstarttoinfinity} follows from \prettyref{eq:stttoinfty}.
\end{proof}

The following proposition contains the statement of \prettyref{thm:minimizers-exist}.
\begin{prop}
\label{prop:globalminimizersexist}Suppose that $\mu\in\widetilde{\Omega}_{1}$.
For each $\theta,t\in\RR$ and $x\in\TT$, the set $\scM_{t,x}^{\theta}$
is nonempty. In particular, we have
\begin{equation}
\{X|_{[T_{*}(t),t]}\st X\in\scM_{t,x}^{\theta}\}=\scM_{T_{*}(t),y_{*}(t)\mid t,x}^{\theta}.\label{eq:Mstarisrestriction}
\end{equation}
Moreover, $\scM_{t,x}^{\theta}$ has minimal and maximal elements
under the partial order $\preccurlyeq$.
\end{prop}

\begin{proof}
Since $\scM_{T_{*}(t),y_{*}(t)\mid t,x}^{\theta}$ is nonempty by
\prettyref{prop:basicproperties}(\ref{enu:straightlines}), to prove
that the set $\scM_{t,x}^{\theta}$ is nonempty it suffices to prove
\prettyref{eq:Mstarisrestriction}. The ``$\subseteq$'' direction
is an immediate consequence of \prettyref{prop:basicproperties}(\ref{enu:subpathalsominimizes})
and the definition of $\scM_{t,x}^{\theta}$, so we turn our attention
to the ``$\supseteq$'' direction. In other words, given $X_{1}\in\scM_{T_{*}(t),y_{*}(t)\mid t,x}^{\theta}$,
we seek to extend $X_{1}$ to an element $X$ of $\scM_{t,x}^{\theta}$.
Let $t_{0}=t$, and let $t_{k}=T_{*}(t_{k-1})$ and $y_{k}=y_{*}(t_{k-1})$
for $k\ge1$. We note that, since $t_{k}\le t_{k-1}-M(\mu)$, we have
\begin{equation}
\bigsqcup_{k=1}^{\infty}(t_{k},t_{k-1}]=(-\infty,t].\label{eq:tksfillup}
\end{equation}
For $k\ge2$, let $X_{k}\in\mathscr{M}_{t_{k},y_{k}\mid t_{k-1},y_{k-1}}^{\theta}$.
(This inclusion is satisfied for $k=1$ as well, but $X_{1}$ has
already been chosen.) Now for $s\in[t_{k},t_{k-1})$, define $X(s)\coloneqq X_{k}(s)$,
so $X$ is defined as an element of $H_{\mathrm{loc}}^{1}((-\infty,t];\TT)$
by \prettyref{eq:tksfillup} and the fact that $X_{k}(t_{k-1})=y_{k-1}=X_{k-1}(t_{k-1})$
by definition.

We claim that $X\in\scM_{t,x}^{\theta}$. Let $s<t$ and let $k$
be large enough that $t_{k}\le s$. Suppose that $z\in\TT$ and that
$Y\in\scM_{t_{k}-M,z\mid t,x}^{\theta}$. Then, by \prettyref{eq:allpassthrough},
we have $Y(t_{j})=y_{j}=X(t_{j})$ whenever $j\le k$. This means
that $Y|_{[t_{j},t_{j-1}]}\in\scM_{t_{k},y_{k}\mid t_{k-1},y_{k-1}}^{\theta}$.
Since the same is true for $X|_{[t_{j},t_{j-1}]}$, we have $\mathcal{A}_{\theta,t_{j},t_{j-1}}[X]=\mathcal{A}_{\theta,t_{j},t_{j-1}}[Y]$.
Summing this up over all $j$, we obtain 
\begin{equation}
\mathcal{A}_{\theta,t_{k},t}[Y]=\mathcal{A}_{\theta,t_{k},t}[X].\label{eq:YXsame}
\end{equation}
Now since $Y\in\scM_{t_{k},y_{k}\mid t,x}^{\theta}$ by \prettyref{prop:basicproperties}(\ref{enu:subpathalsominimizes}),
\prettyref{eq:YXsame} means that $X\in\scM_{t_{k},y_{k}\mid t,x}^{\theta}$
as well. But since $t_{k}\le s$, we can apply \prettyref{prop:basicproperties}(\ref{enu:subpathalsominimizes})
once again to see that $X\in\scM_{s,X(s)\mid t,x}^{\theta}$ Since
this is true for any $s<t$, we conclude that $X\in\scM_{t,x}^{\theta}$.

To show that $\mathscr{M}_{t,x}^{\theta}$ has minimal/leftmost and
maximal/rightmost elements under $\preccurlyeq$, we observe that
each $\scM_{t_{k},y_{k}\mid t_{k-1},y_{k-1}}^{\theta}$ is finite
and has leftmost and rightmost elements (as can be seen using \prettyref{prop:basicproperties}(\ref{enu:pathsurgery})
to build a path that is weakly to the left/right of any other minimizer).
Then we note that a concatenation of these leftmost and rightmost
elements in a similar manner to the above argument will yield leftmost
and rightmost elements of $\mathscr{M}_{t,x}^{\theta}$. 
\end{proof}
The following corollary emphasizes how times of the form $T_{*}(t)$
serve as regeneration times such that the behavior of minimizers before
and after them is independent.
\begin{cor}
\label{cor:justdependsonendstuff}If $\mu\in\widetilde{\Omega}_{1}$,
$\theta\in\RR$, $x\in\TT$, and $s\le t$, then
\begin{equation}
\mathscr{M}_{t,x}^{\theta}=\mathscr{M}_{T_{*}(s),y_{*}(s)}^{\theta}\odot_{T_{*}(s)}\scM_{T_{*}(s),y_{*}(s)\mid t,x}^{\theta}=\{X\odot_{T_{*}(s)}Y\st X\in\mathscr{M}_{T_{*}(s),y_{*}(s)}^{\theta}\text{ and }Y\in\scM_{T_{*}(s),y_{*}(s)\mid t,x}^{\theta}\}.\label{eq:Misconcat}
\end{equation}
In particular, if we make the shorthand definition
\begin{equation}
\mathscr{M}_{*\mid t,x}^{\theta}\coloneqq\mathscr{M}_{T_{*}(t),y_{*}(t)\mid t,x}^{\theta},\label{eq:Mstartx}
\end{equation}
then we have
\[
\mathscr{M}_{t,x}^{\theta}=\mathscr{M}_{T_{*}(s),y_{*}(s)}^{\theta}\odot_{T_{*}(s)}\scM_{*\mid t,x}^{\theta}.
\]
\end{cor}

\begin{proof}
This follows from \prettyref{prop:globalminimizersexist}, \prettyref{prop:basicproperties}(\ref{enu:pathsurgery}),
and induction.
\end{proof}
\FloatBarrier{}

\subsection{Minimizers and shocks}

Now that we have established the existence of global minimizers, the
definition \prettyref{eq:Sdef} of the shock set $\ssS$ makes sense.
Here we establish a few basic properties about minimizers and shocks
that will be useful in the sequel.
\begin{lem}
\label{lem:whencanyousplit}Suppose that $\mu\in\widetilde{\Omega}_{1}$.
Let $\theta\in\RR$, $s<t$, and $X_{1},X_{2}\in\scM_{t}^{\theta}$.
If $r\in[s,t]$ is such that $X_{1}(r)=X_{2}(r)$ but $X_{1}'(r-)\ne X_{2}'(r-)$,
then $(r,X_{i}(r))\in\ssS_{\theta}$. If we moreover assume that $r<t$,
then $(r,X_{i}(r))\in\ssN$ as well.
\end{lem}

\begin{proof}
It follows from the definitions and \prettyref{prop:basicproperties}(\ref{enu:subpathalsominimizes})
that 
\begin{equation}
X_{i}|_{(-\infty,r]}\in\scM_{r,X_{i}(r)}^{\theta}\qquad\text{for each }i\in\{1,2\},\label{eq:minimizelater}
\end{equation}
so $(\theta,r,X_{i}(r))\in\ssS$ by the definition of $\ssS$.

Now we assume that $r<t$ and prove that $(r,X_{i}(r))\in\ssN$. We
note that the restrictions of $X_{1}$ and $X_{2}$ are both elements
of $\scM_{T_{*}(r),y_{*}(r)\mid r,X_{i}(r)}^{\theta}$ by \prettyref{eq:minimizelater}
and \prettyref{eq:Mstarisrestriction}, so \prettyref{prop:basicproperties}(\ref{enu:pathsurgery})
and the assumption that $X_{1}'(r-)\ne X_{2}'(r-)$ imply that there
is an element of $\scM_{T_{*}(r),y_{*}(r)\mid t,X_{1}(t)}^{\theta}$
that changes direction at $(r,X_{1}(r))$. Hence \prettyref{prop:basicproperties}(\ref{enu:straightlines})
implies that $(r,X_{i}(r))\in\ssN$.
\end{proof}

\begin{prop}
\label{prop:allshocksmerge}Suppose that $\mu\in\widetilde{\Omega}_{1}$.
For any $\theta,t_{0}\in\RR$, there is a $t>t_{0}$ such that $\#\GS_{\theta,t}=\#\ssS_{\theta,t}=1$.
\end{prop}

\begin{proof}
Let $M$ be as in \prettyref{def:smallnoisezones}. Using \prettyref{eq:stttoinfty},
we can find an $s\ge t_{0}$ such that \prettyref{eq:intersectiononly}
and \prettyref{eq:youshouldtake} hold. By \prettyref{prop:gothroughy},
we have $X(s)=y$ for any $X\in\scM_{t}^{\theta}$. Thus all minimizers
at time $s+M$ begin with straight line segments to $(s,y)$, and
there is a single shock at the point where the line segments switch
the direction they go around the torus, as shown in 
\begin{figure}[t]
\begin{centering}
\input{figures/single_shock_picture.pgf}
\par\end{centering}
\caption[Sufficient criterion to have a single global shock at a given time]{When all minimizers at time $s+M$ start with straight line segments
to $(s,y)$, then there is a single shock, which is in fact a global
shock, at time $s+M$.\label{fig:single_shock_picture}}
\end{figure}
\prettyref{fig:single_shock_picture}.
\end{proof}

\section{Continuity of minimizers with respect to parameters}

\label{sec:Continuity-of-minimizers}We now want to explore how the
sets $\scM_{t,x}^{\theta}$ change as we vary $\theta$, $t$, and
$x$. The main result of this section is that, if $\theta,t,x$ are
perturbed only slightly, then each new minimizer uses the same forcing
points as one of the original minimizers.

For a path $X$ in $H_{\mathrm{loc}}^{1}((-\infty,t])$, $t_{0}\le t$,
$\tau\in(t_{0}-t,+\infty)$, and $\eta\in\RR$, if $X|_{[t_{0},t]}$
is linear, we define a new path $\mathcal{T}_{t_{0},\tau,\eta}X\in H_{\mathrm{loc}}^{1}((-\infty,t+\tau])$
by
\[
\mathcal{T}_{t_{0},\tau,\eta}X(s)\coloneqq\begin{cases}
X(s), & s\le t_{0};\\
X(t_{0})+(s-t_{0})\cdot\frac{(t-t_{0})X'(t-)+\eta}{t+\tau-t_{0}}, & s\in[t_{0},t+\tau].
\end{cases}
\]
This means in particular that
\[
\mathcal{T}_{t_{0},\tau,\eta}X(t+\tau)=X(t_{0})+(t-t_{0})X'(t-)+\eta=X(t)+\eta
\]
(here we use the assumption that $X|_{[t_{0},t]}$ is linear) and
\begin{equation}
(\mathcal{T}_{t_{0},\tau,\eta}X)'(t+\tau)=\frac{(t-t_{0})X'(t-)+\eta}{t+\tau-t_{0}}.\label{eq:TXderiv}
\end{equation}
We also define
\begin{equation}
t_{-;1}[X]\coloneqq\min\{r<t\st(r,X(r))\in\ssN\}\label{eq:tminus1}
\end{equation}
and
\begin{equation}
t_{-;2}[X]\coloneqq\min\{r\le t\st(r,X(r))\in\ssN\}\label{eq:tminus2}
\end{equation}
and, for $i=1,2$,
\begin{equation}
\mathcal{T}_{\tau,\eta}^{i}X=\mathcal{T}_{t_{-;i}[X],\tau,\eta}X.\label{eq:Tidef}
\end{equation}

We note that if $(t,X(t))\not\in\ssN$, then $\mathcal{T}_{\tau,\eta}^{i}X$
does not depend on $i$. If $(t,X(t))\in\ssN$, then both $\mathcal{T}_{\tau,\eta}^{1}$
and $\mathcal{T}_{\tau,\eta}^{2}$ move the endpoint of $X$ to $(t+\tau,X(t)+\eta)$,
but $\mathcal{T}_{\tau,\eta}^{2}X$ keeps the forcing point at $(t,X(t))$,
while $\mathcal{T}_{\tau,\eta}^{1}X$ skips over it. See \prettyref{fig:move_endpoint_picture}.
\begin{figure}[t]
\begin{centering}
\input{figures/move_endpoint_picture.pgf}
\par\end{centering}
\caption[Distinction between different notions of perturbing a minimizer]{The distinction between $\mathcal{T}_{\tau,\eta}^{1}X$ and $\mathcal{T}_{\tau,\eta}^{2}X$
when $(t,X(t))\in\protect\ssN$.\label{fig:move_endpoint_picture}}
\end{figure}

Now we can state our proposition.
\begin{prop}
\label{prop:main-continuity}Suppose that $\mu\in\widetilde{\Omega}_{1}$.
Fix $\theta,t\in\RR$ and $x\in\TT$. There exists an $\eps=\eps(\mu,\theta,t,x)\in(0,t-t_{-;1}[X])$
such that whenever $\zeta,\tau,\eta\in(-\eps,\eps)$, we have
\begin{equation}
\mathscr{M}_{t+\tau,x+\eta}^{\theta+\zeta}\subseteq\begin{cases}
\mathcal{T}_{\tau,\eta}^{1}\mathscr{M}_{t,x}^{\theta} & \text{if }(t,x)\not\in\ssN\text{ or }\tau\le0;\\
\mathcal{T}_{\tau,\eta}^{1}\mathscr{M}_{t,x}^{\theta}\cup\mathcal{T}_{\tau,\eta}^{2}\mathscr{M}_{t,x}^{\theta} & \text{otherwise.}
\end{cases}\label{eq:main-continuity}
\end{equation}
\end{prop}

The reason for the two cases on the right side of \prettyref{eq:main-continuity}
is that $\mathcal{T}_{\tau,\eta}^{2}X$ is not defined for $\tau\le0$,
and if $(t,x)\not\in\ssN$, then $\mathcal{T}_{\tau,\eta}^{1}=\mathcal{T}_{\tau,\eta}^{2}$
by definition.
\begin{proof}
First we choose $\eps_{0}>0$ small enough that
\begin{equation}
\left\{ (s,X(s))\st\zeta,\tau,\eta\in(-\eps_{0},\eps_{0}),X\in\mathscr{M}_{t+\tau,x+\eta}^{\theta+\zeta},s\in[t-\eps_{0},t+\eps_{0}]\right\} \cap\ssN\subseteq\{(t,x)\}.\label{eq:eps0cond}
\end{equation}
This is possible since $\ssN$ is discrete: first we choose $\eps_{0}$
small enough that any element $(s,y)\in\ssN\cap([t-\eps_{0},t+\eps_{0}]\times\TT)$
must have $s=t$, and then we can make $\eps_{0}$ even smaller if
necessary to ensure that it is not advantageous for a minimizer in
$\mathscr{M}_{t+\tau,x+\eta}^{\theta+\zeta}$ to use any forcing point
$(t,y)$ with $y\ne x$.

Now let $T_{0}=T_{*}(t-\eps_{0})$ and $y_{0}=y_{*}(t-\eps_{0})$.
By \prettyref{eq:Mstarisrestriction}, whenever $\zeta,\tau,\eta\in(-\eps_{0},\eps_{0})$,
we have $\scM_{t+\tau,x+\eta}^{\theta+\zeta}=\mathscr{M}_{T_{0},y_{0}}^{\theta+\zeta}\odot_{T_{0}}\mathscr{M}_{T_{0},y_{0}\mid t+\tau,x+\eta}^{\theta+\zeta}$,
so to complete the proof of the proposition it suffices to show that,
when $\zeta,\tau,\eta$ are sufficiently small, we have
\begin{equation}
\mathscr{M}_{T_{0},y_{0}\mid t+\tau,x+\eta}^{\theta+\zeta}\subseteq\begin{cases}
\mathcal{T}_{\tau,\eta}^{1}\mathscr{M} & \text{if }(t,x)\not\in\ssN\text{ or }\tau\le0;\\
\mathcal{T}_{\tau,\eta}^{1}\mathscr{M}\cup\mathcal{T}_{\tau,\eta}^{2}\mathscr{M} & \text{otherwise,}
\end{cases}\label{eq:goal}
\end{equation}
where we have defined 
\begin{equation}
\mathscr{M}\coloneqq\mathscr{M}_{T_{0},y_{0}\mid t,x}^{\theta}.\label{eq:Mdef}
\end{equation}

Fix $X\in\mathscr{M}$. Let $\mathscr{N}$ be the set of \emph{all}
paths $Y$ with $Y(T_{0})=y_{0}$, $Y(t)=x$, and $Y$ consisting
of straight line segments connecting a subset of the points of $\{(t,x),(T_{0},y_{0})\}\cup\mathsf{N}\cap([T_{0},t+\eps_{0}]\times\TT)$.
We see from the definitions and \prettyref{prop:basicproperties}(\ref{enu:straightlines})
that
\begin{equation}
\mathscr{M}_{T_{0},y_{0}\mid t+\tau,x+\eta}^{\theta+\zeta}\subseteq\begin{cases}
\mathcal{T}_{\tau,\eta}^{1}\mathscr{N} & \text{if }(t,x)\not\in\ssN\text{ or }\tau\le0;\\
\mathcal{T}_{\tau,\eta}^{1}\mathscr{N}\cup\mathcal{T}_{\tau,\eta}^{2}\mathscr{N} & \text{otherwise.}
\end{cases}.\label{eq:cangetallfromshifting}
\end{equation}

The discreteness of $\ssN$ implies that
\begin{equation}
\min_{Y\in\mathscr{N}\setminus\scM}\mathcal{A}_{\theta,T_{0},t}[Y]-\mathcal{A}_{\theta,T_{0},t}[X]>0.\label{eq:deltadef}
\end{equation}
Now \prettyref{eq:eps0cond} implies that, for each fixed $Y\in\mathscr{N}$,
the map 
\[
(\zeta,\tau,\eta)\mapsto\mathcal{A}_{\theta,T,t+\tau}[\mathcal{T}_{\tau,\eta}^{1}Y]
\]
is continuous on $[-\eps_{0}/2,\eps_{0}/2]^{3}$. Combining this observation
with \prettyref{eq:deltadef} and using the discreteness of $\ssN$
again, we can find an $\eps\in(0,\eps_{0}/2)$ such that, if $\zeta,\tau,\eta\in(-\eps,\eps)$,
then
\[
\min_{Y\in\mathscr{N}\setminus\scM}\mathcal{A}_{\theta,T_{0},t+\tau}[\mathcal{T}_{\tau,\eta}^{1}Y]-\mathcal{A}_{\theta,T_{0},t+\tau}[\mathcal{T}_{\tau,\eta}^{1}X]>0,
\]
and so
\begin{equation}
\{\mathcal{T}_{\tau,\eta}^{1}Y\st Y\in\mathscr{N}\setminus\mathscr{M}\}\cap\mathscr{M}_{T_{0},y_{0}\mid t+\tau,x+\eta}^{\theta+\zeta}=\varnothing.\label{eq:nootherguysT1}
\end{equation}
Combining this with \prettyref{eq:cangetallfromshifting} and recalling
that $\mathcal{T}_{\tau,\eta}^{1}Y=\mathcal{T}_{\tau,\eta}^{2}Y$
whenever $(t,x)\not\in\ssN$, we that \prettyref{eq:goal} is proved
in the case when $(t,x)\not\in\ssN$, and also when $\tau\le t$.

Thus we now have to consider the case when $(t,x)\in\ssN$ and $\tau>t$.
By \prettyref{eq:nootherguysT1}, we see that
\begin{equation}
\mathscr{M}_{T_{0},y_{0}\mid t+\tau,x+\eta}^{\theta+\zeta}\subseteq\mathcal{T}_{\tau,\eta}^{1}\mathscr{M}\cup\{\mathcal{T}_{\tau,\eta}^{2}Y\st Y\in\mathscr{N}\}.\label{eq:cangetallfromshifting-1}
\end{equation}
If $Z\in\mathscr{M}_{T_{0},y_{0}\mid t+\tau,x+\eta}^{\theta+\zeta}\cap\{\mathcal{T}_{\tau,\eta}^{2}Y\st Y\in\mathscr{N}\}$,
then in particular $Z(t)=x$, so by \prettyref{prop:basicproperties}(\ref{enu:subpathalsominimizes}),
we have $Z|_{[T_{0},t]}\in\mathscr{M}_{T_{0},y_{0}\mid t,x}^{\theta+\zeta}$.
On the other hand, by the case $\tau\le t$ considered above, we have
$\mathscr{M}_{T_{0},y_{0}\mid t,x}^{\theta+\zeta}\subseteq\mathscr{M}$.
Therefore, we have $Z=\mathcal{T}_{\tau,\eta}^{2}(Z|_{[T_{0},t]})\in\mathcal{T}_{\tau,\eta}^{2}\mathscr{M}$.
This completes the proof of \prettyref{eq:goal}.
\end{proof}
We will apply \prettyref{prop:main-continuity} many times throughout
the paper. The following application is quite simple.
\begin{cor}
\label{cor:minimizers-movex}Suppose that $\mu\in\widetilde{\Omega}_{1}$.
Fix $\theta,t\in\RR$ and $x\in\TT$. There exists an $\eps=\eps(\mu,\theta,t,x)>0$
such that whenever $\eta\in[0,\eps)$,
\begin{equation}
\mathscr{M}_{t,x+\eta}^{\theta}=\{\mathcal{T}_{0,\eta}^{1}X\st X\in\scM_{t,x}^{\theta}\text{ and }X'(t-)=X_{\theta,t,x,\Right}'(t-)\}\label{eq:Mx+eta}
\end{equation}
and
\begin{equation}
\mathscr{M}_{t,x-\eta}^{\theta}=\{\mathcal{T}_{0,-\eta}^{1}X\st X\in\scM_{t,x}^{\theta}\text{ and }X'(t-)=X_{\theta,t,x,\Left}'(t-)\}.\label{eq:Mx-eta}
\end{equation}
In particular,
\begin{equation}
X_{\theta,t,x+\eta,\Right}=\mathcal{T}_{0,\eta}^{1}X_{\theta,t,x,\Right}\label{eq:Xx+eta}
\end{equation}
and
\begin{equation}
X_{\theta,t,x-\eta,\Left}=\mathcal{T}_{0,-\eta}^{1}X_{\theta,t,x,\Left}.\label{eq:Xx-eta}
\end{equation}
\end{cor}

\begin{proof}
We prove the first statement, as the proof of the second is symmetrical,
and the third and fourth statements follow immediately from the first
and second, respectively. Let $\eps$ be as in \prettyref{prop:main-continuity}
(although we will make it even smaller shortly). By \prettyref{cor:justdependsonendstuff},
$c\coloneqq\mathcal{A}_{\theta,T_{*}(t),x}[X]$ does not depend on
$X\in\scM_{t,x}^{\theta}$. Thus, for any $\eta\in[0,\eps)$ and any
$X\in\mathscr{M}_{t,x}^{\theta}$, we can compute using the definitions
that
\begin{align*}
\mathcal{A}_{\theta,T_{*}(t),x+\eta}[\mathcal{T}_{0,\eta}^{1}X]-c & =\frac{1}{2}(t-t_{-;1}[X])\left[\left(X'(t-)-\theta+\frac{\eta}{t-t_{-;1}[X]}\right)^{2}-\left(X'(t-)-\theta\right)^{2}\right]\\
 & =\eta(X'(t-)-\theta)+\frac{\eta^{2}}{2(t-t_{-;1}[X])}.
\end{align*}
From this expression we see that, when $\eta$ is sufficiently small
and nonnegative, the action $\mathcal{A}_{\theta,T_{*}(t),x+\eta}[\mathcal{T}_{0,\eta}^{1}X]$
is minimized over $X\in\scM_{t,x}^{\theta}$ exactly when $X'(t-)=X_{\theta,t,x,\Right}'(t-)$
(as all such $X$ will have the same value of $t_{-;1}[X]$). Using
this along with \prettyref{prop:main-continuity} yields the conclusion.
\end{proof}

\section{Topology of global shocks}

\label{sec:global-shocks}In this section we explore the properties
of the global shock set $\GS$ introduced in \prettyref{def:globalshocks}.

\subsection{Basic properties}

We begin with two refinements of the definition of global shock. We
recall that a global shock is one for which the left and right minimizers,
when they merge back together, do so with different winding numbers.
We now show that the difference in winding number is always exactly
$1$.
\begin{prop}
\label{prop:exactlyoneoff}Suppose that $\mu\in\Omega$. If $(\theta,t,x)\in\GS$,
then 
\begin{equation}
\int_{T_{\vee}(\theta,t,x)}^{t}X'_{\theta,t,x,\Left}(s)\,\dif s=\int_{T_{\vee}(\theta,t,x)}^{t}X'_{\theta,t,x,\Right}(s)\,\dif s+1.\label{eq:exactlyoneoff}
\end{equation}
\end{prop}

\begin{proof}
For $\Square\in\LR$, let $I_{\Square}\coloneqq\int_{T_{\vee}(\theta,t,x)}^{t}X'_{\theta,t,x,\Square}(s)\,\dif s$.
Using the fact that
\[
X_{\theta,t,x,\Left}(T_{\vee}(\theta,t,x))=X_{\theta,t,x,\Right}(T_{\vee}(\theta,t,x))
\]
along with \prettyref{eq:shockslopesdiffernt}, we see that $I_{\Left}-I_{\Right}$
is a positive integer. On the other hand, if $I_{\Left}-I_{\Right}>1$,
then the intermediate value theorem yields an $s\in(T_{\vee}(\theta,t,x),t)$
such that
\[
\int_{s}^{t}X'_{\theta,t,x,\Left}(r)\,\dif r-\int_{s}^{t}X'_{\theta,t,x,\Right}(r)\,\dif r=1,
\]
and then we would have $X_{\theta,t,x,\Left}(s)=X_{\theta,t,x,\Right}(s)$,
contradicting the definition \prettyref{eq:Tveedef} of $T_{\vee}(\theta,t,x)$.
Thus we must have $I_{\Left}-I_{\Right}=1$, which was the claim.
\end{proof}
The next proposition says that, in order to check that a shock is
a global shock, it is sufficient to find two minimizers that have
different slopes at time $t$ and a positive winding number upon their
first reconnection---they do not have to be the left and right minimizers.
\begin{prop}
\label{prop:XLXXR}Suppose that $\mu\in\widetilde{\Omega}_{1}$. If
there are $X,Y\in\mathscr{M}_{t,x}^{\theta}$ such that $X'(t)>Y'(t)$
and, with $T_{X,Y}\coloneqq\max\{s<t\st X(s)=Y(s)\}$, we assume that
\begin{equation}
\int_{T_{X,Y}}^{t}X'(s)\,\dif s\ne\int_{T_{X,Y}}^{t}Y'(s)\,\dif s,\label{eq:XYdiscrep}
\end{equation}
then $(\theta,t,x)\in\GS$.
\end{prop}

\begin{proof}
We first note that
\begin{equation}
\int_{r}^{t}(X'-Y')(s)\,\dif s>0\text{ for all }r\in[T_{X,Y},t]\label{eq:XYdiscrep-pos}
\end{equation}
by the intermediate value theorem, the assumption that $X'(t)>Y'(t)$,
and \prettyref{eq:XYdiscrep}, since if there were a point $r\in(T_{X,Y},t]$
such that $\int_{r}^{t}(X'-Y')(s)\,\dif s=0$ then we would have $X(r)=Y(r)$
and hence $T_{X,Y}\ge r$. In particular, this means that \prettyref{eq:XYdiscrep}
can be refined to
\begin{equation}
\int_{T_{X,Y}}^{t}X'(s)\,\dif s>\int_{T_{X,Y}}^{t}Y'(s)\,\dif s\label{eq:XYdiscrep-sign}
\end{equation}
We also note that
\[
X_{\theta,t,x,\Left}'(t-)\ge X'(t-)>Y'(t-)\ge X_{\theta,t,x,\Right}'(t),
\]
so $(\theta,t,x)\in\ssS$. Furthermore,
\[
\int_{T_{X,Y}}^{t}(X_{\theta,t,x,\Left}'-X_{\theta,t,x,\Right}')(s)\,\dif s\ge\int_{T_{X,Y}}^{t}(X'-Y')(s)\,\dif s\ge1
\]
(with the last inequality by \prettyref{eq:XYdiscrep-sign} and the
fact that the last integral must be an integer), and so $T_{\vee}(\theta,t,x)\ge T_{X,Y}$.
(See \prettyref{fig:XLXXR_picture}.)
\begin{figure}[t]
\begin{centering}
\input{figures/XLXXR_picture.pgf}
\par\end{centering}
\caption[Illustration of the situation in \prettyref{prop:XLXXR}]{Illustration of the situation in \prettyref{prop:XLXXR}.\label{fig:XLXXR_picture}}
\end{figure}
 But then we must in fact have
\[
\int_{T_{\vee}(\theta,t,x)}^{t}(X_{\theta,t,x,\Left}'-X_{\theta,t,x,\Right}')(s)\,\dif s\ge\int_{T_{\vee}(\theta,t,x)}^{t}(X'-Y')(s)\,\dif s\overset{\prettyref{eq:XYdiscrep-pos}}{>}0,
\]
and the proof is complete.
\end{proof}

\subsection{Existence of global shocks}

Now we show that there must be at least one global shock for each
$\theta$ and $t$. The reason for this is topological: since all
minimizers for a given $\theta$ and $t$ must merge, there must be
some point at which the minimizers switch the direction they travel
around the torus.
\begin{prop}
\label{prop:findglobalshock}Suppose that $\mu\in\widetilde{\Omega}$.
For each $\theta,t\in\RR$, we have $\#\GS_{\theta,t}\ge1$.
\end{prop}

\begin{proof}
Fix $\theta,t\in\RR$. For $x\in\TT$ and $\Square\in\LR$, let
\[
f_{\Square}(x)=\int_{T_{*}(t)}^{t}X_{\theta,t,x,\Square}'(s)\,\dif s.
\]
Now since $y\coloneqq X_{\theta,t,x,\Square}(T_{*}(t))$ depends neither
on $x$ nor on $\Square$, we know that $f_{\Square}(x)\in x-y+\ZZ$
for each $\Square,x$.\footnote{Since $x,y\in\TT$, $x-y$ is not well-defined, but $x-y+\mathbb{Z}$
is.} This implies that there exists an $x\in\TT$ such that $f_{\Left}(x)-f_{\Right}(x)\in\mathbb{Z}\setminus\{0\}$.
Let
\begin{equation}
s_{\otimes}=\min\left\{ s\in[T_{*}(t),t]\st\begin{array}{c}
X_{\theta,t,x,\Left}(s)=X_{\theta,t,x,\Right}(s)\text{ and}\\
\int_{T_{*}(t)}^{s}(X_{\theta,t,x,\Left}'-X_{\theta,t,x,\Right}')(s)\,\dif s\in\ZZ\setminus\{0\}
\end{array}\right\} .\label{eq:s0def}
\end{equation}
Let $y_{\otimes}=X_{\theta,t,x,\Left}(s_{\otimes})=X_{\theta,t,x,\Right}(s_{\otimes})$.
Then we can see from the definitions and \prettyref{lem:whencanyousplit}
that $(s_{\otimes},y_{\otimes})\in\GS_{\theta}$. If $s_{\otimes}=t$,
then we are done, so assume that $s_{\otimes}<t$. In this case, \prettyref{lem:whencanyousplit}
further implies that $(s_{\otimes},y_{\otimes})\in\ssN$. (Thus we
have $(s_{\otimes},y_{\otimes})=(s_{\otimes}(\theta),y_{\otimes}(\theta))$
as in \prettyref{def:Omegatilde2}, justifying the notation.) On the
other hand, the definition \prettyref{eq:s0def} implies that $s_{\otimes}>T_{*}(t)$,
so in fact
\[
\eta\coloneqq\min\{\xi\ge0\st\text{there is an }X\in\scM_{t,x+\xi}^{\theta}\text{ such that }X(s_{\otimes})\ne y_{\otimes}\}
\]
exists. Define $x'=x+\eta$. By \prettyref{prop:main-continuity},
there are $X,Y\in\scM_{t,x'}^{\theta}$ such that 
\begin{equation}
X(s_{\otimes})=y_{\otimes}\ne Y(s_{\otimes}).\label{eq:Xs0Ys0}
\end{equation}
Moreover, since $(s_{\otimes},y_{\otimes})\in\GS_{\theta}$, we can
use \prettyref{prop:basicproperties}(\ref{enu:pathsurgery}) to modify
$X$ on $[T_{*}(t),s_{\otimes}]$ if necessary to ensure that
\begin{equation}
\int_{T_{*}(t)}^{t}(X'-Y')(s)\,\dif s\in\ZZ\setminus\{0\}.\label{eq:winding}
\end{equation}

Let
\[
s_{1}\coloneqq\sup\{s\in[T_{*}(t),t]\st X(s)\ne Y(s)\}
\]
and
\[
s_{2}\coloneqq\sup\{s\in[T_{*}(t),t)\st X(s)=Y(s)\}.
\]
We claim that
\begin{equation}
s_{1}=t\qquad\text{and}\qquad\int_{s_{2}}^{t}(X'-Y')(s)\,\dif s\in\ZZ\setminus\{0\}.\label{eq:claim}
\end{equation}
Indeed, if either of these conditions fail, then by \prettyref{lem:whencanyousplit}
and \prettyref{eq:winding}, we can find an $s\in[T_{*}(t),t]$ such
that $(s,X(s))=(s,Y(s))\in\ssS_{\theta}\cap\ssN$. But then \prettyref{eq:Xs0Ys0}
implies that $(s,X(s))\ne(s_{\otimes},y_{\otimes})$, and since $(s_{\otimes},y_{\otimes})\in\ssN$
this contradicts the assumption that $\mu\in\widetilde{\Omega}_{2}$
(recalling \prettyref{def:Omegatilde2}). Thus we can conclude that
\prettyref{eq:claim} holds. Using \prettyref{prop:XLXXR}, we see
that this implies that $(t,x')\in\GS_{\theta}$, and hence the proposition
is proved.
\end{proof}

\subsection{Closedness of the global shock set}

It will be helpful to know that the set $\GS$ is closed.
\begin{prop}
\label{prop:GSclosed}Suppose that $\mu\in\widetilde{\Omega}_{1}$.
The set $\GS$ is a closed subset of $\TT\times\RR\times\RR$.
\end{prop}

\begin{proof}
Suppose that we have $(\theta_{n},t_{n},x_{n})\to(\theta,t,x)$ and
$(\theta_{n},t_{n},x_{n})\in\GS$ for each $n\in\NN$. We claim that
$(\theta,t,x)\in\GS$. By \prettyref{prop:main-continuity}, we can
find $n\in\NN$, $\zeta,\tau,\eta_{\Square}\in\RR$, $i_{\Square}\in\{1,2\}$,
and $Z_{\Square}\in\mathscr{M}_{t,x}^{\theta}$ (for $\Square\in\LR$)
such that $\theta_{n}=\theta+\zeta$, $t_{n}=t+\tau$, $x_{n}=x+\eta_{\Square}$,
$\tau>t_{-;i_{\Square}}[Z_{\Square}]-t$, and
\[
X_{\theta_{n},t_{n},x_{n},\Square}=\mathcal{T}_{\tau,\eta_{\Square}}^{i_{\Square}}Z_{\Square}
\]
for $\Square\in\LR$. In fact, by taking $n$ sufficiently large,
it is further possible to ensure that $\eta_{\Left}=\eta_{\Right}\eqqcolon\eta$
(which will be very small), so we have
\[
X_{\theta_{n},t_{n},x_{n},\Square}=\mathcal{T}_{\tau,\eta}^{i_{\Square}}Z_{\Square}\overset{\prettyref{eq:Tidef}}{=}\mathcal{T}_{t_{-;i_{\Square}}[Z_{\Square}],\tau,\eta}Z_{\Square}.
\]
Now if $Z_{\Left}'(t-)=Z_{\Right}'(t-)$, then $T_{\vee}(\theta_{n},t_{n},x_{n})\ge t_{-;1}[Z_{\Left}]=t_{-;1}[Z_{\Right}]$
and it is impossible that $(\theta_{n},t_{n},x_{n})\in\GS$. Therefore,
$Z_{\Left}'(t-)\ne Z_{\Right}'(t-)$. From this we see that $T_{\vee}(\theta_{n},t_{n},x_{n})=\sup\{s<t\st Z_{\Left}(s)\ne Z_{\Right}(s)\}$
and that
\[
1=\int_{T_{\vee}(\theta_{n},t_{n},x_{n})}^{t_{n}}(X_{\theta_{n},t_{n},x_{n},\Left}'-X_{\theta_{n},t_{n},x_{n},\Right}')(s)\,\dif s=\int_{T_{\vee}(\theta_{n},t_{n},x_{n})}^{t}(Z_{\Left}'-Z_{\Right}')(s)\,\dif s
\]
which means by \prettyref{prop:XLXXR} that $(\theta,t,x)\in\GS$.
\end{proof}
\begin{rem}
The set $\ssS$ of shocks is \emph{not} a closed subset of $\TT\times\RR\times\RR$.
Indeed, each forcing point generates two shocks extending from it
forward in time (see \prettyref{prop:advanceint-N} below), but most
forcing points are not shocks. But \prettyref{prop:GSclosed} says
that a forcing point with a \emph{global} shock coming from it must
be a global shock. This is illustrated in \prettyref{fig:global_shock_split_picture}.
(In fact, the closure of $\ssS$ is $\ssS\cup(\RR\times\ssN)$, but
we do not prove this since we do not need this fact in the paper.)
\end{rem}

\subsection{Multiple global shocks}

We have shown in \prettyref{prop:findglobalshock} that there is at
least one global shock for every $\theta$ and $t$. In this section
we show that there cannot be more than two such global shocks. The
key ingredients for this are \prettyref{lem:whencanyousplit}, which
says that minimizers can only cross each other at an element of $\ssS_{\theta}\cap\ssN$,
and \prettyref{def:Omegatilde2}, which says that if $\mu\in\widetilde{\Omega}_{2}$
then there can be at most one element of $\ssS_{\theta}\cap\ssN$.
Since, in order to have multiple global shocks at the same time, there
must be some crossing of minimizers in order to satisfy the topological
conditions in \prettyref{def:globalshocks}, these conditions restrict
the structure of multiple global shocks.

Define
\[
\GS^{\Left}=\{(\theta,t,x)\in\GS\st\text{there is an }s\in(T_{\vee}(\theta,t,x),t)\text{ such that }(s,X_{\theta,t,x,\Right}(s))\in\GS_{\theta}\cap\ssN\}
\]
and similarly
\[
\GS^{\Right}=\{(\theta,t,x)\in\GS\st\text{there is an }s\in(T_{\vee}(\theta,t,x),t)\text{ such that }(s,X_{\theta,t,x,\Left}(s))\in\GS_{\theta}\cap\ssN\}.
\]
Define $\GS_{\theta}^{\Square}=\{(t,x)\st(\theta,t,x)\in\GS^{\Square}\}$
and $\GS_{\theta,t}^{\Square}=\{x\st(\theta,t,x)\in\GS^{\Square}\}$.
\begin{lem}
\label{lem:GSLnotintersect}If $\mu\in\widetilde{\Omega}$, then $\GS^{\Left}\cap\GS^{\Right}=\varnothing$.
\end{lem}

\begin{proof}
If there were a $(\theta,t,x)\in\GS^{\Left}\cap\GS^{\Right},$then
since $X_{\theta,t,x,\Left}$ and $X_{\theta,t,x,\Right}$ do not
intersect on the time interval $(T_{\vee}(\theta,t,x),t)$ by definition,
there would have to be two distinct elements of $\GS_{\theta}\cap\ssN$,
contradicting the assumption that $\mu\in\widetilde{\Omega}_{2}$.
\end{proof}
\begin{prop}
\label{prop:findsplitpoint}Suppose that $\mu\in\widetilde{\Omega}$.
Let $\theta,t\in\RR$ and $x_{1},x_{2}\in\GS_{\theta,t}$ be such
that $x_{1}\ne x_{2}$. Then there are $i_{\Left},i_{\Right}\in\{1,2\}$
such that $\{i_{\Left},i_{\Right}\}=\{1,2\}$ and $x_{i_{\Square}}\in\GS^{\Square}$
for $\Square\in\LR$. 
\end{prop}

\begin{proof}
Define $T_{i}=T_{\vee}(\theta,t,x_{i})$, and assume without loss
of generality that 
\begin{equation}
T_{1}\ge T_{2}.\label{eq:T1geT2}
\end{equation}
Let $\hat{x}_{i}\in\RR$ be such that 
\begin{equation}
\hat{x}_{1}\le\hat{x}_{2}\le\hat{x}_{1}+1\label{eq:xhatordering}
\end{equation}
and $\pi(\hat{x}_{i})=x_{i}$ for $i\in\{1,2\}$. (Recall that $\pi\colon\RR\to\TT$
is the projection map.) Now for $i\in\{1,2\}$ and $\Square\in\LR$,
define
\[
\hat{X}_{i,\Square}(s)=\hat{x}_{i}+\int_{t}^{s}X_{\theta,t,x_{i},\Square}'(r)\,\dif r+\mathbf{1}\{(i,\Square)=(1,\Left)\},
\]
so $\pi\circ\hat{X}_{i,\Square}=X_{\theta,t,x_{i},\Square}$. In particular,
we have
\begin{equation}
\hat{X}_{1,\Left}(t)-1=\hat{X}_{1,\Right}(t)=\hat{x}_{1}.\label{eq:Xhat1soffset}
\end{equation}
and
\begin{equation}
\hat{X}_{2,\Left}(t)=\hat{X}_{2,\Right}(t)=\hat{x}_{2}.\label{eq:Xhat2smatch}
\end{equation}
Using these together with \prettyref{eq:xhatordering}, we summarize
that
\begin{equation}
\underbrace{\hat{X}_{1,\Left}(t)-1=\hat{X}_{1,\Right}(t)}_{=\hat{x}_{1}}\le\underbrace{\hat{X}_{2,\Left}(t)=\hat{X}_{2,\Right}(t)}_{=\hat{x}_{2}}\le\underbrace{\hat{X}_{1,\Left}(t)}_{=\hat{x}_{1}+1}.\label{eq:initialordering}
\end{equation}
Also, \prettyref{eq:Xhat1soffset} along with \prettyref{eq:exactlyoneoff}
imply that
\begin{equation}
\hat{X}_{1,\Left}(T_{1})=\hat{X}_{1,\Right}(T_{1}),\label{eq:Xhat1satS2}
\end{equation}
while \prettyref{eq:Xhat2smatch} along with \prettyref{eq:exactlyoneoff}
imply that
\begin{equation}
\hat{X}_{2,\Left}(T_{2})=\hat{X}_{2,\Right}(T_{2})-1.\label{eq:Xhat2satS2}
\end{equation}
Finally, we have
\[
|\hat{X}_{i,\Left}(s)-\hat{X}_{i,\Right}(s)|<1\qquad\text{for all }i\in\{1,2\}\text{ and }s\in(S_{i},t)
\]
by the definition of $T_{i}$.

We claim that 
\begin{equation}
\hat{X}_{2,\Left}(T_{1})<\hat{X}_{2,\Right}(T_{1}).\label{eq:matchcontrad}
\end{equation}
Indeed, if not, then we must have 
\begin{equation}
\hat{X}_{2,\Left}(T_{1})=\hat{X}_{2,\Right}(T_{1})\label{eq:whatifmatch}
\end{equation}
by the definitions, and since we assumed in \prettyref{eq:T1geT2}
that $T_{1}\ge T_{2}$, this implies that $T_{1}=T_{2}$, and then
\prettyref{eq:whatifmatch} contradicts \prettyref{eq:Xhat2satS2}.
Therefore, we conclude that \prettyref{eq:matchcontrad} holds.

Define
\[
s_{1}\coloneqq\min\{s\in(T_{1},t)\st\hat{X}_{2,\Left}(s)\ge\hat{X}_{1,\Right}(s)\}
\]
and
\[
s_{2}\coloneqq\min\{s\in(T_{1},t)\st\hat{X}_{2,\Right}(s)\le\hat{X}_{2,\Left}(s)\}.
\]
Now \prettyref{eq:matchcontrad} and \prettyref{eq:initialordering}
along with the intermediate value theorem imply that $s_{1}\vee s_{2}>T_{1}$.
If $s_{1}>S_{1}$, then \prettyref{lem:whencanyousplit} implies that
\[
\hat{X}_{2,\Left}(s_{1})=\hat{X}_{1,\Right}(s_{1})\eqqcolon\hat{y}_{1}
\]
and $(s_{1},\pi(\hat{y}_{1}))\in\ssS_{\theta}\cap\ssN$, while if
$s_{2}>S_{1}$, then \prettyref{lem:whencanyousplit} similarly implies
that
\[
\hat{X}_{2,\Right}(s_{2})=\hat{X}_{1,\Left}(s_{2})\eqqcolon\hat{y}_{2}
\]
and $(s_{2},\pi(\hat{y}_{2}))\in\ssS_{\theta}\cap\ssN$. Using the
fact that $\mu\in\widetilde{\Omega}_{2}$, we conclude that in fact
there is an $i\in\{1,2\}$ such that $(s_{i},\pi(\hat{y}_{i}))\in\GS_{\theta}\cap\ssN$,
and the proof is complete.
\end{proof}
\begin{cor}
\emph{\label{cor:atmosttwoshocks}}Suppose that $\mu\in\widetilde{\Omega}$.
For any $\theta,t\in\RR$, we have $\#\GS_{\theta,t}\le2$. 
\end{cor}

\begin{proof}
Suppose for the sake of contradiction that we have distinct $x_{1},x_{2},x_{3}\in\GS_{\theta,t}$.
By \prettyref{prop:findsplitpoint}, we can assume without loss of
generality that $x_{1}\in\GS_{\theta,t}^{\Left}$ and $x_{2}\in\GS_{\theta,t}^{\Right}$.
But then applying \prettyref{prop:findsplitpoint} twice more we see
that $x_{3}\in\GS_{\theta,t}^{\Left}\cap\GS_{\theta,t}^{\Right}$,
contradicting \prettyref{lem:GSLnotintersect}.
\end{proof}
\prettyref{prop:findsplitpoint} and \prettyref{cor:atmosttwoshocks}
allow us to make the following definition.
\begin{defn}
\label{def:sLsRdefs}Suppose that $\mu\in\widetilde{\Omega}$. We
define functions $\sfs_{\Left},\sfs_{\Right}\colon\RR\times\RR\to\TT$
as follows. For $(\theta,t)\in\RR\times\RR$, if $\GS_{\theta,t}=\{x\}$,
then we define $\sfs_{\Left}(\theta,t)=\sfs_{\Right}(\theta,t)=x$.
If $\GS_{\theta,t}=\{x_{\Left},x_{\Right}\}$ with $x_{\Square}\in\GS_{\theta,t}^{\Square}$
for $\Square\in\LR$, then we define $\sfs_{\Square}(\theta,t)=x_{\Square}$.
\end{defn}

Then the following proposition follows immediately from the definitions
and \prettyref{prop:findsplitpoint}.
\begin{prop}
\label{prop:ifnotequalthenwherecollide}Suppose that $\mu\in\widetilde{\Omega}$.
If $\sfs_{\Left}(\theta,t)\ne\sfs_{\Right}(\theta,t)$, then $\theta\in\Theta_{\otimes}$,
$s_{\otimes}(t)>T_{\vee}(\theta,t,x)\ge T_{*}(t)$, and
\[
X_{\theta,t,\sfs_{\Left}(\theta,t),\Right}(\sfs_{\otimes}(\theta))=X_{\theta,t,\sfs_{\Right}(\theta,t),\Left}(\sfs_{\otimes}(\theta))=y_{\otimes}(\theta).
\]
\end{prop}

\subsection{Classifying minimizers by winding number}

For certain purposes, it will be helpful to have a finer-grained classification
of minimizers according to their winding number. Recall the definition
\prettyref{eq:Mstartx} of $\scM_{*\mid t,x}^{\theta}$.
\begin{defn}
\label{def:equiv-relation}Suppose that $\mu\in\widetilde{\Omega}_{1}$.
Let $\theta,t\in\RR$ and $x\in\TT$. We define an equivalence relation
$\sim$ on $\scM_{*\mid t,x}^{\theta}$ by $X\sim Y$ whenever
\[
\int_{T_{*}(t)}^{t}(X'-Y')(r)\,\dif r=0.
\]
We define $[X]$ to be the equivalence class of $X$ under $\sim$.
We also define, for $\Square\in\LR$,
\[
\mathscr{M}_{*|t,x,\Square}^{\theta}\coloneqq[X_{\theta,t,x,\Square}].
\]
The partial order $\preccurlyeq$ is defined in \prettyref{eq:partialorderdef}
as a partial order on $\mathscr{M}_{t,x}^{\theta}$. We extend its
definition to $\mathscr{M}_{*\mid t,x}^{\theta}$ in the obvious way.
For $\Square\in\LR$, we define $X_{\theta,t,x,\Square,\Left}$ and
$X_{\theta,t,x,\Square,\Right}$ to be the minimal and maximal elements,
respectively, of $\mathscr{M}_{*\mid t,x,\Square}^{\theta}$ under
$\preccurlyeq$. Finally, for $\Square,\Diamond\in\LR$, we define
\begin{equation}
u_{\theta,\Square,\Diamond}(t,x)=X_{\theta,t,x,\Square,\Diamond}'(t-).\label{eq:usquarediamonddef}
\end{equation}
\end{defn}

We note that 
\begin{equation}
X_{\theta,t,x,\Square,\Square}=X_{\theta,t,x,\Square}\qquad\text{and}\qquad u_{\theta,\Square,\Square}=u_{\theta,\Square}\qquad\text{for }\Square\in\LR.\label{eq:leftleftrightright}
\end{equation}

First we show that global shocks lead to multiple equivalence classes.
\begin{prop}
\label{prop:GS-equiv-class}Suppose that $\mu\in\widetilde{\Omega}_{1}$.
If $(\theta,t,x)\in\GS$, then $\mathscr{M}_{*\mid t,x,\Left}^{\theta}\ne\mathscr{M}_{*\mid t,x,\Right}^{\theta}$.
\end{prop}

\begin{proof}
If $(\theta,t,x)\in\GS$, then we have
\[
\int_{T_{*}(t)}^{t}(X_{\theta,t,x,\Left}'-X_{\theta,t,x,\Right}')(s)\,\dif s\ge\int_{T_{\vee}(\theta,t,x)}^{t}(X_{\theta,t,x,\Left}'-X_{\theta,t,x,\Right}')(s)\,\dif s>0,
\]
so $X_{\theta,t,x,\Left}\not\sim X_{\theta,t,x,\Right}$ and hence
$\mathscr{M}_{*\mid t,x,\Left}^{\theta}\ne\mathscr{M}_{*\mid t,x,\Right}^{\theta}$.
\end{proof}
The converse of \prettyref{prop:GS-equiv-class} is false, even if
we additionally assume that $(\theta,t,x)\in\ssS$. Indeed, if $\sfs_{\Left}(\theta,t)\ne\sfs_{\Right}(\theta,t)$,
then any minimizer that passes through $(s_{\otimes}(\theta),y_{\otimes}(\theta))$
will split into two minimizers with different winding numbers on $[T_{*}(t),t]$.
There may be additional non-global shocks at time $t$ whose minimizers
pass through $(s_{\otimes}(\theta),y_{\otimes}(\theta))$ after re-merging
for the first time. However, we do have the following:
\begin{prop}
\label{prop:globalshocknojump}Suppose that $\mu\in\widetilde{\Omega}_{1}$.
If $(\theta,t,x)\not\in\GS$, then $u_{\theta,\Left,\Square}(t,x)=u_{\theta,\Right,\Square}(t,x)$
for each $\Square\in\LR$.
\end{prop}

\begin{proof}
We prove this for $\Square=\Left$; the proof for $\Square=\Right$
is symmetrical. We abbreviate $X_{\Diamond}\coloneqq X_{\theta,t,x,\Diamond,\Left}$,
so (recalling \prettyref{eq:usquarediamonddef}) we have $u_{\theta,\Diamond,\Left}(t,x)=X_{\Diamond}'(t-)$.
Suppose for the sake of contradiction that $X_{\Left}'(t-)\ne X_{\Right}'(t-)$.
Let $T\coloneqq\max\{s<t\st X_{\Left}(s)=X_{\Right}(s)\}$ and let
$y\coloneqq X_{\Left}(T)=X_{\Right}(T)$. We note that
\begin{equation}
\int_{T}^{t}(X_{\Left}'-X_{\Right}')(s)\,\dif s=0,\label{eq:samewindingwhentheyhit}
\end{equation}
since otherwise we would have $(\theta,t,x)\in\GS$ by \prettyref{prop:XLXXR}.
Now we consider two cases (which are in fact symmetrical, but we prove
both for clarity).
\begin{casenv}
\item Suppose first that $X_{\Left}'(t-)<X_{\Right}'(t-)$. Define $Y=X_{\Left}\odot_{T}X_{\Right}$.
Then we have $Y\ne X_{\Left}$ and $Y\preccurlyeq X_{\Left}$. Also,
we have $\int_{T_{*}(t)}^{t}(Y'-X_{\Left}')(s)\,\dif s=0$ by \prettyref{eq:samewindingwhentheyhit},
so $Y\in\mathscr{M}_{*|t,x,\Left}^{\theta}$. This contradicts the
definition of $X_{\Left}$ as the leftmost element of $\mathscr{M}_{*|t,x,\Left}^{\theta}$.
\item Suppose now that $X_{\Left}'(t-)>X_{\Right}'(t-)$. Define $Y=X_{\Right}\odot_{T}X_{\Left}$.
Then we have $Y\ne X_{\Right}$ and $Y\preccurlyeq X_{\Right}$. Also,
we have $\int_{T_{*}(t)}^{t}(Y'-X_{\Right}')(s)\,\dif s=0$ by \prettyref{eq:samewindingwhentheyhit},
so $Y\in\mathscr{M}_{*|t,x,\Right}^{\theta}$. This contradicts the
definition of $X_{\Right}$ as the leftmost element of $\mathscr{M}_{*|t,x,\Right}^{\theta}$.\qedhere
\end{casenv}
\end{proof}
The following lemma will also be useful.
\begin{lem}
\label{lem:onlytwoguys}Suppose that $\mu\in\widetilde{\Omega}_{1}$.
Let $\theta,t\in\RR$ and $x\in\TT$. If
\begin{equation}
\int_{T_{*}(t)}^{t}[X_{\theta,t,x,\Left}'(s)-X_{\theta,t,x,\Right}'(s)]\,\dif s\ge2,\label{eq:twowraps}
\end{equation}
then there is some $(s_{\otimes},y_{\otimes})\in[(T_{*}(t),t)\times\TT]\cap\GS_{\theta}\cap\ssN$
such that $X_{\theta,t,x,\Left}(s_{\otimes})=X_{\theta,t,x,\Right}(s_{\otimes})=y_{\otimes}$.
\end{lem}

\begin{proof}
For \prettyref{eq:twowraps} to hold, the minimizers $X_{\theta,t,x,\Left}$
and $X_{\theta,t,x,\Right}$ must cross in the interval $(T_{*}(t),t)$,
so the conclusion is a consequence of \prettyref{lem:whencanyousplit}.
\end{proof}
We now use the equivalence relation introduced in \prettyref{def:equiv-relation}
to give a further characterization of the left and right global shocks,
when they differ.
\begin{prop}
\label{prop:findleftrightshocks}Suppose that $\mu\in\widetilde{\Omega}$
and that $(\theta,t,x)\in\GS$. If $u_{\theta,\Left,\Right}(t,x)=u_{\theta,\Right,\Right}(t,x)$,
then 
\begin{equation}
\sfs_{\Left}(\theta,t)=x\ne\sfs_{\Right}(\theta,t).\label{eq:shocksdifferent-xleft}
\end{equation}
Similarly, if $u_{\theta,\Left,\Left}(t,x)=u_{\theta,\Right,\Left}(t,x)$,
then
\begin{equation}
\sfs_{\Left}(\theta,t)\ne x=\sfs_{\Right}(\theta,t).\label{eq:shocksdifferent-xright}
\end{equation}
\end{prop}

\begin{proof}
We will only prove the first statement, as the second is symmetrical.
So assume that $u_{\theta,\Left,\Right}(t,x)=u_{\theta,\Right,\Right}(t,x)$.
Let $s_{\otimes}=\sup\{r\le t\st X_{\theta,t,x,\Left,\Right}(r)\ne X_{\theta,t,x,\Right,\Right}(r)\}$.
Since $u_{\theta,\Left,\Right}(t,x)=u_{\theta,\Right,\Right}(t,x)$,
we have $s_{\otimes}<t$. Also, we cannot have $X_{\theta,t,x,\Left,\Right}=X_{\theta,t,x,\Right,\Right}$
since these paths are in different equivalence classes of $\sim$,
so we must have $s_{\otimes}>-\infty$. Define
\begin{equation}
y_{\otimes}=X_{\theta,t,x,\Right}(s_{\otimes}).\label{eq:y0def}
\end{equation}
By \prettyref{lem:whencanyousplit}, we have
\begin{equation}
(s_{\otimes},y_{\otimes})\in\ssS_{\theta}\cap\ssN.\label{eq:s0y0inSN}
\end{equation}

We claim that 
\begin{equation}
s_{\otimes}>T_{\vee}(\theta,t,x).\label{eq:s0gt}
\end{equation}
Suppose for the sake of contradiction that $s_{\otimes}\le T_{\vee}(\theta,t,x)$.
Define the path
\[
Y(s)\coloneqq\left(X_{\theta,t,x,\Left,\Right}\odot_{T_{\vee}(\theta,t,x)}X_{\theta,t,x,\Left}\right)(s)=\begin{cases}
X_{\theta,t,x,\Left}(s) & \text{if }s\ge T_{\vee}(\theta,t,x);\\
X_{\theta,t,x,\Left,\Right}(s), & \text{if }s\in[T_{*}(t),T_{\vee}(\theta,t,x)].
\end{cases}
\]
This path is continuous since $X_{\theta,t,x,\Left}(T_{\vee}(\theta,t,x))=X_{\theta,t,x,\Right}(T_{\vee}(\theta,t,x))=X_{\theta,t,x,\Left,\Right}(T_{\vee}(\theta,t,x))$.
The first identity by the definition \prettyref{eq:Tveedef} of $T_{\vee}$
and the second identity is because, since $s_{\otimes}\le T_{\vee}(\theta,t,x)$,
we see that, whenever $s\ge T_{\vee}(\theta,t,x)$, we have $X_{\theta,t,x,\Left,\Right}(s)=X_{\theta,t,x,\Right}(s)$.
This last observation moreover allows us to compute
\begin{align*}
\int_{T_{*}(t)}^{t}(X_{\theta,t,x,\Left}'-X_{\theta,t,x,\Left,\Right}')(s)\,\dif s & \ge\int_{T_{*}(t)}^{t}(Y'-X_{\theta,t,x,\Left,\Right}')(s)\,\dif s\\
 & =\int_{T_{\vee}(\theta,t,x)}^{t}(X'_{\theta,t,x,\Left}-X'_{\theta,t,x,\Right})(s)\,\dif s\overset{\prettyref{eq:shockslopesdiffernt}}{>}0,
\end{align*}
but this contradicts the definition of $X_{\theta,t,x,\Left,\Right}$.
Thus we conclude that \prettyref{eq:s0gt} holds.

Now \prettyref{eq:y0def} and \prettyref{eq:s0gt} imply that $X_{\theta,t,x,\Left}(s_{\otimes})\ne y_{\otimes}$.
In fact, using \prettyref{lem:whencanyousplit}, \prettyref{eq:s0y0inSN},
and the fact that $\mu\in\widetilde{\Omega}_{2}$, we see that
\begin{equation}
X(s_{\otimes})\ne y_{\otimes}\qquad\text{for all }X\in\mathscr{M}_{t,x}^{\theta}\text{ with }X'(t-)=X_{\theta,t,x,\Left}'(t-).\label{eq:noty0forabit}
\end{equation}
Therefore, if we define $\xi\coloneqq\inf\{\eta\ge0\st X_{\theta,t,x+\eta,\Right}(s_{\otimes})\ne y_{\otimes}\}$,
then $\prettyref{eq:y0def}$ and \prettyref{eq:Xx+eta} imply that
$\xi>0$, while \prettyref{eq:noty0forabit} and \prettyref{eq:Mx-eta}
imply that $\xi<1$. Therefore, 
\begin{equation}
x'\coloneqq x+\xi\ne x.\label{eq:xprimedef}
\end{equation}
Similar logic implies that
\begin{equation}
X_{\theta,t,x',\Left}(s_{\otimes})=y_{\otimes}\ne X_{\theta,t,x',\Right}(s_{\otimes})\label{eq:lefthitsy0}
\end{equation}
which means (again using \prettyref{lem:whencanyousplit} and the
assumption that $\mu\in\widetilde{\Omega}_{2}$) that
\begin{equation}
T_{\vee}(\theta,t,x')<s_{\otimes}.\label{eq:dontintersectbys0}
\end{equation}
Now \prettyref{eq:lefthitsy0} implies that
\begin{equation}
X_{\theta,t,x',\Left}(s)=X_{\theta,s_{\otimes},y_{\otimes},\Left}(s)\qquad\text{for all }s\in[T_{*}(t),s_{\otimes}].\label{eq:follows0y0}
\end{equation}
Also, again using \prettyref{lem:whencanyousplit}, in light of the
fact that the path $X_{\theta,t,x',\Right}$ does not cross any element
of $\ssS_{\theta}\cap\ssN$ by \prettyref{eq:lefthitsy0} and the
assumption that $\mu\in\widetilde{\Omega}_{2}$, we see that
\[
X_{\theta,t,x',\Right}(T_{\vee}(\theta,t,x))=X_{\theta,s_{\otimes},y_{\otimes},\Right}(T_{\vee}(\theta,t,x))
\]
and thus that
\[
X_{\theta,t,x',\Right}(s)=X_{\theta,s_{\otimes},y_{\otimes},\Right}(s)\qquad\text{for all }s\in[T_{*}(t),T_{\vee}(\theta,t,x)].
\]
Combined with \prettyref{eq:dontintersectbys0}, this implies that
\[
T_{\vee}(\theta,t,x')=T_{\vee}(\theta,s_{\otimes},y_{\otimes}).
\]
See \prettyref{fig:global_shock_minimizers_picture}
\begin{figure}[t]
\begin{centering}
\input{figures/global_shock_minimizers_picture.pgf}
\par\end{centering}
\caption[Illustration of the situation in the proof of \prettyref{prop:findleftrightshocks}]{Illustration of the situation in the proof of \prettyref{prop:findleftrightshocks}.\label{fig:global_shock_minimizers_picture}}
\end{figure}
 for an illustration.

Now if we define
\begin{equation}
Z\coloneqq X_{\theta,s_{\otimes},y_{\otimes},\Right}\odot_{s_{\otimes}}X_{\theta,t,x,\Left},\label{eq:Zdef}
\end{equation}
then $Z\preccurlyeq X_{\theta,t,x,\Right}$ and so
\begin{align*}
\int_{T_{\vee}(\theta,t,x')}^{t}(X_{\theta,t,x',\Left}'-X_{\theta,t,x',\Right}')(r)\,\dif r & \ge\int_{T_{\vee}(\theta,t,x')}^{t}(X_{\theta,t,x',\Left}'-Z')(r)\,\dif r\\
 & =\int_{T_{\vee}(\theta,t,x')}^{s_{\otimes}}(X_{\theta,s_{\otimes},y_{\otimes},\Left}'-X_{\theta,s_{\otimes},y_{\otimes},\Right}')(r)\,\dif r=1,
\end{align*}
with the first identity by \prettyref{eq:follows0y0} and \prettyref{eq:Zdef}
and the second by \prettyref{prop:exactlyoneoff}. This implies that
$x'\in\GS_{\theta,t}$. But since \prettyref{eq:follows0y0} implies
in particular that $X_{\theta,t,x',\Left}(s_{\otimes})=y_{\otimes}$,
and $(s_{\otimes},y_{\otimes})\in\ssS_{\theta}\cap\ssN$ and hence
in $\GS_{\theta}\cap\ssN$ by \prettyref{def:Omegatilde2}, we actually
have $x'\in\GS_{\theta,t}^{\Right}$, and hence that $x\in\GS_{\theta,t}^{\Left}$
by \prettyref{prop:findsplitpoint} and \prettyref{lem:GSLnotintersect}.
Therefore, we have $\sfs_{\Left}(\theta,t)=x\ne x'=\sfs_{\Right}(\theta,t)$
by \prettyref{def:sLsRdefs} and \prettyref{eq:xprimedef}. This completes
the proof.
\end{proof}
\begin{prop}
\label{prop:converse}Suppose that $\mu\in\widetilde{\Omega}$. Suppose
that $\theta,t\in\RR$ are such that $\sfs_{\Left}(\theta,t)\ne\sfs_{\Right}(\theta,t)$.
Then we have 
\[
u_{\theta,\Left,\Right}(t,\sfs_{\Left}(\theta,t))=u_{\theta,\Right,\Right}(t,\sfs_{\Left}(\theta,t))
\]
and similarly 
\[
u_{\theta,\Left,\Left}(t,\sfs_{\Right}(\theta,t))=u_{\theta,\Right,\Left}(t,\sfs_{\Right}(\theta,t)).
\]
\end{prop}

\begin{proof}
Again, we just prove the first statement, relying on symmetry to prove
the second. It is trivial that $u_{\theta,\Left,\Right}(t,\sfs_{\Left}(\theta,t))\ge u_{\theta,\Right,\Right}(t,\sfs_{\Left}(\theta,t))$,
so it suffices to prove that
\begin{equation}
u_{\theta,\Left,\Right}(t,\sfs_{\Left}(\theta,t))\le u_{\theta,\Right,\Right}(t,\sfs_{\Left}(\theta,t)).\label{eq:goaluLRuRR}
\end{equation}
Because of the assumption $\sfs_{\Left}(\theta,t)\ne\sfs_{\Right}(\theta,t)$,
we know from \prettyref{prop:findsplitpoint} that there is some $(s_{\otimes},y_{\otimes})\in\GS_{\theta}\cap\ssN$
such that 
\[
y_{\otimes}=X_{\theta,t,\sfs_{\Left}(\theta,t),\Right}(s_{\otimes})=X_{\theta,t,\sfs_{\Right}(\theta,t),\Left}(s_{\otimes}).
\]
Now define, for $\Square\in\LR$,
\[
Y_{\Square}\coloneqq X_{\theta,s_{\otimes},y_{\otimes},\Square}\odot_{s_{\otimes}}X_{\theta,t,\sfs_{\Left}(\theta,t),\Right}.
\]
Then we see that $Y_{\Square}\in\scM_{t,\sfs_{\Left}(\theta,t)}^{\theta}$
by \prettyref{prop:basicproperties}(\ref{enu:pathsurgery}). We also
note that 
\begin{equation}
Y_{\Right}=X_{\theta,t,\sfs_{\Left}(\theta,t),\Right},\label{eq:YRisXR}
\end{equation}
while
\begin{align}
\int_{T_{*}(t)}^{t}[Y_{\Left}'(s)-Y_{\Right}'(s)]\,\dif s & =\int_{T_{*}(t)}^{s_{\otimes}}[X'_{\theta,s_{\otimes},y_{\otimes},\Left}(s)-X'_{\theta,s_{\otimes},y_{\otimes},\Left}(s)]\,\dif s\nonumber \\
 & \ge\int_{T_{\vee}(\theta,s_{\otimes},y_{\otimes})}^{s_{\otimes}}[X'_{\theta,s_{\otimes},y_{\otimes},\Left}(s)-X'_{\theta,s_{\otimes},y_{\otimes},\Left}(s)]\,\dif s=1\label{eq:diffgtr1}
\end{align}
since $(s_{\otimes},y_{\otimes})\in\GS_{\theta}$. On the other hand,
since $X_{\theta,t,\sfs_{\Left}(\theta,t),\Left}(s_{\otimes})\ne y_{\otimes}$
by \prettyref{lem:GSLnotintersect}, \prettyref{lem:onlytwoguys}
tells us that
\[
\int_{T_{*}(t)}^{t}[X_{\theta,t,\sfs_{\Left}(\theta,t),\Left}'-X_{\theta,t,\sfs_{\Left}(\theta,t),\Right}'](s)\,\dif s=1,
\]
so \prettyref{eq:YRisXR} and \prettyref{eq:diffgtr1} tell us that
$Y_{\Left}\in\mathscr{M}_{*\mid t,\sfs_{\Left}(\theta,t),\Left}^{\theta}$.
Therefore, $X_{\theta,t,\sfs_{\Left}(\theta,t),\Left,\Right}$ must
lie to the right of $Y_{\Left}$, so 
\[
u_{\theta,\Left,\Right}(t,\sfs_{\Left}(\theta,t))=X_{\theta,t,\sfs_{\Left}(\theta,t),\Left,\Right}'(t-)\le Y_{\Left}'(t-)=X_{\theta,t,\sfs_{\Left}(\theta,t),\Right}'(t)=u_{\theta,\Right,\Right}(\theta,\sfs_{\Left}(\theta,t)),
\]
and so \prettyref{eq:goaluLRuRR} is proved.
\end{proof}
\begin{prop}
\label{prop:findshocktostay}Let $\mu\in\widetilde{\Omega}$ and $\theta,t\in\RR$.
If $\sfs_{\Left}(\theta,t)\ne\sfs_{\Right}(\theta,t)$, then we have
\begin{equation}
u_{\theta,\Right,\Left}(t,\sfs_{\Right}(\theta,t))>u_{\theta,\Right,\Right}(t,\sfs_{\Right}(\theta,t))\label{eq:rightdiff}
\end{equation}
 and
\begin{equation}
u_{\theta,\Left,\Left}(t,\sfs_{\Left}(\theta,t))>u_{\theta,\Left,\Right}(t,\sfs_{\Left}(\theta,t)).\label{eq:leftdiff}
\end{equation}
\end{prop}

\begin{proof}
We limit ourselves to proving \prettyref{eq:rightdiff}, as the proof
of \prettyref{eq:leftdiff} is symmetrical. By \prettyref{prop:ifnotequalthenwherecollide},
we have $\theta\in\Theta_{\otimes}$, 
\begin{equation}
s_{\otimes}(t)>T_{\vee}(\theta,t,\sfs_{\Right}(\theta,t)),\label{eq:mergeabove}
\end{equation}
and
\[
X_{\theta,t,\sfs_{\Right}(\theta,t),\Left}(s_{\otimes}(\theta))=y_{\otimes}(\theta).
\]
Define 
\begin{equation}
\tilde{X}(s)=X_{\theta,s_{\otimes}(\theta),y_{\otimes}(\theta),\Right}\odot_{s_{\otimes}(\theta)}X_{\theta,t,\sfs_{\Right}(\theta,t),\Left},\label{eq:Xtildedef}
\end{equation}
so $\tilde{X}\in\mathscr{M}_{t,x}^{\theta}$ by \prettyref{prop:basicproperties}(\ref{enu:pathsurgery}).
Suppose for the sake of contradiction that 
\[
\int_{T_{*}(t)}^{t}(\tilde{X}'-X_{\theta,t,\sfs_{\Right}(\theta,t),\Right}')(s)\,\dif s>0.
\]
Then we would have
\[
\int_{T_{*}(t)}^{t}(X'_{\theta,t,\sfs_{\Right}(\theta,t),\Left}-X_{\theta,t,\sfs_{\Right}(\theta,t),\Right}')(s)\,\dif s>1,
\]
which implies by \prettyref{lem:onlytwoguys} that $X_{\theta,t,\sfs_{\Right}(\theta,t),\Right}(s_{\otimes}(\theta))=y_{\otimes}(\theta)$,
contradicting \prettyref{eq:mergeabove}. Therefore, we must have
$\int_{T_{*}(t)}^{t}(\tilde{X}'-X_{\theta,t,\sfs_{\Right}(\theta,t),\Right}')(s)\,\dif s=0$
(as it is clearly nonnegative), and hence $\tilde{X}\in\mathscr{M}_{*|t,x,\Right}^{\theta}$.
But this means that 
\begin{align*}
u_{\theta,\Right,\Left}(t,\sfs_{\Right}(\theta,t))=X_{\theta,t,\sfs_{\Right}(\theta,t),\Right,\Left}'(t-)\ge\tilde{X}'(t-)\overset{\prettyref{eq:Xtildedef}}{=}X_{\theta,t,\sfs_{\Right}(\theta,t),\Left}'(t-) & >X_{\theta,t,\sfs_{\Right}(\theta,t),\Right}'(t-)\\
 & =u_{\theta,\Right,\Right}(t,\sfs_{\Right}(\theta,t)),
\end{align*}
and the proof is complete.
\end{proof}

\section{Movement of shocks as time is varied}

\label{sec:Movement-of-shocks}In this section we study how shocks,
and in particular in \prettyref{subsec:movement-global} global shocks,
move as $t$ changes.

\subsection{General case}

\begin{figure}[t]
\begin{centering}
\subfloat[Three minimizers from a single shock at time $t$.\label{fig:shockmoveintime_picture0}]{\input{figures/shockmoveintime_picture0.pgf}

}\hfill{}\subfloat[\prettyref{prop:advanceint}: just after $t$, there is a single shock.
The minimizers follow the left and right time-$t$ minimizers. \label{fig:global_shock_minimizers_picture-1-2}]{\input{figures/shockmoveintime_picture+1.pgf}

}\hfill{}\subfloat[\prettyref{prop:backwards}: just before $t$, there is one shock
for each pair of time-$t$ minimizers with \textquotedblleft adjacent\textquotedblright{}
slopes.\label{fig:shockmoveintime_picture1}]{\input{figures/shockmoveintime_picture-1.pgf}

}
\par\end{centering}
\caption{Dynamics of a shock as time changes when there are more than two minimizers
coming from a single shock.}
\end{figure}
\prettyref{prop:main-continuity} tells us that, if $(t,x)\in\ssS_{\theta}$
and we perturb $t$ and $x$ by $\tau$ and $\eta$, respectively,
then the minimizers will be perturbations of a subset of the original
minimizers starting at $t$ and $x$. Thus, we can seek another shock
at $(t+\tau,x+\eta)$ by trying to solve for $\eta$ such that there
are multiple distinct slopes of minimizers starting at $(t+\tau,x+\eta)$.
If there are only two minimizers starting at $(t,x)$, then this procedure
is relatively straightforward since in that case those two minimizers
must be the ones that we perturb to find minimizers at $(t+\tau,x+\eta)$.
In this case we simply recover the usual Rankine--Hugoniot condition
(\prettyref{eq:RHcondition} below). However, it could also be the
case that there are more than two minimizers meeting at the same shock.
(See \prettyref{fig:shockmoveintime_picture0} for an example.) In
this case, when time is moved slightly forward, only the greatest
and least slopes of minimizers at time $t$ persist, as illustrated
in \prettyref{fig:global_shock_minimizers_picture-1-2} and proved
in \prettyref{prop:advanceint}. On the other hand, when time is moved
slightly backward, if the slopes of minimizers are ordered, then every
pair of adjacent slopes corresponds to a separate shock, as illustrated
in \prettyref{fig:shockmoveintime_picture1} and proved in \prettyref{prop:backwards}.

We begin by considering the perturbative theory for single pairs of
minimizers. Recall the definition \prettyref{eq:Tidef} of the path
perturbation operator $\mathcal{T}_{\tau,\eta}^{1}$.
\begin{lem}
\label{lem:differentiate-A}Suppose that $\mu\in\widetilde{\Omega}_{1}$.
Suppose that $\theta\in\RR$ and $(t,x)\in\ssS_{\theta}\setminus\ssN$.
If $X_{1},X_{2}\in\scM_{t,x}^{\theta}$ are such that
\begin{equation}
X_{1}'(t-)\ne X_{2}'(t-),\label{eq:derivsnotthesame}
\end{equation}
then there is an $\eps>0$ and a unique function $\sfr_{X_{1},X_{2}}\colon(-\eps,\eps)\to\RR$
such that $\sfr_{X_{1},X_{2}}(0)$ and, for all $\tau\in(-\eps,\eps)$,
we have, defining $T_{*}=T_{*}(t-\eps)$, that
\begin{equation}
\mathcal{A}_{\theta,T_{*},t+\tau}[\mathcal{T}_{\tau,\sfr_{X_{1},X_{2}}(\tau)}^{1}X_{i}]\text{ is equal for }i\in\{1,2\}.\label{eq:Aequal}
\end{equation}
Moreover, we have
\begin{equation}
\sfr_{X_{1},X_{2}}'(0)=\frac{1}{2}(X_{1}'(t-)+X_{2}'(t-)).\label{eq:RHcondition}
\end{equation}
\end{lem}

\begin{proof}
For any piecewise-linear path $X\colon[T_{*},t]\to\TT$, we first
compute from \prettyref{eq:TXderiv} that
\begin{align}
(\mathcal{T}_{\tau,\eta}^{1}X)'(t+\tau-)-X'(t-) & =\frac{\eta-\tau X'(t-)}{t+\tau-t_{-;1}[X]}.\label{eq:derivdiff}
\end{align}
Thus we can compute, using the definition \prettyref{eq:Adef} as
well as \prettyref{eq:derivdiff}, that
\begin{equation}
\mathcal{A}_{\theta,T_{*},t+\tau}[\mathcal{T}_{\tau,\eta}^{1}X]=\frac{1}{2}(t+\tau-t_{-;1}[X])\left(X'(t-)+\frac{\eta-\tau X'(t-)}{t+\tau-t_{-;1}[X]}-\theta\right)^{2}+\mathcal{A}_{\theta,T_{*},t_{-;1}[X]}[X].\label{eq:AX}
\end{equation}
Using this with $X=X_{1},X_{2}$ and differentiating with respect
to $\eta$, we see that 
\begin{equation}
\frac{\dif}{\dif\eta}\left[\mathcal{A}_{\theta,T_{*},t+\tau}[\mathcal{T}_{\tau,\eta}^{1}X_{2}]-\mathcal{A}_{\theta,T_{*},t+\tau}[\mathcal{T}_{\tau,\eta}^{1}X_{1}]\right]_{\eta=\tau=0}=X_{2}'(t-)-X_{1}'(t-).\label{eq:etaderiv}
\end{equation}
Now the assumption \prettyref{eq:derivsnotthesame} let us use the
implicit function theorem to find an $\eps>0$ and a unique function
$\sfr_{X_{1},X_{2}}\colon(-\eps,\eps)\to\RR$ such that $\sfr_{X_{1},X_{2}}(0)=0$
and \prettyref{eq:Aequal} holds for $\tau\in(-\eps,\eps)$. We can
also compute from \prettyref{eq:AX} that 
\begin{equation}
\frac{\dif}{\dif\tau}\left[\mathcal{A}_{\theta,T_{*},t+\tau}[\mathcal{T}_{\tau,\eta}^{1}X_{2}]-\mathcal{A}_{\theta,T_{*},t+\tau}[\mathcal{T}_{\tau,\eta}^{1}X_{1}]\right]_{\eta=\tau=0}=-\frac{1}{2}X_{2}'(t-)^{2}+\frac{1}{2}X_{1}'(t)^{2},\label{eq:tauderiv}
\end{equation}
so \prettyref{eq:RHcondition} follows by implicit differentiation
along with \prettyref{eq:Aequal}, \prettyref{eq:etaderiv}, and \prettyref{eq:tauderiv}.
\end{proof}
\begin{rem}
Generalizing the computations in the last proof, we can compute, for
\emph{any} $\sff\colon(-\eps,\eps)\to\RR$, that
\begin{align}
\frac{\dif}{\dif\tau} & \left[\mathcal{A}_{\theta,T_{*},t+\tau}[\mathcal{T}_{\tau,\sff(\tau)}^{1}X_{2}]-\mathcal{A}_{\theta,T_{*},t+\tau}[\mathcal{T}_{\tau,\sff(\tau)}^{1}X_{1}]\right]_{\tau=0}\nonumber \\
 & =-\frac{1}{2}X_{2}'(t-)^{2}+\frac{1}{2}X_{1}'(t-)^{2}+\sff'(0)(X_{2}'(t-)-X_{1}'(t-))\nonumber \\
 & =\left(\sff'(0)-\frac{X_{2}'(t)+X_{1}'(t)}{2}\right)(X_{2}'(t)-X_{1}'(t)).\label{eq:generalfderiv}
\end{align}
\end{rem}

Now we will seek to identify the actual minimizers that occur after
a perturbation. First we see what happens immediately after time $t$.
The following proposition justifies the picture in \prettyref{fig:global_shock_minimizers_picture-1-2}.
\begin{prop}
\label{prop:advanceint}Suppose that $\mu\in\widetilde{\Omega}_{1}$.
Suppose that $\theta\in\RR$ and $(t,x)\in\ssS_{\theta}\setminus\ssN$.
There is an $\eps>0$ and a function $\sfr\colon[0,\eps)\to\RR$ such
that $\sfr(0)=0$, 
\begin{equation}
\sfr'(0+)=\frac{1}{2}(X_{\theta,t,x,\Left}'(t-)+X_{\theta,t,x,\Right}'(t-)),\label{eq:advanceint-rderiv}
\end{equation}
and, for all $\tau\in(0,\eps)$, we have
\begin{equation}
\mathscr{M}_{t+\tau,x+\sfr(\tau)}^{\theta}=\{\mathcal{T}_{\tau,\sfr(\tau)}^{1}X\st X\in\scM_{t,x}^{\theta}\text{ and }X'(t-)\in\{X_{\theta,t,x,\Square}'(t-)\st\Square\in\LR\}\}.\label{eq:usesameminimum}
\end{equation}
In particular, $(t+\tau,x+\sfr(\tau))\in\ssS_{\theta}$. Moreover,
if $\eta\in(-\eps,\eps)\setminus\{\sfr(\tau)\}$, then $(t+\tau,x+\eta)\not\in\ssS_{\theta}$.
\end{prop}

\begin{proof}
We choose $\eps>0$, $T_{*}=T_{*}(t-\eps)$, and $\sfr\coloneqq\sfr_{X_{\theta,t,x,\Left},X_{\theta,t,x,\Right}}$
as defined in \prettyref{lem:differentiate-A}, so \prettyref{eq:advanceint-rderiv}
holds by \prettyref{eq:RHcondition}, and we have
\begin{equation}
\mathcal{A}_{\theta,T_{*},t+\tau}[\mathcal{T}_{\tau,\sfr(\tau)}^{1}X_{\theta,t,x,\Left}]\overset{\prettyref{eq:Aequal}}{=}\mathcal{A}_{\theta,T_{*},t+\tau}[\mathcal{T}_{\tau,\sfr(\tau)}^{1}X_{\theta,t,x,\Right}]\qquad\text{for all }\tau\in(-\eps,\eps).\label{eq:Aequal-apply}
\end{equation}
If $X\in\scM_{t,x}^{\theta}$ is such that
\begin{equation}
X_{\theta,t,x,\Left}'(t-)>X'(t-)>X_{\theta,t,x,\Right}'(t-),\label{eq:intermediate}
\end{equation}
then we have
\begin{equation}
\begin{aligned}\frac{\dif}{\dif\tau} & \left[\mathcal{A}_{\theta,T_{*},t+\tau}[\mathcal{T}_{\tau,\sfr(\tau)}^{1}X]-\mathcal{A}_{\theta,T_{*},t+\tau}[\mathcal{T}_{\tau,\sfr(\tau)}^{1}X_{\theta,t,x,\Square}]\right]_{\tau=0}\\
\oversetx{\prettyref{eq:generalfderiv}} & =\left(\sfr'(0)-\frac{X'(t-)+X_{\theta,t,x,\Square}'(t)}{2}\right)(X'(t)-X_{\theta,t,x,\Square}'(t)).
\end{aligned}
\label{eq:derivA}
\end{equation}
Now, since $(t,x)\in\ssS_{\theta}$, we have $X_{\theta,t,x,\Left}'(t-)>X_{\theta,t,x,\Right}'(t-)$,
and hence must either have
\[
\sfr'(0)<\frac{X'(t-)+X_{\theta,t,x,\Left}'(t-)}{2}\qquad\text{or}\qquad\sfr'(0)>\frac{X'(t-)+X_{\theta,t,x,\Right}'(t-)}{2}.
\]
In the first case, we see (also using \prettyref{eq:intermediate})
that the right side of \prettyref{eq:derivA} is strictly positive
when $\Square=\Left$, and in the second case the right side of \prettyref{eq:derivA}
is strictly positive when $\Square=\Right$. From this and \prettyref{prop:main-continuity},
we see that the ``$\subseteq$'' direction of \prettyref{eq:usesameminimum}
holds.

On the other hand, if $X'(t)\in\{X_{\theta,t,x,\Square}'(t-)\st\Square\in\LR\}$,
then it is clear from the definitions and \prettyref{eq:Aequal-apply}
that
\begin{equation}
\mathcal{A}_{\theta,T_{*},t+\tau}[\mathcal{T}_{\tau,\sfr(\tau)}^{1}X]=\mathcal{A}_{\theta,T_{*},t+\tau}[\mathcal{T}_{\tau,\sfr(\tau)}^{1}X_{\theta,t,x,\Square}]\qquad\text{for all }\tau\in(-\eps,\eps)\text{ and }\Square\in\LR.\label{eq:sidessame}
\end{equation}
Using this observation along with \prettyref{prop:main-continuity},
we conclude the equality in \prettyref{eq:usesameminimum}. The last
assertion of the proposition then follows from the uniqueness assertion
in \prettyref{lem:differentiate-A}.
\end{proof}
We also consider the case when $(t,x)\in\ssN$.
\begin{prop}
\label{prop:advanceint-N}Suppose that $\mu\in\widetilde{\Omega}_{1}$.
Suppose that $\theta\in\RR$ and $(t,x)\in\ssN$. There is an $\eps>0$
and two continuous functions $\sfr_{\Left},\sfr_{\Right}\colon[0,\eps)\to\RR$
such that $\sfr_{\Square}(0)=0$ and for all $\tau\in(0,\eps)$, we
have $\sfr_{\Left}(\tau)<\sfr_{\Right}(\tau)$ and
\begin{equation}
\begin{aligned}\mathscr{M}_{t+\tau,x+\sfr_{\Square}(\tau)}^{\theta} & =\left\{ \mathcal{T}_{\tau,\sfr_{\Square}(\tau)}^{1}X\st X\in\scM_{t,x}^{\theta}\text{ and }X'(t-)=X_{\theta,t,x,\Square}'(t-)\right\} \\
 & \qquad\cup\left\{ \mathcal{T}_{\tau,\sfr_{\Square}(\tau)}^{2}X\st X\in\scM_{t,x}^{\theta}\right\} .
\end{aligned}
\label{eq:usesameminimums-split}
\end{equation}
In particular, this implies that if $(t,x)\in\GS_{\theta}\cap\ssN$,
then $(t+\tau,x+\sfr_{\Square}(\tau))\in\GS_{\theta}$ for all $\tau\in[0,\eps)$
and $\Square\in\LR$. 
\end{prop}

\begin{proof}
Let $X_{\Square}\coloneqq X_{\theta,t,x,\Square}$, $m_{\Square}\coloneqq X_{\Square}'(t-)-\theta$,
and $s_{\Square}\coloneqq t-t_{-;1}[Y_{\Square}]$. We have, for $\tau>0$
and $\eta\in\RR$, that
\begin{align*}
\mathcal{A}_{\theta,T_{*}(t),t+\tau} & [\mathcal{T}_{\tau,\eta}^{2}X_{\Square}]-\mathcal{A}_{\theta,T_{*}(t),t+\tau}[\mathcal{T}_{\tau,\eta}^{1}X_{\Square}]\\
 & =\frac{1}{2}\tau\left(\frac{\eta}{\tau}-\theta\right)^{2}-\mu(\{(t,x)\})+\mathcal{A}_{\theta,T_{*}(t),t}[Y_{\Square}]-\mathcal{A}_{\theta,T_{*}(t),t+\tau}[\mathcal{T}_{\tau,\eta}^{1}X_{\Square}]\\
\oversetx{\prettyref{eq:AX}} & =\frac{(\eta-\tau\theta)^{2}}{2\tau}-\mu(\{(t,x)\})+\frac{1}{2}s_{\Square}m_{\Square}^{2}-\frac{1}{2}(s_{\Square}+\tau)\left(m_{\Square}+\frac{\eta-\tau(m_{\Square}+\theta)}{s_{\Square}+\tau}\right)^{2}\\
 & =\frac{s_{\Square}}{2\tau(s_{\Square}+\tau)}\left((\eta-\tau\theta)^{2}+2\tau m_{\Square}(\eta-\tau\theta)+\tau^{2}m_{\Square}^{2}\right)-\mu(\{(t,x)\})\\
 & =\frac{s_{\Square}}{2\tau(s_{\Square}+\tau)}\left(\eta-\tau\theta+\tau m_{\Square}\right)^{2}-\mu(\{(t,x)\}).
\end{align*}
Therefore, if we define
\[
\sfr_{\Square}(\tau)\coloneqq\tau(\theta-m_{\Square})\pm_{\Square}\sqrt{2\tau(1+\tau/s_{\Square})\mu(\{(t,x)\})}
\]
(with $\pm_{\Square}$ as in \prettyref{eq:pmbookkeeping}), then
the fact that $\sfr_{\Left}(\tau)<\sfr_{\Right}(\tau)$ is clear,
and we moreover have for sufficiently small $\tau>0$ that
\begin{equation}
\mathcal{A}_{\theta,T_{*}(t),t+\tau}[\mathcal{T}_{\tau,\sfr_{\Square}(\tau)}^{2}X_{\Square}]=\mathcal{A}_{\theta,T_{*}(t),t+\tau}[\mathcal{T}_{\tau,\sfr_{\Square}(\tau)}^{1}X_{\Square}]\qquad\text{for }\Square\in\LR.\label{eq:setequal}
\end{equation}

We note that 
\begin{equation}
\mathcal{A}_{\theta,T_{*}(t),t+\tau}[\mathcal{T}_{\tau,\sfr_{\Square}(\tau)}^{2}X]\text{ is independent of }\Square\in\LR.\label{eq:independentof}
\end{equation}
Also, if $X\in\mathscr{M}_{t,x}^{\theta}$ and $X'(t-)>X'_{\theta,t,x,\Right}(t-)$,
then \prettyref{eq:generalfderiv} implies that, for sufficiently
small $\tau>0$, we have
\begin{equation}
\mathcal{A}_{\theta,T_{*}(t),t+\tau}[\mathcal{T}_{\tau,\sfr_{\Square}(\tau)}^{1}X]>\mathcal{A}_{\theta,T_{*}(t),t+\tau}[\mathcal{T}_{\tau,\sfr_{\Square}(\tau)}^{1}X_{\theta,t,x,\Right}],\label{eq:biggerthanright}
\end{equation}
and similarly if $X'(t-)<X'_{\theta,t,x,\Left}(t-)$, then for sufficiently
small $\tau>0$ we have
\begin{equation}
\mathcal{A}_{\theta,T_{*}(t),t+\tau}[\mathcal{T}_{\tau,\sfr_{\Square}(\tau)}^{1}X]>\mathcal{A}_{\theta,T_{*}(t),t+\tau}[\mathcal{T}_{\tau,\sfr_{\Square}(\tau)}^{1}X_{\theta,t,x,\Left}].\label{eq:biggerthanleft}
\end{equation}
Now using \prettyref{prop:main-continuity} along with \prettyref{eq:setequal},
\prettyref{eq:independentof}, \prettyref{eq:biggerthanright}, and
\prettyref{eq:biggerthanleft}, we conclude \prettyref{eq:usesameminimums-split}.

The last claim of the proposition statement is thus clear, since in
this case, at least if $\tau$ is sufficiently small, then $\mathcal{T}_{\tau,\sfr_{\Left}(\tau)}^{1}X_{\theta,t,x,\Left}$
and $\mathcal{T}_{\tau,\sfr_{\Left}(\tau)}^{2}X_{\theta,t,x,\Right}$
will satisfy the conditions of \prettyref{def:globalshocks}, as will
$\mathcal{T}_{\tau,\sfr_{\Right}(\tau)}^{1}X_{\theta,t,x,\Right}$
and $\mathcal{T}_{\tau,\sfr_{\Right}(\tau)}^{2}X_{\theta,t,x,\Left}$.
\end{proof}
Now we look at what happens just before time $t$. The following proposition
justifies the picture in \prettyref{fig:shockmoveintime_picture1}.
\begin{prop}
\label{prop:backwards}Suppose that $\mu\in\widetilde{\Omega}_{1}$
and that $\theta\in\RR$, $(t,x)\in\ssS_{\theta}$, and $X_{1},X_{2}\in\scM_{t,x}^{\theta}$
are such that $X_{1}'(t-)>X_{2}'(t-)$ and, for any $X\in\scM_{t,x}^{\theta}$,
we have $X'(t-)\not\in(X_{2}'(t-),X_{1}'(t-))$. Then there is an
$\eps=\eps(\theta,t,x)>0$ and a function $\sfr_{X_{1},X_{2}}\colon(-\eps,0]\to\TT$
such that $\sfr(0)=0$,
\begin{equation}
\sfr_{X_{1},X_{2}}'(0-)=\frac{X_{1}'(t-)+X_{2}'(t-)}{2},\label{eq:rderiv-LR-1}
\end{equation}
and, for all $\tau\in(-\eps,0)$,
\begin{equation}
\scM_{t+\tau,x+\sfr_{X_{1},X_{2}}(\tau)}^{\theta}=\left\{ \mathcal{T}_{\tau,\sfr_{X_{1},X_{2}}(\tau)}^{1}X\st X\in\scM_{t,x}^{\theta}\text{ and }X'(t-)\in\{X_{1}'(t-),X_{2}'(t-)\}\right\} .\label{eq:useminimizers-before}
\end{equation}
\end{prop}

\begin{proof}
We choose $\eps>0$, $T_{*}=T_{*}(t-\eps)$, and $\sfr\coloneqq\sfr_{X_{1},X_{2}}$
as in \prettyref{lem:differentiate-A}, so \prettyref{eq:rderiv-LR-1}
is simply \prettyref{eq:RHcondition}. If $X\in\scM_{t,x}^{\theta}$
is such that $X'(t-)\not\in[X_{2}'(t-),X_{1}'(t-)]$, then we can
use \prettyref{eq:generalfderiv} and \prettyref{eq:rderiv-LR-1}
to obtain
\begin{align*}
\frac{\dif}{\dif\tau}\left[\mathcal{A}_{\theta,T_{*},t+\tau}[\mathcal{T}_{\tau,\sfr(\tau)}^{1}X]-\mathcal{A}_{\theta,T_{*},t+\tau}[\mathcal{T}_{\tau,\sfr(\tau)}^{1}X_{1}]\right]_{\tau=0} & =\frac{1}{2}(X_{2}'(t-)-X'(t-))(X'(t-)-X_{1}'(t-))\\
 & <0.
\end{align*}
On the other hand, if $X'(t-)\in\{X_{1}'(t-),X_{2}'(t-)\}$, then
it is clear from the definitions and \prettyref{eq:Aequal} that 
\[
\mathcal{A}_{\theta,T_{*},t+\tau}[\mathcal{T}_{\tau,\sfr(\tau)}^{1}X]=\mathcal{A}_{\theta,T_{*},t+\tau}[\mathcal{T}_{\tau,\sfr(\tau)}^{1}X_{1}]=\mathcal{A}_{\theta,T_{*},t+\tau}[\mathcal{T}_{\tau,\sfr(\tau)}^{1}X_{2}]\qquad\text{for all }\tau\in(-\eps,\eps).
\]
Thus we can take $\eps$ smaller if necessary and apply \prettyref{prop:main-continuity}
along with the last two displays to see that \prettyref{eq:useminimizers-before}
holds for all $\tau\in(-\eps,0)$. 
\end{proof}

\subsection{Global shock movement}

\label{subsec:movement-global}We now use the results of the previous
subsection to describe the movement in time of $\sfs_{\Left}(\theta,\cdot)$
and $\sfs_{\Right}(\theta,\cdot)$. The first part of the following
proposition is \prettyref{thm:shockstructure}(\ref{enu:shockcontinuity}).
At the end of this section, we also prove \prettyref{thm:shockstructure}(\ref{enu:splitatforcing}).
\begin{prop}
\label{prop:smovement}Suppose that $\mu\in\widetilde{\Omega}$. Fix
$\theta\in\RR$ and $\Square\in\LR$.
\begin{enumerate}
\item \label{enu:scts}The function $t\mapsto\sfs_{\Square}(\theta,t)$
is continuous in $t$.
\item \label{enu:sderiv}The function $t\mapsto\sfs_{\Square}(\theta,t)$
is right-differentiable at every $t\in\RR$ such that $(t,\sfs_{\Square}(\theta,t))\not\in\ssN$,
and for each such $t$ we have
\begin{equation}
\partial_{t}\sfs_{\Square}(\theta,t+)=\frac{1}{2}\sum_{\Diamond\in\LR}u_{\theta,\Diamond}(t,\sfs_{\Square}(\theta,t)).\label{eq:sderivintime}
\end{equation}
\end{enumerate}
\end{prop}

\begin{proof}
Let $t\in\RR$. We consider two cases.
\begin{casenv}
\item First we consider the case when $\sfs_{\Left}(\theta,t)=\sfs_{\Right}(\theta,t)$.
Let $\Square\in\LR$. Suppose for the sake of contradiction that there
is a sequence $t_{k}\to t$ such that $\sfs_{\Square}(\theta,t_{k})$
does not converge to $\sfs_{\Square}(\theta,t)$ as $k\to\infty$.
By the compactness of $\TT$, we can pass to a subsequence to assume
that $x\coloneqq\lim\limits_{k\to\infty}\sfs_{\Square}(\theta,t_{k})\ne\sfs_{\Square}(\theta,t)$
exists. But then \prettyref{prop:GSclosed} implies that $x\in\GS_{\theta,t}$,
whereas the assumption that $\sfs_{\Left}(\theta,t)=\sfs_{\Right}(\theta,t)$
implies that $\GS_{\theta,t}=\{\sfs_{\Square}(\theta,t)\}$, a contradiction.
Therefore, we must have
\[
\lim_{t'\to t}\sfs_{\Square}(\theta,t')=\sfs_{\Square}(\theta,t),
\]
and hence $\sfs_{\Square}(\theta,\cdot)$ is continuous at $t$. If
we now moreover assume that $(t,\sfs_{\Square}(\theta,t))\not\in\ssN$,
then \prettyref{prop:advanceint} implies that, if we take $\eps$
and $\sfr$ as in that proposition, we have $\sfs_{\Square}(\theta,t+\tau)=\sfs_{\Square}(\theta,t)+\sfr(\tau)$
for $\tau\in[0,\eps)$, and then \prettyref{eq:advanceint-rderiv}
implies \prettyref{eq:sderivintime}.
\item Now we consider the case when $\sfs_{\Right}(\theta,t)\ne\sfs_{\Left}(\theta,t)$.
Note that this implies in particular that $\theta\in\Theta_{\otimes}$.
We proceed in two steps.
\begin{thmstepnv}
\item First we address continuity from above. For $\Square\in\LR$, define
$\eps_{\Square},\sfr_{\Square}$ as in \prettyref{prop:advanceint}
with $x=\sfs_{\Square}(\theta,t)$, and put $\eps=\eps_{\Left}\wedge\eps_{\Right}$,
so we know from that proposition that, for all $\tau\in(0,\eps)$,
we have
\begin{equation}
X_{\theta,t+\tau,\sfs_{\Square}(\theta,t)+\sfr_{\Square}(\tau),\Diamond}=\mathcal{T}_{\tau,\sfr_{\Square}(\tau)}^{1}X_{\theta,t,\sfs_{\Square}(\theta,t),\Diamond}\qquad\text{for }\Square,\Diamond\in\LR.\label{eq:itsatranslation}
\end{equation}
From this we see that, for $\Square\in\LR$, we have
\[
\int_{T_{\vee}(\theta,t+\tau,\sfs_{\Square}(\theta,t)+\sfr_{\Square}(\tau))}^{t+\tau}\left(X_{\theta,t+\tau,\sfs_{\Square}(\theta,t)+\sfr_{\Square}(\tau),\Left}'-X_{\theta,t+\tau,\sfs_{\Square}(\theta,t)+\sfr_{\Square}(\tau),\Right}'\right)(s)\,\dif s\ne0,
\]
so $(t+\tau,\sfs_{\Square}(\theta,t)+\sfr_{\Square}(\tau))\in\GS_{\theta}$.
Also, using \prettyref{eq:itsatranslation} and \prettyref{prop:ifnotequalthenwherecollide},
we have
\[
X_{\theta,t+\tau,\sfs_{\Right}(\theta,t)+\sfr_{\Right}(\tau),\Left}(s_{\otimes}(\theta))=X_{\theta,t,\sfs_{\Right}(\theta,t),\Left}(s_{\otimes}(\theta))=y_{\otimes}(\theta)
\]
and
\[
X_{\theta,t+\tau,\sfs_{\Left}(\theta,t)+\sfr_{\Left}(\tau),\Right}(s_{\otimes}(\theta))=X_{\theta,t,\sfs_{\Left}(\theta,t),\Right}(s_{\otimes}(\theta))=y_{\otimes}(\theta),
\]
which means that $\sfs_{\Square}(\theta,t+\tau)=\sfs_{\Square}(\theta,t)+\sfr_{\Square}(\tau)$
for each $\Square\in\LR$. Then the continuity from above follows
from the continuity of $\sfr_{\Square}(\tau)$. Moreover, if we assume
that $(t,\sfs_{\Square}(\theta,t))\not\in\ssN$ (which is in fact
guaranteed in this case since $s_{\otimes}(\theta)<t$ by \prettyref{prop:findsplitpoint}),
then \prettyref{eq:sderivintime} follows from \prettyref{eq:rderiv-LR-1}.
\item Now we address continuity from below. We will prove this for $\sfs_{\Right}$;
the proof for $\sfs_{\Left}$ is symmetrical. Continuity from below
is somewhat more delicate than continuity from above because multiple
shocks may be merging at time $t$, and so we need to figure out which
of the merging shocks is the global shock to be followed backward
in time. (At most one can be a global shock since we have assumed
that $\sfs_{\Left}(\theta,t)\ne\sfs_{\Right}(\theta,t)$, so we cannot
be at a point where global shocks merge.)

Let $\tilde{Z}_{\Left}$ be the rightmost element $Z$ of $\mathscr{M}_{t,\sfs_{\Right}(\theta,t)}^{\theta}$
such that 
\[
Z(s_{\otimes}(\theta))=X_{\theta,t,\sfs_{\Right}(\theta,t),\Left}(s_{\otimes}(\theta))=y_{\otimes}(\theta),
\]
and define the concatenated paths
\begin{equation}
Z_{\Left}\coloneqq X_{\theta,s_{\otimes}(\theta),y_{\otimes}(\theta),\Left}\odot_{s_{\otimes}(\theta)}\tilde{Z}_{\Left}\label{eq:ZLdef}
\end{equation}
and
\begin{equation}
Y\coloneqq X_{\theta,s_{\otimes}(\theta),y_{\otimes}(\theta),\Right}\odot_{s_{\otimes}(\theta)}\tilde{Z}_{\Left},\label{eq:Ydef}
\end{equation}
recalling the definition \prettyref{eq:odotdef} of $\odot$. Now
let 
\[
m_{0}=\max\{X'(t-)\st X\in\mathscr{M}_{t,\sfs_{\Right}(\theta,t)}^{\theta}\text{ and }X'(t-)<Z_{\Left}'(t-)\},
\]
and let $Z_{\Right}$ be the leftmost element $Z$ of $\mathscr{M}_{t,\sfs_{\Right}(\theta,t)}^{\theta}$
such that $Z'(t)=m_{0}$. (See \prettyref{fig:find_Zs_picture} for
an illustration.)
\begin{figure}[t]
\begin{centering}
\input{figures/find_Zs_picture.pgf}
\par\end{centering}
\caption[Identifying the global shock out of several merging shocks]{Three shocks, including the right global shock $\protect\sfs_{\protect\Right}(\theta,\cdot)$,
are merging at time $t$. In backward time, the right global shock
continues between $Z_{\protect\Left}$ and $Z_{\protect\Right}$ (drawn
with thicker black dashed lines).\label{fig:find_Zs_picture}}
\end{figure}

From the definitions, we have $Z_{\Left}'(t-)>Z_{\Right}'(t-)$ and,
if $X\in\scM_{t,\sfs_{\Right}(\theta,t)}^{\theta}$, then $X'(t-)\not\in(Z_{\Right}'(t-),Z_{\Left}'(t-))$.
Thus, \prettyref{prop:backwards} applies, so taking $\eps>0$ and
$\sfr=\sfr_{Z_{\Left},Z_{\Right}}$ as in that proposition, we see
that for any $\tau\in(-\eps,0)$, we have
\begin{equation}
\scM_{t+\tau,\sfs_{\Right}(\theta,t)+\sfr(\tau)}^{\theta}=\left\{ \mathcal{T}_{\tau,\sfr(\tau)}^{1}X\st X\in\mathscr{M}_{t,\sfs_{\Right}(\theta,t)}^{\theta}\text{ and }X'(t-)\in\{Z_{\Left}'(t-),Z_{\Right}'(t-)\}\right\} .\label{eq:applybackward}
\end{equation}

Now we observe that 
\begin{equation}
X\preccurlyeq Z_{\Right}\text{ for all }X\in\mathscr{M}_{t,\sfs_{\Right}(\theta,t)}^{\theta}\text{ such that }X'(t-)=Z_{\Right}'(t-),\label{eq:ZRRightmost}
\end{equation}
which is clear from the definition of $Z_{\Right}$, and also that
\begin{equation}
X\succcurlyeq Z_{\Left}\text{ for all }X\in\mathscr{M}_{t,\sfs_{\Right}(\theta,t)}^{\theta}\text{ such that }X'(t-)=Z_{\Left}'(t-).\label{eq:ZLleftmost}
\end{equation}
To see \prettyref{eq:ZLleftmost}, we note that if there is some $r\le t$
such that $X(r)=Z_{\Left}(r)$ but $X'(r-)\ne Z_{\Left}'(r-)$, then
by \prettyref{lem:whencanyousplit} we must have $r=s_{\otimes}(\theta)$
and $X(r)=Z_{\Left}(r)=y_{\otimes}(\theta)$. This means that $X|_{[s_{\otimes}(\theta),t]}=Z_{\Left}|_{[s_{\otimes}(\theta),t]}$.
Then the fact that $X\succcurlyeq Z_{\Left}$ follows from the fact
that $Z_{\Left}|_{(-\infty,s_{\otimes}(\theta))}$ is the leftmost
minimizer from $(s_{\otimes}(\theta),y_{\otimes}(\theta))$ by the
definition \prettyref{eq:ZLdef}. Using \prettyref{eq:ZRRightmost}
and \prettyref{eq:ZLleftmost} in \prettyref{eq:applybackward}, we
see that for any $\tau\in(-\eps,0)$, we have
\begin{equation}
X_{\theta,t+\tau,\sfs_{\Right}(\theta,t)+\sfr(\tau),\Square}=\mathcal{T}_{\tau,\sfr(\tau)}^{1}Z_{\Square}\qquad\text{for }\Square\in\LR.\label{eq:leftrightisZ}
\end{equation}

We also note that
\begin{equation}
\int_{r}^{s_{\otimes}(\theta)}Z_{\Right}'(s)\,\dif s\le\int_{r}^{s_{\otimes}(\theta)}Y'(s)\,\dif s\le\int_{r}^{s_{\otimes}(\theta)}Z_{\Left}'(s)\,\dif s\quad\text{for all }r\in[T_{*}(t),t].\label{eq:Zordering}
\end{equation}
The second inequality is clear from the definitions, and the first
is a restatement of \prettyref{eq:ZLleftmost}. Using this, we see
that if we define
\[
T\coloneqq\sup\{s<t\st Z_{\Left}(s)=Z_{\Right}(s)\},
\]
then
\begin{equation}
\int_{T}^{t}(Z_{\Left}'-Z_{\Right}')(s)\,\dif s=1.\label{eq:Zsintegral}
\end{equation}
Indeed, the fact that this integral is positive follows from \prettyref{eq:Zordering},
the definitions \prettyref{eq:ZLdef} and \prettyref{eq:Ydef}, and
the fact that $(s_{\otimes}(\theta),y_{\otimes}(\theta))\in\GS_{\theta}$,
and then the definition of $T$ implies that it is an integer and
at most $1$. From \prettyref{eq:Zsintegral} we conclude that, for
any $\tau\in(-\eps,0)$, we have
\[
\int_{T}^{t}((\mathcal{T}_{\tau,\sfr(\tau)}^{1}Z_{\Left})'-(\mathcal{T}_{\tau,\sfr(\tau)}^{1}Z_{\Right})')(s)\,\dif s=1
\]
and hence that $(t+\tau,\sfs_{\Right}(\theta,t)+\sfr(\tau))\in\GS_{\theta}$
in light of \prettyref{eq:leftrightisZ}. Moreover, since 
\[
X_{\theta,t+\tau,\sfs(\theta,t)+\sfr(\tau),\Left}(s_{\otimes}(\theta))=Z_{\Left}(s_{\otimes}(\theta))=y_{\otimes}(\theta)
\]
by \prettyref{eq:leftrightisZ} and the definition of $Z_{\Left}$,
we in fact have $\sfs_{\Right}(t+\tau)=\sfs_{\Right}(\theta,t)+\sfr(\tau)$
for all $\tau\in(-\eps,0)$. Thus the continuity from below of $\sfs_{\Right}$
at $t$ follows from the continuity of $\sfr$.\qedhere
\end{thmstepnv}
\end{casenv}
\end{proof}
Now we can prove \prettyref{thm:shockstructure}(\ref{enu:splitatforcing}). 
\begin{proof}[Proof of \prettyref{thm:shockstructure}\textup{(\ref{enu:splitatforcing})}.]
If $\theta\in\RR\setminus\Theta_{\otimes}$, then $\sfs_{\Left}(\theta,t)=\sfs_{\Right}(\theta,t)$
for all $t\in\RR$ by \prettyref{prop:ifnotequalthenwherecollide}.

Thus we now assume that $\theta\in\Theta_{\otimes}$. Let $\mathsf{T}\coloneqq\{t\in\RR\st\sfs_{\Left}(\theta,t)\ne\sfs_{\Right}(\theta,t)\}$.
By \prettyref{prop:smovement}(\ref{enu:scts}), $\mathsf{T}$ is
an open subset of $\RR$. Also, \prettyref{prop:findsplitpoint} tells
us that $\mathsf{T}\subseteq(s_{\otimes}(\theta),\infty)$.

We claim that if $t_{k}\in\mathsf{T}$ and $t_{k}\downarrow t\not\in\mathsf{T}$
as $k\to\infty$, then in fact $t=s_{\otimes}(\theta)$. The continuity
of $\sfs_{\Left}(\theta,\cdot)$ and $\sfs_{\Right}(\theta,\cdot)$
implies that $\lim\limits_{k\to\infty}\sfs_{\Left}(\theta,t_{k})=\lim\limits_{k\to\infty}\sfs_{\Right}(\theta,t_{k})\eqqcolon y$.
But this is implies that $(t,y)\in\ssN$ by the last statement in
\prettyref{prop:advanceint}, and since the continuity of $\sfs_{\Square}(\theta,\cdot)$
implies that $(t,y)\in\GS_{\theta}$ as well, we have $t=s_{\otimes}(\theta)$
by the assumption that $\mu\in\widetilde{\Omega}_{2}$.

Thus we can conclude that $\mathsf{T}=(s_{\otimes}(\theta),s_{\wedge}(\theta))$
for some $s_{\wedge}(\theta)\in[s_{\otimes}(\theta),\infty]$. But
we know that $s_{\wedge}(\theta)<\infty$ by  \prettyref{prop:allshocksmerge},
and that $s_{\wedge}(\theta)>s_{\otimes}(\theta)$ by \prettyref{prop:advanceint-N}.
This completes the proof.
\end{proof}

\section{The dependence of \texorpdfstring{$u$}{u} on \texorpdfstring{$\theta$}{θ}}

\label{sec:udependenceontheta}Now we begin our study of how the structure
of the Burgers flow changes as $\theta$ changes, and ultimately prove
\prettyref{thm:mainthm}. The key observation is the following.
\begin{lem}
Suppose that $\mu\in\Omega$. If $s<t$, $X\colon[s,t]\to\TT$ and
$\theta_{1},\theta_{2}\in\RR$, then
\begin{align}
\mathcal{A}_{\theta_{2},s,t}[X]-\mathcal{A}_{\theta_{1},s,t}[X] & =(\theta_{1}-\theta_{2})\int_{s}^{t}X'(r)\,\dif r+\frac{1}{2}(t-s)(\theta_{2}^{2}-\theta_{1}^{2}).\label{eq:diffAtheta1theta2}
\end{align}
Also, for $\theta\in\RR$, we have
\begin{equation}
\frac{\dif}{\dif\theta}\mathcal{A}_{\theta,s,t}[X]=-\int_{s}^{t}X'(r)\,\dif r+\theta(t-s).\label{eq:difAXderiv}
\end{equation}
\end{lem}

\begin{proof}
This follows immediately from the definition \prettyref{eq:Adef}.
\end{proof}
The remainder of this section is divided into two parts. In \prettyref{subsec:only-change-at-global-shock},
we show that, as $\theta$ changes, $u_{\theta,\Square}(t,x)$ is
constant except when $(\theta,t,x)\in\GS$, and at such values of
$\theta$ the size of the jump is $u_{\theta,\Left,\Square}(t,x)-u_{\theta,\Right,\Square}(t,x)$.
See \prettyref{prop:differenceintermsoftopology} below. In \prettyref{subsec:Movement-of-shocks},
we relate $u_{\theta,\Left,\Square}(t,x)-u_{\theta,\Right,\Square}(t,x)$
to $\partial_{\theta}\sfs_{\square}(\theta,\cdot)$ (see \prettyref{prop:relatetooneoverderiv}
below), and then complete the remaining proofs.

\subsection{Jumps of \texorpdfstring{$u$}{u} occur at global shocks}

\label{subsec:only-change-at-global-shock}Recall the definition \prettyref{eq:udef}
of $u_{\theta,\Square}$, and the definitions \prettyref{eq:Mstartx}
of $\scM_{*\mid t,x}^{\theta}$ and \prettyref{eq:usquarediamonddef}
of $u_{\theta,\Square,\Diamond}(t,x)$. 
\begin{prop}
\label{prop:GS'discrete}Suppose that $\mu\in\widetilde{\Omega}_{1}$.
For each fixed $t\in\RR$ and $x\in\TT$, the set
\begin{equation}
\GS'_{t,x}\coloneqq\{\theta\in\RR\st(\theta,t,x)\in\GS\}\label{eq:GS'def}
\end{equation}
is discrete.
\end{prop}

It is important in the statement of \prettyref{prop:GS'discrete}
that we do \emph{not} assume that $\mu\in\widetilde{\Omega}_{2}$,
because we use \prettyref{prop:GS'discrete} in the proof of \prettyref{prop:omega21}
below.
\begin{proof}
To have $(\theta,t,x)\in\GS$, we must have paths $X,Y\colon[T_{*}(t),t]\to\TT$,
each consisting of straight line segments connecting points of $\ssN$,
such that
\begin{equation}
\int_{T_{*}(t)}^{t}X'(s)\,\dif s>\int_{T_{*}(t)}^{t}Y'(s)\,\dif s\label{eq:Xprimebigger}
\end{equation}
and
\begin{equation}
\mathcal{A}_{\theta,T_{*}(t),t}[X]=\mathcal{A}_{\theta,T_{*}(t),t}[Y].\label{eq:Asequal}
\end{equation}
By \prettyref{eq:difAXderiv}, we have
\[
\frac{\dif}{\dif\theta}\left[\mathcal{A}_{\theta,T_{*}(t),t}[X]-\mathcal{A}_{\theta,T_{*}(t),t}[Y]\right]=\int_{T_{*}(t)}^{t}(Y'-X')(s)\,\dif s\overset{\prettyref{eq:Xprimebigger}}{<}0,
\]
which means that for a given $X$ and $Y$, there is only at most
a single value of $\theta$ for which \prettyref{eq:Asequal} can
hold. Since, for $\theta$ in any bounded set, there are only finitely
many $X$ and $Y$ that can achieve \prettyref{eq:Asequal}, this
means that the set $\{\theta\in\RR\st(\theta,t,x)\in\GS\}$ must be
discrete.
\end{proof}
\begin{prop}
\label{prop:locallyconstant}Suppose that $\mu\in\widetilde{\Omega}_{1}$.
For each $\theta,t\in\RR$, $x\in\TT$, and $\Square\in\LR$, there
is an $\eps=\eps(\theta,t,x,\Square)>0$ such that, whenever $\theta-\eps<\theta_{-}<\theta<\theta_{+}<\theta+\eps$,
we have
\begin{equation}
\mathscr{M}_{*\mid t,x}^{\theta'}=\begin{cases}
\scM_{*\mid t,x,\Right}^{\theta} & \text{if }\theta'\in(\theta-\eps,\theta);\\
\scM_{*\mid t,x,\Left}^{\theta} & \text{if }\theta'\in(\theta,\theta+\eps),
\end{cases}\label{eq:howdominimizerschange}
\end{equation}
and hence that
\begin{equation}
u_{\theta_{+},\Square}(t,x)=u_{\theta,\Left,\Square}(t,x)\qquad\text{and}\qquad u_{\theta_{-},\Square}(t,x)=u_{\theta,\Right,\Square}(t,x).\label{eq:jumpsize-mo}
\end{equation}
\end{prop}

\begin{proof}
It is a consequence of \prettyref{prop:main-continuity} that there
is an $\eps>0$ such that if $|\theta'-\theta|<\eps$, then
\[
\mathscr{M}_{*\mid t,x}^{\theta'}=\left\{ X\in\scM_{*\mid t,x}^{\theta}\st\mathcal{A}_{\theta',T_{*}(t),t}[X]\le\mathcal{A}_{\theta',T_{*}(t),t}[Y]\text{ for all }Y\in\mathscr{M}_{*\mid t,x}^{\theta}\right\} .
\]
We have by \prettyref{eq:diffAtheta1theta2} that, for $X,Y\in\mathscr{M}_{*\mid t,x}^{\theta'}$,
\begin{align*}
 & \frac{\mathcal{A}_{\theta',T_{*}(t),t}[X]-\mathcal{A}_{\theta',T_{*}(t),t}[Y]}{\theta'-\theta}.\\
 & =\frac{\mathcal{A}_{\theta',T_{*}(t),t}[X]-\mathcal{A}_{\theta,T_{*}(t),t}[X]-\left(\mathcal{A}_{\theta',T_{*}(t),t}[Y]-\mathcal{A}_{\theta,T_{*}(t),t}[Y]\right)}{\theta'-\theta}\\
 & =-\int_{T_{*}(t)}^{t}(X'-Y')(s)\,\dif s\mathrel{\begin{cases}
=0 & \text{if and only if }X\sim Y;\\
\ge0 & \text{if }X\in\mathscr{M}_{*\mid t,x,\Right}^{\theta};\\
\le0 & \text{if }X\in\scM_{*\mid t,x,\Left}^{\theta},
\end{cases}}
\end{align*}
where we recalled \prettyref{def:equiv-relation} of $\sim$ and $\scM_{*\mid t,x,\Square}^{\theta}$.
Combining the last two displays, we conclude that \prettyref{eq:howdominimizerschange}
holds, and then \prettyref{eq:jumpsize-mo} follows from \prettyref{eq:howdominimizerschange}
and the definitions.
\end{proof}
The next statement gives us our first formula for $u_{\theta_{2},\Square}(t,x)-u_{\theta_{1},\Square}(t,x)$,
which is the subject of \prettyref{thm:mainthm}. Recall that $\llbracket\theta_{1},\theta_{2}\rrbracket_{\Square}$
was defined in \prettyref{eq:pmbookkeeping}.
\begin{prop}
\label{prop:differenceintermsoftopology}Suppose that $\mu\in\widetilde{\Omega}_{1}$.
Fix $t\in\RR$ and $x\in\TT$. We have, for $\Square\in\LR$, that
\begin{equation}
u_{\theta_{2},\Square}(t,x)-u_{\theta_{1},\Square}(t,x)=\sum_{\theta\in\llbracket\theta_{1},\theta_{2}\rrbracket_{\Square}\cap\GS'_{t,x}}\left(u_{\theta,\Left,\Square}(t,x)-u_{\theta,\Right,\Square}(t,x)\right).\label{eq:differenceintermsoftopology}
\end{equation}
\end{prop}

\begin{proof}
First we note that if $\theta\in\RR\setminus\GS'_{t,x}$, then $u_{\theta,\Left,\Diamond}(t,x)=u_{\theta,\Right,\Diamond}(t,x)$
for $\Diamond\in\LR$ by \prettyref{prop:globalshocknojump}, which
means by \prettyref{prop:locallyconstant} that there is some $\eps>0$
such that if $|\theta'-\theta|<\eps$, then
\[
u_{\theta',\Square}(t,x)=u_{\theta,\Left,\Square}(t,x)=u_{\theta,\Right,\Square}(t,x)=u_{\theta,\Square}(t,x)
\]
(with the last identity by \prettyref{eq:leftleftrightright}). This,
for each $\Square\in\LR$, the map $\theta\mapsto u_{\theta,\Square}(t,x)$
is constant on each connected component of $\mathbb{R}\setminus\GS'_{t,x}$.
Moreover, for each $\theta\in\GS'_{t,x}$, we have again by \prettyref{prop:locallyconstant}
that
\[
\lim_{\theta'\downarrow\theta}u_{\theta',\Square}(t,x)-\lim_{\theta'\uparrow\theta}u_{\theta',\Square}(t,x)=u_{\theta,\Left,\Square}(t,x)-u_{\theta,\Right,\Square}(t,x).
\]
So, for any $\Square\in\LR$, we have
\begin{align*}
u_{\theta_{2},\Square}(t,x)-u_{\theta_{1},\Square}(t,x) & =u_{\theta_{2},\Square}(t,x)-u_{\theta_{2},\Right,\Square}(t,x)+\sum_{\theta\in(\theta_{1},\theta_{2})\cap\GS'_{t,x}}\left(u_{\theta,\Left,\Square}(t,x)-u_{\theta,\Right,\Square}(t,x)\right)\\
 & \qquad+u_{\theta_{1},\Left,\Square}(t,x)-u_{\theta_{1},\Square}(t,x).
\end{align*}
We observe that the sum comprises only finitely many terms by \prettyref{prop:GS'discrete}.
For $\Square=\Right$, we note that $u_{\theta,\Right}(t,x)=u_{\theta,\Right,\Right}(t,x)$
by \prettyref{eq:leftleftrightright}, so we obtain
\[
u_{\theta_{2},\Right}(t,x)-u_{\theta_{1},\Right}(t,x)=\sum_{\theta\in[\theta_{1},\theta_{2})\cap\GS'_{t,x}}\left(u_{\theta,\Left,\Square}(t,x)-u_{\theta,\Right,\Square}(t,x)\right).
\]
Similarly, for $\Square=\Left$, we have $u_{\theta,\Left}(t,x)=u_{\theta,\Left,\Left}(t,x)$
by \prettyref{eq:leftleftrightright}, so we obtain
\[
u_{\theta_{2},\Left}(t,x)-u_{\theta_{1},\Left}(t,x)=\sum_{\theta\in(\theta_{1},\theta_{2}]\cap\GS'_{t,x}}\left(u_{\theta,\Left,\Left}(t,x)-u_{\theta,\Right,\Left}(t,x)\right).
\]
The conclusion \prettyref{eq:differenceintermsoftopology} is a summary
of the last two displays using the notation \prettyref{eq:pmbookkeeping}.
\end{proof}
\begin{prop}
\label{prop:shocksdiffdiscrete}Suppose that $\mu\in\widetilde{\Omega}$.
For any fixed $t\in\RR$, the set $\{\theta\in\RR\st\sfs_{\Left}(\theta,t)\ne\sfs_{\Right}(\theta,t)\}$
is discrete.
\end{prop}

\begin{proof}
Let $\mathsf{N}_{*}(t)\coloneqq([T_{*}(t),t]\times\TT)\cap\ssN$.
We know from \prettyref{prop:findsplitpoint} and \prettyref{eq:Tveedef}
that
\begin{equation}
\text{if }\sfs_{\Left}(\theta,t)\ne\sfs_{\Right}(\theta,t),\text{ then }\mathsf{N}_{*}(t)\cap\GS_{\theta}\ne\varnothing.\label{eq:implication}
\end{equation}
By \prettyref{prop:GS'discrete}, we see that for each $(s,y)\in\mathsf{N}_{*}(t)$,
the set $\{\theta\in\RR\st(s,y)\in\GS_{\theta}\}$ is discrete. But
since $\ssN$ is discrete, the set $\mathsf{N}_{*}(t)$ is finite,
and so
\[
\{\theta\in\RR\st\mathsf{N}_{*}(t)\cap\GS_{\theta}\ne\varnothing\}=\bigcup_{(s,y)\in\ssN(t)}\{\theta\in\RR\st(s,y)\in\GS_{\theta}\}
\]
is also discrete. Hence, using \prettyref{eq:implication}, we conclude
that $\{\theta\in\RR\st\sfs_{\Left}(\theta,t)\ne\sfs_{\Right}(\theta,t)\}$
is discrete, as claimed.
\end{proof}
We can also complete the proof of \prettyref{thm:mainthm}(\ref{enu:switching}).
\begin{proof}[Proof of \prettyref{thm:mainthm}\textup{(\ref{enu:switching})}]
We will just prove \prettyref{eq:rightstillshock}, as the proof
of \prettyref{eq:leftstillshock} is symmetrical. By \prettyref{prop:locallyconstant},
there is an $\eps>0$ such that, if $\theta'\in(\theta-\eps,\theta]$,
then $\mathscr{M}_{*\mid t,\sfs_{\Right}(\theta,t)}^{\theta'}=\scM_{*\mid t,\sfs_{\Right}(\theta,t),\Right}^{\theta}$,
and so in particular
\begin{equation}
u_{\theta',\Square}(t,\sfs_{\Right}(\theta,t))=u_{\theta,\Right,\Square}(t,\sfs_{\Right}(\theta,t))\qquad\text{for }\Square\in\LR.\label{eq:theta'equalstheta}
\end{equation}
By \prettyref{prop:findshocktostay}, we have 
\[
u_{\theta,\Right,\Left}(t,\sfs_{\Right}(\theta,t))>u_{\theta,\Right,\Right}(t,\sfs_{\Right}(\theta,t)),
\]
so \prettyref{eq:theta'equalstheta} implies that
\[
u_{\theta',\Left}(t,\sfs_{\Right}(\theta,t))>u_{\theta',\Right}(t,\sfs_{\Right}(\theta,t)),
\]
so $(\theta',t,\sfs_{\Right}(\theta,t))\in\ssS$ by the definition
\prettyref{eq:Sdef}.
\end{proof}

\subsection{Movement of shocks as \texorpdfstring{$\theta$}{θ} is varied}

\label{subsec:Movement-of-shocks}In this subsection, we will prove
the following proposition, which we will combine with \prettyref{prop:differenceintermsoftopology}
to complete the proof of \prettyref{thm:mainthm}(\ref{enu:derivintheta}).
\begin{prop}
\label{prop:relatetooneoverderiv}Suppose that $\mu\in\widetilde{\Omega}$
and fix $\theta,t\in\RR$. We have
\begin{equation}
u_{\theta,\Left,\Right}(t,\sfs_{\Right}(\theta,t))-u_{\theta,\Right,\Right}(t,\sfs_{\Right}(\theta,t))=\frac{1}{\partial_{\theta}\sfs_{\Right}(\theta+,t)}\label{eq:rightside}
\end{equation}
and
\begin{equation}
u_{\theta,\Left,\Left}(t,\sfs_{\Left}(\theta,t))-u_{\theta,\Right,\Left}(t,\sfs_{\Left}(\theta,t))=\frac{1}{\partial_{\theta}\sfs_{\Left}(\theta-,t)}.\label{eq:leftside}
\end{equation}
 In particular, the one-sided derivatives on the right sides of \prettyref{eq:rightside}
and \prettyref{eq:leftside} both exist.
\end{prop}

In this section we will only prove \prettyref{eq:rightside}, as the
proof of \prettyref{eq:leftside} is symmetrical. To simplify notation,
we make the abbreviation
\[
\mathcal{T}_{\eta}\coloneqq\mathcal{T}_{0,\eta}^{1},
\]
with $\mathcal{T}_{0,\eta}^{1}$ defined as in \prettyref{eq:Tidef},
and we also recall the definition \prettyref{eq:tminus1} of $t_{-;1}[X]$. 
\begin{lem}
Suppose that $\mu\in\widetilde{\Omega}$. Let $\theta,t\in\RR$ and
$x\in\TT$. There is an $\eps=\eps(\theta,t,x)>0$ such that, if $0<|\zeta|,|\beta\zeta|<\eps$
and $X\in\mathscr{M}_{*\mid t,x}^{\theta}$, then 
\begin{equation}
\begin{aligned}\mathcal{A}_{\theta+\zeta,T_{*}(t),t} & [\mathcal{T}_{\beta\zeta}X]-\mathcal{A}_{\theta,T_{*}(t),t}[X]-\frac{1}{2}\zeta^{2}(t-T_{*}(t)-2\beta)+\mu(\{(t,x)\})\\
 & =\zeta\left[\beta(X'(t-)-\theta)-\int_{T_{*}(t)}^{t}(X'(s)-\theta)\,\dif s\right]+\frac{\beta^{2}\zeta^{2}}{2(t-t_{-;1}[X])}.
\end{aligned}
\label{eq:Adifflem}
\end{equation}
\end{lem}

An immediate consequence of \prettyref{eq:Adifflem} is that if $X_{1},X_{2}\in\mathscr{M}_{*\mid t,x}^{\theta}$,
then
\begin{equation}
\begin{aligned}\mathcal{A}_{\theta+\zeta,T_{*}(t),t} & [\mathcal{T}_{\beta\zeta}X_{2}]-\mathcal{A}_{\theta+\zeta,T_{*}(t),t}[\mathcal{T}_{\beta\zeta}X_{1}]\\
 & =\zeta\left[\beta\left(X_{2}'(t-)-X'_{1}(t-)\right)-\int_{T_{*}(t)}^{t}(X_{2}'(s)-X_{1}'(s))\,\dif s\right]\\
 & \qquad+\frac{1}{2}\zeta^{2}\beta^{2}\left(\frac{1}{t-t_{-;1}[X_{2}]}-\frac{1}{t-t_{-;1}[X_{1}]}\right).
\end{aligned}
\label{eq:Adiff-double}
\end{equation}

\begin{prop}
\label{prop:findYs}Suppose that $\mu\in\widetilde{\Omega}$. Let
$\theta,t\in\RR$ and $x\in\TT$, and let $Y_{1}\coloneqq X_{\theta,t,x,\Left,\Right}$
and $\tilde{Y}_{2}\coloneqq X_{\theta,t,x,\Right,\Right}$. Let $s_{0}\coloneqq\sup\{s<t\st Y(s)=\tilde{Y}_{2}(s)\}$
and $Y_{2}\coloneqq Y_{1}\odot_{s_{0}}\tilde{Y}_{2}$. Then we have
\begin{equation}
\int_{T_{*}(t)}^{t}(Y_{1}'-Y_{2}')(s)\,\dif s=\int_{s_{0}}^{t}(Y_{1}'-Y_{2}')(s)\,\dif s\in\{0,1\}\label{eq:Y1Y2diff}
\end{equation}
 and, for any $Z\in\mathscr{M}_{*\mid t,x}^{\theta}$, there is an
$i\in\{1,2\}$ such that
\begin{equation}
Z'(t-)-Y_{i}'(t-)\ge0\qquad\text{and}\qquad\int_{T_{*}(t)}^{t}(Z'-Y_{i}')(s)\,\dif s\le0.\label{eq:ourconstruct}
\end{equation}
 
\end{prop}

\begin{proof}
It is clear from the definitions that the first two expressions in
\prettyref{eq:Y1Y2diff} are equal, and moreover that they are in
$\{-1,0,1\}$. Moreover, if they were equal to $-1$, then we would
have
\[
-1=\int_{T_{*}(t)}^{t}(Y_{1}'-Y_{2}')(s)\,\dif s=\int_{T_{*}(t)}^{t}(X_{\theta,t,x,\Left}'-X'_{\theta,t,x,\Right})(s)\,\dif s\ge0,
\]
a contradiction. Thus we conclude that \prettyref{eq:Y1Y2diff} holds.

If $Z\in\mathscr{M}_{*\mid t,x,\Left}^{\theta}$, then $\int_{T_{*}(t)}^{t}(Z'-Y_{1}')(s)\,\dif s=0$
and $Z'(t-)\ge Y_{1}'(t-)$, so \prettyref{eq:ourconstruct} holds
with $i=1$ in this case. Thus we now assume that $Z\not\in\mathscr{M}_{*\mid t,x,\Left}^{\theta}$.
This implies in particular that
\[
\int_{T_{*}(t)}^{t}(Z'-Y_{1}')(s)\,\dif s\le-1.
\]
Combining this with \prettyref{eq:Y1Y2diff}, we see that
\[
\int_{T_{*}(t)}^{t}(Z'-Y_{2}')(s)\,\dif s\le0.
\]
Since it is also clear that $Z'(t-)\ge\tilde{Y}_{2}'(t-)=Y_{2}'(t-)$,
we see that \prettyref{eq:ourconstruct} holds with $i=2$.
\end{proof}
\begin{prop}
Suppose that $\mu\in\widetilde{\Omega}_{1}$. Suppose that $(\theta,t,x)\in\GS$
and
\begin{equation}
u_{\theta,\Left,\Right}(t,x)>u_{\theta,\Right,\Right}(t,x).\label{eq:canmoveshock}
\end{equation}
Then there is an $\eps=\eps(\theta,t,x)>0$ and a function $\sfr\colon[0,\eps)\to\TT$
such that $\sfr(0)=x$, $(\theta+\zeta,t,\sfr(\zeta))\in\GS$ for
each $\zeta\in[0,\eps)$, and
\begin{equation}
\sfr'(0+)=\frac{1}{u_{\theta,\Left,\Right}(t,x)-u_{\theta,\Right,\Right}(t,x)}.\label{eq:rderiv}
\end{equation}
\end{prop}

\begin{proof}
The proof proceeds in several steps.
\begin{thmstepnv}
\item By \prettyref{prop:main-continuity}, there is an $\eps>0$ such that,
if $|\beta\zeta|<\eps$ and $|\theta'-\theta|<\eps$, then $x+\beta\zeta\in\GS_{\theta',t}$
if there exist $Y_{1},Y_{2}\in\mathscr{M}_{*\mid t,x}^{\theta}$ such
that
\begin{equation}
(\mathcal{T}_{\beta\zeta}Y_{2})'(t-)>(\mathcal{T}_{\beta\zeta}Y_{1})'(t-),\qquad\int_{T_{*}(t)}^{t}[(\mathcal{T}_{\beta\zeta}Y_{2})'-(\mathcal{T}_{\beta\zeta}Y_{1})'](s)\,\dif s>0,\label{eq:globalshock}
\end{equation}
and, for all $Z\in\mathscr{M}_{*\mid t,x}^{\theta}$, we have 
\begin{equation}
\mathcal{A}_{\theta',T_{*}(t),t}[\mathcal{T}_{\beta\zeta}Z]\ge\mathcal{A}_{\theta',T_{*}(t),t}[\mathcal{T}_{\beta\zeta}Y_{1}]=\mathcal{A}_{\theta',T_{*}(t),t}[\mathcal{T}_{\beta\zeta}Y_{2}].\label{eq:globalshockm}
\end{equation}

We also assume that $\eps$ is chosen small enough that, if $X\in\mathscr{M}_{*\mid t,x}^{\theta}$
and $|\beta\zeta|<\eps$, then there are no points of $\ssN$ on $\mathcal{T}_{\beta\zeta}X$
that are not also on $X$; this is possible by the discreteness of
$\ssN$.
\item We select $Y_{1},\tilde{Y}_{2},Y_{2}$ as in the statement of \prettyref{prop:findYs}.
We note that 
\begin{equation}
Y_{2}'(t-)<Y_{1}'(t-)\label{eq:Yderivorder}
\end{equation}
by the assumption \prettyref{eq:canmoveshock}. We claim that
\begin{equation}
\int_{T_{*}(t)}^{t}(Y_{1}'-Y_{2}')(s)\,\dif s=1.\label{eq:diffisL}
\end{equation}
In light of \prettyref{eq:Y1Y2diff}, it suffices to show that 
\[
\int_{T_{*}(t)}^{t}(Y_{1}'-Y_{2}')(s)\,\dif s\ne0,
\]
but this is true because otherwise we would have $Y_{2}\in\scM_{*\mid t,x,\Left}^{\theta}$
and then \prettyref{eq:Yderivorder} would contradict the definition
of $Y_{1}$ as the rightmost element of $\mathscr{M}_{*\mid t,x,\Left}^{\theta}$.
Now \prettyref{eq:Yderivorder} and \prettyref{eq:diffisL} imply
that \prettyref{eq:globalshock} holds.
\item By \prettyref{eq:Adiff-double} and \prettyref{eq:diffisL}, we have
\begin{equation}
\mathcal{A}_{\theta+\zeta,T_{*}(t),t}[\mathcal{T}_{\beta\zeta}Y_{2}]-\mathcal{A}_{\theta+\zeta,T_{*}(t),t}[\mathcal{T}_{\beta\zeta}Y_{1}]=\zeta\left(\beta\left(Y_{2}'(t-)-Y'_{1}(t-)\right)-1+\beta^{2}Q\zeta/2\right),\label{eq:Adiff-double-1}
\end{equation}
where we have defined
\[
Q\coloneqq\frac{1}{t-t_{-;1}[Y_{2}]}-\frac{1}{t-t_{-;1}[Y_{1}]}.
\]
Now define, as long as $\zeta$ is sufficiently small,
\begin{equation}
\beta_{\zeta}\coloneqq\begin{cases}
(Y_{2}'(t-)-Y_{1}'(t-))\cdot\frac{\left(1+2Q\zeta\left(Y_{2}'(t-)-Y_{1}'(t-)\right)^{-2}\right)^{1/2}-1}{Q\zeta}, & \text{if }Q\zeta\ne0;\\
\left(Y_{2}'(t-)-Y_{1}'(t-)\right)^{-1}, & \text{if }Q\zeta=0.
\end{cases}\label{eq:betazetadef}
\end{equation}
Using \prettyref{eq:betazetadef} in \prettyref{eq:Adiff-double-1},
we see (for $\zeta$ small enough that \prettyref{eq:betazetadef}
is well-defined) that
\[
\mathcal{A}_{\theta+\zeta,T_{*}(t),t}[\mathcal{T}_{\beta_{\zeta}\zeta}Y_{2}]=\mathcal{A}_{\theta+\zeta,T_{*}(t),t}[\mathcal{T}_{\beta_{\zeta}\zeta}Y_{1}],
\]
which verifies the identity in \prettyref{eq:globalshockm}. We also
observe that, for small $\zeta$, we have the Taylor expansion
\begin{equation}
\beta_{\zeta}=\frac{1}{Y_{2}'(t-)-Y_{1}'(t-)}-\frac{Q\zeta}{2\left(Y_{2}'(t-)-Y_{1}'(t-)\right)^{3}}+O(\zeta^{2}).\label{eq:betazetaTaylor}
\end{equation}
\item Now we want to verify the inequality in \prettyref{eq:globalshockm}.
Let $Z\in\mathscr{M}_{*\mid t,x}^{\theta}$. By \prettyref{prop:findYs},
we can find an $i\in\{1,2\}$ such that \prettyref{eq:ourconstruct}
holds. Using \prettyref{eq:Adiff-double}, we can compute, for $\zeta\ge0$,
that
\begin{align}
 & \mathcal{A}_{\theta+\zeta,T_{*}(t),t}[\mathcal{T}_{\beta_{\eta}\eta}Z]-\mathcal{A}_{\theta+\zeta,T_{*}(t),t}[\mathcal{T}_{\beta_{\eta}\eta}Y_{i}]\nonumber \\
 & \ =\zeta\beta_{\zeta}\left(Z'(t-)-Y_{i}'(t-)\right)-\zeta\int_{T_{*}(t)}^{t}(Z'-Y_{i}')(s)\,\dif s+\frac{1}{2}\zeta^{2}\beta_{\zeta}^{2}\left(\frac{1}{t-t_{-;1}[Z]}-\frac{1}{t-t_{-;1}[Y_{i}]}\right)\nonumber \\
 & \ \ge\zeta\beta_{\zeta}\left(Z'(t-)-Y_{i}'(t-)\right)+\frac{1}{2}\zeta^{2}\beta_{\zeta}^{2}\left(\frac{1}{t-t_{-;1}[Z]}-\frac{1}{t-t_{-;1}[Y_{i}]}\right),\label{eq:derivZY}
\end{align}
where the inequality is by the second inequality in \prettyref{eq:ourconstruct}.
Now if the first inequality in \prettyref{eq:ourconstruct} is strict,
then the right side of \prettyref{eq:derivZY} is (strictly) positive
for sufficiently small $\zeta>0$ (and is zero for $\zeta=0$). On
the other hand, if $Z'(t-)=Y_{i}'(t-)$, then $t_{-;1}[Z]=t_{-;1}[Y_{i}]$
as well (recalling the definition \prettyref{eq:tminus1}), and so
in this case the right side of \prettyref{eq:derivZY} is zero. Therefore,
the inequality in \prettyref{eq:globalshockm} is verified.
\item We now define $\mathsf{r}(\zeta)=x+\beta_{\zeta}\zeta$. For $\zeta\ge0$,
the conditions \prettyref{eq:globalshock} and \prettyref{eq:globalshockm}
have been verified and so we have $(\theta+\zeta,t,\mathsf{r}(\zeta))\in\GS$.
The derivative \prettyref{eq:rderiv} follows from \prettyref{eq:betazetaTaylor},
and so the proof is complete.\qedhere
\end{thmstepnv}
\end{proof}
We can now complete the proof of \prettyref{prop:relatetooneoverderiv}
and in fact simultaneously prove \prettyref{eq:sLsRareleftrightlimits}.
\begin{proof}[Proof of \prettyref{prop:relatetooneoverderiv} and \prettyref{eq:sLsRareleftrightlimits}]
As noted above, we will only prove \prettyref{eq:rightside} of \prettyref{prop:relatetooneoverderiv},
since the proof of \prettyref{eq:leftside} is symmetrical. Similarly,
for the proof of \prettyref{eq:sLsRareleftrightlimits}, we will only
prove that
\begin{equation}
\sfs_{\Right}(\theta,t)=\lim_{\theta'\downarrow\theta}\sfs_{\Diamond}(\theta,t),\qquad\Diamond\in\LR,\label{eq:sctsgoal}
\end{equation}
as the proof of the other limit is again symmetrical. By \prettyref{prop:findleftrightshocks},
we have
\[
u_{\theta,\Left,\Right}(t,\sfs_{\Right}(\theta,t))\ne u_{\theta,\Right,\Right}(t,\sfs_{\Right}(\theta,t)),
\]
which means that \prettyref{eq:canmoveshock} is satisfied and \prettyref{prop:findYs}
applies with $x=\sfs_{\Right}(\theta,t)$. Therefore, we have an $\eps>0$
and a function $\sfr\colon[0,\eps)\to\TT$ such that
\begin{equation}
\sfr(0)=\sfs_{\Right}(\theta,t),\label{eq:r0issR}
\end{equation}
\begin{equation}
(\theta+\zeta,t,\sfr(\zeta))\in\GS\qquad\text{for each }\zeta\in[0,\eps),\label{eq:inGS}
\end{equation}
and
\begin{align}
\sfr'(0+) & =\frac{1}{u_{\theta,\Left,\Right}(t,\sfs_{\Right}(\theta,t))-u_{\theta,\Right,\Right}(t,\sfs_{\Right}(\theta,t))}.\label{eq:rprime0}
\end{align}
Now \prettyref{prop:shocksdiffdiscrete} and \prettyref{eq:inGS}
tell us that, by reducing $\eps$ if necessary, we can assume that
$\sfr(\zeta)=\sfs_{\Left}(\theta+\zeta,t)=\sfs_{\Right}(\theta+\zeta,t)$
for all $\zeta\in(0,\eps)$. Combining this with \prettyref{eq:r0issR},
we see that in fact $\sfr(\zeta)=\sfs_{\Right}(\theta+\zeta,t)$ for
all $\zeta\in[0,\eps)$. The continuity of $\sfr$ then implies \prettyref{eq:sctsgoal},
and the conclusion \prettyref{eq:rightside} of \prettyref{prop:relatetooneoverderiv}
is now simply \prettyref{eq:rprime0}.
\end{proof}
Finally, we can prove \prettyref{thm:mainthm}(\ref{enu:derivintheta}).
\begin{proof}[Proof of \prettyref{thm:mainthm}\textup{(\ref{enu:derivintheta})}]
We assume that $\Square=\Right$; the proof in the case $\Square=\Left$
is analogous. We have by \prettyref{prop:differenceintermsoftopology}
that, for any $\theta_{1}<\theta_{2}$, $t\in\RR$, and $x\in\TT$,
we have
\begin{equation}
u_{\theta_{2},\Right}(t,x)-u_{\theta_{1},\Right}(t,x)=\sum_{\substack{\theta\in[\theta_{1},\theta_{2})\\
(\theta,t,x)\in\GS
}
}\left(u_{\theta,\Left,\Right}(t,x)-u_{\theta,\Right,\Right}(t,x)\right).\label{eq:applyuthetadiff}
\end{equation}
Now if $(\theta,t,x)\in\GS$, then $x\in\{s_{\Square}(\theta,t)\st\Square\in\LR\}$.
But if $x=\sfs_{\Left}(\theta,t)\ne\sfs_{\Right}(\theta,t)$, then
\prettyref{prop:converse} tells us that $u_{\theta,\Left,\Right}(t,x)=u_{\theta,\Right,\Right}(t,x)$,
so the contribution to the right side of \prettyref{eq:applyuthetadiff}
is zero. Thus, we obtain
\[
u_{\theta_{2},\Right}(t,x)-u_{\theta_{1},\Right}(t,x)=\sum_{\substack{\theta\in[\theta_{1},\theta_{2})\\
\sfs_{\Right}(\theta,t)=x
}
}\left(u_{\theta,\Left,\Right}(t,x)-u_{\theta,\Right,\Right}(t,x)\right)=\sum_{\substack{\theta\in[\theta_{1},\theta_{2})\\
\sfs_{\Right}(\theta,t)=x
}
}\frac{1}{\partial_{\theta}\sfs_{\Right}(\theta+,t)},
\]
with the last identity by \prettyref{prop:relatetooneoverderiv},
and thus we obtain \prettyref{eq:utheta2theta1}.
\end{proof}

\section{Verifying the hypotheses for compound Poisson forcing}

\label{sec:poisson}In this section we show that compound Poisson
processes are in $\widetilde{\Omega}$ (defined in \prettyref{def:Omegatilde2})
with probability $1$, proving \prettyref{thm:poisson-ok}. In this
section, we let $\mathbb{P}$ be a probability measure under which
$\mu$ is the measure associated to a nonnegative homogeneous compound
Poisson process on $\mathbb{R}\times\mathbb{T}$. \prettyref{thm:poisson-ok}
is a combination of the three propositions in this section. 

First we address $\widetilde{\Omega}_{1}$, defined in \prettyref{def:smallnoisezones}.
\begin{prop}
\label{prop:POmega11}We have $\mathbb{P}(\widetilde{\Omega}_{1})=1$.
\end{prop}

\begin{proof}
Let $\mathbb{U}$ be the weight distribution of $\mathbb{P}$, so
$\mathbb{U}$ is a probability measure on $(0,\infty)$. This means
that
\begin{equation}
\mathbb{U}(B)=\mathbb{P}\left(\mu(A)\in B\mid\#(\mathsf{N}\cap A)=1\right)\qquad\text{for any Borel }A\subseteq\RR\times\TT\text{ and }B\subseteq(0,\infty).\label{eq:Udef}
\end{equation}
Choose $M>0$ large enough that $\mathbb{U}\left(\left(\frac{1}{4M},\infty\right)\right)>0$.
Defining $A_{r,N}\coloneqq[r-N,r+N]\times\TT$, we see that
\begin{align*}
\mathbb{P} & \left(\#(\ssN\cap A_{r,2M})=\#(\ssN\cap A_{r,M})=1\text{ and }\mu(\ssN\cap A_{r,M})>\frac{1}{4M}\right)\\
 & =\mathbb{P}\left(\#(\ssN\cap(A_{r,2M}\setminus A_{r,M}))=0\right)\cdot\mathbb{P}\left(\#(\ssN\cap A_{r,M})=1\right)\\
 & \qquad\cdot\mathbb{P}\left(\mu(\ssN\cap A_{r,M})>\frac{1}{4M}\ \middle|\ \#(\ssN\cap A_{r,M})=1\right)\\
 & =\mathbb{P}\left(\#(\ssN\cap(A_{r,2M}\setminus A_{r,M}))=0\right)\cdot\mathbb{P}\left(\#(\ssN\cap A_{r,M})=1\right)\cdot\mathbb{U}\left(\left(\frac{1}{4M},\infty\right)\right)>0.
\end{align*}
By the independence and spatial homogeneity of the Poisson process,
this means that there are constants $\rho<1$ and $C<\infty$, independent
of $t$, such that for any $k\ge3$, we have
\begin{align*}
\mathbb{P} & \left(\exists r\in[t-kM,t-2M]\text{ s.t. }\#(\ssN\cap A_{r,2M})=\#(\ssN\cap A_{r,M})=1\text{ and }\mu(\ssN\cap A_{r,M})>\frac{1}{4M}\right)\\
 & \ge1-C\rho^{k}.
\end{align*}
The fact that $\mathbb{P}(\widetilde{\Omega}_{1})=1$ then follows
from the Borel--Cantelli theorem.
\end{proof}
Now we address $\widetilde{\Omega}_{2}$, defined in \prettyref{def:Omegatilde2}.
We note that $\widetilde{\Omega}_{2}=\widetilde{\Omega}_{2;1}\cap\widetilde{\Omega}_{2;2}$,
with $\widetilde{\Omega}_{2;1}$ and $\widetilde{\Omega}_{2;2}$ defined
in \prettyref{prop:omega21} and \prettyref{prop:omega22} below,
respectively.
\begin{prop}
\label{prop:omega21}Let $\widetilde{\Omega}_{2;1}$ be the set of
all $\mu\in\widetilde{\Omega}$ such that $\#(\GS_{\theta}\cap\ssN)\le1$
for all $\theta\in\RR$. We have $\mathbb{P}(\widetilde{\Omega}_{2;1})=1$.
\end{prop}

\begin{proof}
If $\mu\in\widetilde{\Omega}_{1}$ and $t\in\RR$, then we define
\[
\mathcal{R}_{t}\coloneqq\{\theta\in\RR\st\GS_{\theta}\cap\ssN\cap((-\infty,t]\times\TT)\ne\varnothing\}=\bigcup_{(s,x)\in\ssN\cap((-\infty,t]\times\TT)}\GS'_{s,x},
\]
with $\GS'_{s,x}$ defined in \prettyref{eq:GS'def}. By \prettyref{prop:GS'discrete}
and the countability of $\ssN$, we see that $\mathcal{R}_{t}$ is
countable for each $t$, being the countable union of countable sets.
Also, for $t\le t'$, we define
\[
\mathcal{E}_{t,t'}\coloneqq\left\{ \theta\in\RR\st\exists(s,y)\in\ssS_{\theta}\text{ s.t. }[(t,t')\times\TT]\cap\ssN=\{(s,y)\}\right\} .
\]
We note that, if $\mu\not\in\widetilde{\Omega}_{2;1}$, then there
must be some $t,t'\in\mathbb{Q}$ such that $\mathcal{R}_{t}\cap\mathcal{E}_{t,t'}\ne\varnothing$,
which means that
\begin{equation}
\mathbb{P}(\widetilde{\Omega}_{1}\setminus\widetilde{\Omega}_{2;1})\le\sum_{t,t'\in\mathbb{Q}}\mathbb{P}\left(\mu\in\widetilde{\Omega}_{1}\text{ and }\mathcal{R}_{t}\cap\mathcal{E}_{t,t'}\ne\varnothing\right).\label{eq:ctblunion}
\end{equation}

Let $\mathcal{F}_{t,t'}$ be the $\sigma$-algebra generated by $\mu|_{\RR\setminus(t,t')}$.
From \prettyref{def:minimizers}, it is not difficult to check that
$\widetilde{\Omega}_{1}$ is measurable with respect to $\mathcal{F}_{t,t'}$
for any $-\infty<t\le t'<+\infty$. (That is, it is a tail event in
an appropriate sense, although we will not need any zero-one theorems
here.) It is also not difficult to see from \prettyref{prop:basicproperties}(\ref{enu:straightlines})
that, for any fixed $\theta\in\RR$ and $t\le t'$, we have, $\mathbf{1}_{\widetilde{\Omega}_{1}(\mu)}\mathbb{P}(\theta\in\mathcal{E}_{t,t'}\mid\mathcal{F}_{t,t'})=0$.
Thus we can compute 
\begin{align*}
\mathbb{P}\left(\mu\in\widetilde{\Omega}_{1}\text{ and }\mathcal{R}_{t}\cap\mathcal{E}_{t,t'}\ne\varnothing\right) & \le\mathbb{E}[\#(\mathcal{R}_{t}\cap\mathcal{E}_{t,t'});\mu\in\widetilde{\Omega}_{1}]=\mathbb{E}\left[\mathbb{E}[\#(\mathcal{R}_{t}\cap\mathcal{E}_{t,t'})\mid\mathcal{F}_{t,t'}];\mu\in\widetilde{\Omega}_{1}\right]\\
 & =\mathbb{E}\left[\sum_{\theta\in\mathcal{R}_{t}}\mathbb{P}(\theta\in\mathcal{E}_{t,t'}\mid\mathcal{F}_{t,t'});\mu\in\widetilde{\Omega}_{1}\right]=0.
\end{align*}
Using this in \prettyref{eq:ctblunion}, we see that $\mathbb{P}(\widetilde{\Omega}_{1}\setminus\widetilde{\Omega}_{2;1})=0$
and hence that $\mathbb{P}(\widetilde{\Omega}_{2;1})=1$ by \prettyref{prop:POmega11}.
\end{proof}
\begin{prop}
\label{prop:omega22}Let $\widetilde{\Omega}_{2;2}$ be the set of
all $\mu\in\widetilde{\Omega}$ such that $\ssS_{\theta}\cap\ssN\setminus\GS_{\theta}=\varnothing$
for all $\theta\in\RR$. We have $\mathbb{P}(\widetilde{\Omega}_{2;2})=1$.
\end{prop}

\begin{proof}
Suppose that $\mu\in\widetilde{\Omega}_{1}$. If $(t,x)\in\ssS_{\theta}\setminus\GS_{\theta}$,
then we must have $X_{\theta,t,x,\Left}'(t-)\ne X_{\theta,t,x,\Right}'(t-)$
but
\begin{equation}
\int_{T_{\vee}(\theta,t,x)}^{t}X_{\theta,t,x,\Left}'(s)\,\dif s=\int_{T_{\vee}(\theta,t,x)}^{t}X_{\theta,t,x,\Right}'(s)\,\dif s,\label{eq:notglobal}
\end{equation}
as well as
\begin{align*}
0 & =\mathcal{A}_{\theta,T_{\vee}(\theta,t,x),t}[X_{\theta,t,x,\Right}]-\mathcal{A}_{\theta,T_{\vee}(\theta,t,x),t}[X_{\theta,t,x,\Left}]\\
\oversetx{\prettyref{eq:Adef}} & =\frac{1}{2}\int_{T_{\vee}(\theta,t,x)}^{t}\left(X_{\theta,t,x,\Right}'(s)^{2}-X_{\theta,t,x,\Left}'(s)^{2}\right)\,\dif s-\theta\int_{T_{\vee}(\theta,t,x)}^{t}\left(X_{\theta,t,x,\Right}'(s)-X_{\theta,t,x,\Left}'(s)\right)\,\dif s\\
 & \qquad-\mu\left(\{(s,X_{\theta,t,x,\Right}(s))\st s\in(T_{\vee}(\theta,t,x),t)\}\right)+\mu\left(\{(s,X_{\theta,t,x,\Left}(s))\st s\in(T_{\vee}(\theta,t,x),t)\}\right)\\
\oversetx{\prettyref{eq:notglobal}} & =\left[\frac{1}{2}\int_{T_{\vee}(\theta,t,x)}^{t}X_{\theta,t,x,\Square}'(s)^{2}\,\dif s-\mu\left(\{(s,X_{\theta,t,x,\Right}(s))\st s\in(T_{\vee}(\theta,t,x),t)\}\right)\right]_{\Square=\Left}^{\Square=\Right}.
\end{align*}
Fix some arbitrary $y\in\TT$ such that $(t_{1},y)\not\in\ssN$. If
$\mu\not\in\widetilde{\Omega}_{2;2}$, then there must exist $t_{1},t_{2}\in\mathbb{Q}$
with $t_{1}<t_{2}$, paths $Y_{1},Y_{2}$ connecting $(y,t_{1})$
and elements of $\mathsf{N}$ by straight line segments (of which
there are at most countably many), and $\tau,\eta\in\RR$ such that
\begin{equation}
\frac{1}{2}\int_{T_{*}(t)}^{t}(\mathcal{T}_{\tau,\eta}^{1}Y_{1})'(s)^{2}\,\dif s-\mu(\{(s,\mathcal{T}_{\tau,\eta}^{1}Y_{\Square}(s))\st s\in(T_{*}(t),t)\})\quad\text{does not depend on }i\in\{1,2\}.\label{eq:Ydoesnotdepend}
\end{equation}
Now given such $t_{1},t_{2},Y_{1},Y_{2}$, the set $\mathcal{J}(t_{1},t_{2},Y_{1},Y_{2})$
of $(t,x)\in[t_{1},t_{2}]\times\TT$ such that \prettyref{eq:Ydoesnotdepend}
can hold is at most countable. This means that 
\[
\mathbb{P}\left(\widetilde{\Omega}_{1}\setminus\widetilde{\Omega}_{2;2}\right)\le\sum_{\substack{t_{1},t_{2}\in\mathbb{Q}\\
Y_{1},Y_{2}
}
}\sum_{(t,x)\in\mathcal{J}(t_{1},t_{2},Y_{1},Y_{2})}\mathbb{P}((t,x)\in\mathsf{N})=0,
\]
and hence that $\mathbb{P}(\widetilde{\Omega}_{2;2})=1$ by \prettyref{prop:POmega11}.
\end{proof}
\begin{flushleft}
\printbibliography[heading=bibintoc]
\par\end{flushleft}
\end{document}